\documentclass[final,leqno,onefignum]{SIIMS}

\usepackage[]{amsfonts,amsmath, bm}
\usepackage{epsfig}
\usepackage[]{graphicx,subfigure}
\usepackage{algorithm,algorithmic}
\graphicspath{{./}{./figures/}}

\def \Cm {\mathbb{C}}

\def \Rm {\mathbb{R}}
\def \Sm {\mathbb{S}}

\def \Zm {\mathbb{Z}}

\def\C{\mathcal{C}}
\def\D{\mathcal{D}}

\def\K{\mathcal{K}}
\def\L{\mathcal{L}}

\def\R{\mathcal{R}}

\def\diam{ {{\rm diam}}}

\newcommand{\where}{\quad\text{ where }}
\newcommand{\qandq}{\quad\text{ and }\quad}
\newcommand{\bfe}{ {\bf e}}

\newcommand{\cout}[1]{}

\newcommand{\x}{{\mathrm x}}

\newcommand{\tu}{{\tilde{u}}}
\newcommand{\dprod}[2]{\langle #1, #2 \rangle}

\title{Numerical Implementation of Geodesic X-Ray Transforms and Their Inversion\thanks{Received by the editors XXXX; accepted for publication (in revised form) XXXX;
published electronically DATE.
\URL{siims/x-x/xxxxx.html}}}

\author{Fran\c cois Monard\thanks{Department of Mathematics, University of Washington (\email{fmonard@math.washington.edu}). This research is partially supported by NSF grant DMS-1025372.}}

\begin{document}
\maketitle
\newcommand{\slugmaster}{%
\slugger{siims}{xxxx}{xx}{x}{x--x}}
\renewcommand{\thefootnote}{\fnsymbol{footnote}}

\begin{abstract}
    We present a numerical implementation of the geodesic ray transform and its inversion over functions and solenoidal vector fields on two-dimensional Riemannian manifolds. For each problem, inversion formulas previously derived by Pestov and Uhlmann in [Int. Math Research Notices 80 (2004)] then extended by Krishnan in [J. Inv. Ill-Posed Problems {\bf 18} (2010)] are implemented in the case of simple and some non-simple metrics. These numerical tools are also used to better understand and gain intuition about non-simple manifolds, for which injectivity and stability of the corresponding integral geometric problems are still under active study.
\end{abstract}



\pagestyle{myheadings}
\thispagestyle{plain}
\markboth{FRAN\c COIS MONARD}{NUMERICAL FORWARD AND INVERSE GEODESIC RAY TRANSFORMS}

\section{Introduction}

The present article discusses a numerical implementation in {\tt MatLab} of geodesic X-Ray transforms of functions and solenoidal (i.e. divergence-free) vector fields and their inversion on two dimensional Riemannian manifolds with boundary. 

Geodesic X-ray transforms appear in problems of mathematical physics where particles travel along some curves and ``gather information'' along their path. Probably the best known example of geodesic ray transform is that of the Radon transform in two dimensions (also known as the X-Ray transform, as these two transforms coincide in two dimensions), which is the collection of integrals of a given function over all straight lines in the plane, corresponding to the case of a Euclidean metric. Reconstructing a function from its integrals along lines was first considered and solved in \cite{Radon1917} and is now used every day in medical imaging. A thorough account of theoretical and numerical aspects for this transform may be found in \cite{Natterer2001}. In the Euclidean case, solenoidal tensors of any order were also explicitely reconstructed in \cite{Sharafudtinov1994}.

When optical rays propagate in a medium with variable index of refraction, their trajectories, no longer straight lines, come as geodesics of some Riemannian metric. In this framework, the same questions (injectivity, stability, range characterization, reconstruction algorithms, inverse problems with partial data) as for the straight line case are still under active theoretical study. To the author's knowledge, numerical simulations for these transforms remain to be documented. 

When the metric is simple, injectivity over functions was proved in \cite{Mukhometov1977} and injectivity over solenoidal vector fields was established in \cite{Anikonov1978,Anikonov1997}. Under the same simplicity assumption, the problem was recently proved in \cite{Paternain2011a} to be injective over solenoidal tensors (``s-injective'') of any order, and previously in \cite{Dairbekov2006} under assumptions on the curvature. Independently, stability estimates were given in \cite{Stefanov2004} via a microlocal study of the normal operator. 

While s-injectivity is now proved, explicit methods for reconstructing the solenoidal part of tensors of order $\ge 2$ remain to be found. The case of functions and solenoidal vector fields was however tackled by Pestov and Uhlmann in \cite{Pestov2004}, where explicit Fredholm reconstruction formulas were derived for simple metrics. These formulas are exact in the case of constant-curvature metrics, and the Fredholm error was further proved in \cite{Krishnan2010} to be an $L^2$-contraction for metrics with curvature close to constant, leading again to exact reconstruction formulas in the form of Neumann series. 

When the metric is no longer simple, results are known for some geometries with certain symmetries and for smooth metrics. Sharafutdinov established in \cite{Sharafudtinov1997} s-injectivity over tensors of any order on spherically symmetric layers satisfying the Herglotz non-trapping condition. In dimension three or higher, Stefanov and Uhlmann proved in \cite{Stefanov2008} s-injectivity for real-analytic metrics satisfying some additional conditions, and Uhlmann and Vasy proved in \cite{Uhlmanna} local injectivity of the ray transform on manifolds satisfying a foliation condition, including a reconstruction algorithm. While the question of s-injectivity remains open for general domains and metrics, it is shown in \cite{Stefanov2012a} using microlocal analysis that when the metric has caustics of fold type and the manifold is two-dimensional, the singularities of the unknown function that are conormal to a fold can no longer be resolved by the ray transform, thus showing that caustic sets have detrimental effects on the stability of the ray transform. 

On the numerical side, probably one of the most thorough accounts on the two-dimensional Radon transform is found in \cite{Natterer2001}. There, of crucial help is the presence of the {\em parallel geometry} (i.e. the global parameterization of lines in terms of their distance from the center and their direction) and the Fourier Slice Theorem, providing a proper way of constructing regularized reconstruction algorithms as well as efficient FFT-based reconstruction algorithms. On the other hand, {\em fan-beam} data (i.e. direct parameterization from the influx boundary), are processed either via a reparameterization of a formula initially obtained in the parallel geometry, or ``re-binned'' into the parallel geometry in order to use its wealthier machinery. For general metrics (e.g. non-constant curvature, non-spherically symmetric), there seems to be no other obvious parameterization of the data space than the so-called {\em fan-beam} variables. This is because unless the manifold is a well-known one, one does not really knows a geodesic better than locally and cannot parameterize them globally other than from the boundary. 

The present implementation is therefore based on the fan-beam geometry, which despite its disadvantages mentioned above allows to treat the most general case while easily revisiting the well-known ones by means of the Pestov-Uhlmann reconstruction formulas \cite{Pestov2004}. The code presented comes as a handy tool for obtaining a better understanding of two-dimensional geodesic X-Ray transforms and will be used in a forthcoming theoretical and numerical study of the nonsimple case.  

\paragraph{Outline} The rest of the paper is organised as follows. Section \ref{sec:theory} covers the formulation of the problem and notation, some theoretical results of interest established in prior literature as well as the reconstruction formulas recovering functions and solenoidal vector fields from their ray transforms. Section \ref{sec:numerics} covers the numerical implementation, describing the building blocks in \S \ref{sec:buildingblocks}, treating the constant curvature cases in \S \ref{sec:constcurv} and more general cases in \S \ref{sec:nonconstcurv}. Section \ref{sec:conclu} concludes.

\section{Theoretical background} \label{sec:theory}

Let $(M,g)$ be a compact oriented, simply connected, two-dimensional Riemannian manifold with boundary. Here and below, we denote by $SM$ the unit tangent bundle $SM = \{ (\x,v)\in TM; |v| := g(v,v)^{\frac{1}{2}} = 1 \}$. Following the notational conventions in \cite{Pestov2004}, let $\nu_\x$ denote the unit inner normal to $M$ at a point $x\in \partial M$, and define
\begin{align*}
    \partial_{\pm} SM := \{ (\x,v)\in SM;\ \x\in \partial M,\ \pm \dprod{v}{\nu_\x} >0 \}.
\end{align*}
The metric $g$ induces a geodesic flow $\phi_t = (\gamma_{\x,v}(t), \dot\gamma_{\x,v}(t))$ acting on $SM$, locally described by the differential system 
\begin{align}
    \ddot \gamma_{\x,v}^k + \Gamma_{ij}^k(\gamma_{\x,v}(t)) \dot\gamma_{\x,v}^i \dot\gamma_{\x,v}^j = 0, \quad k=1,2, \quad \gamma_{\x,v}(0) = \x, \quad \dot\gamma_{\x,v}(0) = v, 
    \label{eq:geodesic}
\end{align}
defined on the domain
\begin{align}
    \D := \{ (\x,v,t); (\x,v)\in SM, \quad -\tau (\x,-v)< t <\tau(\x,v) \},
    \label{eq:Dset}
\end{align}
where $\tau(\x,v)$ is the first exit time of the geodesic $\gamma_{\x,v}$. In \eqref{eq:geodesic}, the coefficients $\Gamma_{ij}^k$ denote the Christoffel symbols
\begin{align*}
    \Gamma_{ij}^k := \frac{1}{2} g^{kl} \left( \partial_i g_{jl} + \partial_j g_{il} - \partial_l g_{ij} \right). 
\end{align*}
Such a flow can be described by means of a vector field $X$ defined on $SM$, and whose integral curves are precisely the unit-speed geodesics. 
Note that the geodesics going from $\partial M$ into $M$ can be parameterized over $\partial_+ SM$. Given a symmetric covariant $m$-tensor $f$, we define the geodesic X-ray transform of $f$, as follows 
\begin{align}
    If(\x,v) = \int_0^{\tau(\x,v)} f(\gamma_{\x,v}(t), \dot\gamma_{\x,v}(t)^m)\ dt, \quad (\x,v)\in \partial_+ SM.
    \label{eq:geoRT}
\end{align}

The tensor tomography problem consists in reconstructing $f$ (or, rather, its solenoidal part in the sense of Sharafutdinov's decomposition, see \cite[Sec. 3.3]{Sharafudtinov1994}) from $If$. 

\paragraph{Restriction to isotropic metrics}

For simplicity of implementation, we consider the case where the metric is isotropic, that is, $g_{ij} = g \delta_{ij}, \quad i,j=1,2$ for some function $g:M\mapsto [g_0,\infty)$ with $g_0 >0$. In two dimensions, this is not a restrictive assumption because {\em isothermal coordinates} (i.e. coordinates in which the metric tensor becomes isotropic) always exist \cite{Spivak1999}, globally for simply connected surfaces. It is also convenient to define the function $\lambda:= \frac{1}{2}\log g$, i.e. $g = e^{2\lambda}$. In this case, we have $g^{ij} = g^{-1} \delta_{ij}$ and the Christoffel symbols in the isothermal frame $(\partial_1, \partial_2)$ take the following expression
\begin{align}
    \begin{split}
	\Gamma_{11}^1 &= - \Gamma_{22}^1 = \Gamma_{12}^2 = \Gamma_{21}^2 = \frac{1}{2} \partial_1\log g = \partial_1 \lambda, \\
        \Gamma_{22}^2 &= - \Gamma_{11}^2 = \Gamma_{12}^1 = \Gamma_{21}^1 = \frac{1}{2} \partial_2\log g = \partial_2\lambda.
    \end{split}  
    \label{eq:christoffel}
\end{align}
It is easy to establish that geodesics have constant speed modulus, see e.g. \cite[Lemma 5.5 p70]{Lee1997}, and that the ray transform is homogeneous with respect to the geodesic's speed modulus. This is what allows us to restrict the study of this problem to the unit circle bundle $SM$, over which the velocity vector $\dot\gamma$ can be described by an angle function $\theta$ such that 
\begin{align*}
    \dot\gamma_{\x,v} (t) = e^{-\lambda (\gamma_{\x,v}(t))} \bm\theta(t), \quad \bm\theta(t) := \binom{\cos\theta(t)}{\sin\theta(t)}. 
\end{align*} 
In this setting, a geodesic is thus really described by the three scalar coordinates $(x(t),y(t),\theta(t))$, satisfying the ordinary differential system
\begin{align}
    \begin{split}
	\dot x(t) &= e^{-\lambda(x(t),y(t))} \cos\theta(t), \\
        \dot y(t) &= e^{-\lambda(x(t),y(t))} \sin\theta(t),  \\
	\dot\theta(t) &= e^{-\lambda(x(t),y(t))} ( -\sin\theta \partial_1\lambda + \cos\theta\partial_2\lambda).
    \end{split}    
    \label{eq:geodesicSM}
\end{align}
The {\em Gaussian curvature} in this representation is given by 
\begin{align*}
    \kappa = -\frac{1}{2g} (\partial_1^2 + \partial_2^2) \log g = -\Delta \lambda, \qquad x\in M,
\end{align*}
where we have defined the Laplace-Beltrami operator $\Delta := \frac{1}{g} (\partial_1^2 + \partial_2^2)$.

\paragraph{Jacobi fields, conjugate points and simple metrics} For $(\x,v)\in SM$ an initial point, let $\gamma(t)\equiv \gamma_{\x,v}(t)$ be the geodesic with initial conditions $(\x,v)$. With $\R$ denoting the curvature tensor, the following Jacobi field $J(t)$ defined on the geodesic curve above by the equation
\begin{align*}
    D_t^2 J + \R(J,\dot\gamma)\dot\gamma = 0, \quad J(0) = 0, \quad D_t J(0) = v^\perp,
\end{align*}
due to its initial conditions, is such that $J(t)\cdot\dot\gamma(t) = 0$ for all $t$, so there exists a function $b(t)$ such that for every $t$, $J(t) = b(t) \dot\gamma(t)^\perp$. Now it is easy to establish that $b(t) \equiv b_{x,v}(t)$ satisfies the differential equation
\begin{align*}
    \ddot{b} + \kappa(\gamma(t)) b(t) = 0,\quad b(0) = 0, \quad \dot{b}(0) = 1.
\end{align*}
If $0< t < \tau(\x,v)$ is such that $b(t) = 0$, then one says that the points $\x$ and $\gamma_{\x,v}(t)$ are {\em conjugate points}. The points that are conjugate to $\x$ are precisely those points where the exponential map at $\x$ fails to be a diffeomorphism. The metric is said to be {\em simple} if $\partial M$ is {\em strictly convex} in the sense that the second fundamental form is positive definite at the boundary, and if $(M,g)$ is free of conjugate points, i.e. the function $(\x,v,t)\mapsto b_{\x,v}(t)$ never vanishes outside $\{t=0\}$ on the set $\D$ defined in \eqref{eq:Dset}.

\paragraph{Terminator constants}

While very few results are known in the non-simple case, it is unclear whether simplicity alone determines the borderline of injectivity. In that regard, a finer tool is that of the {\em terminator constant} $\beta$, as described for instance in \cite{Paternain2012a}. For given $0<\beta<\infty$, the manifold $(M,g)$ is said to be {\em free of $\beta$-conjugate points} if the function $(\x,v,t)\mapsto b_{\beta}(\x,v,t)$ defined by the modified Jacobi equation
\begin{align*}
    \ddot{b}_\beta (t) + \beta \kappa(\gamma(t)) b_\beta(t) = 0,\quad b_\beta(0) = 0, \quad \dot{b}_\beta(0) = 1
\end{align*}
never vanishes outside $\{t=0\}$ on the set $\D$ defined in \eqref{eq:Dset}. Thanks to results on second-order ODE's, if $(M,g)$ is free of $\beta_0$-conjugate points, then it is also free of $\beta$-conjugate points for any $\beta\le \beta_0$. This allows to define the {\em terminator constant} $\beta_{Ter}$ of $(M,g)$ as 
\begin{align}
    \beta_{Ter} := \sup\{\beta\in [0,\infty]:\ (M,g) \text{ is free of $\beta$-conjugate points} \}.
    \label{eq:betater}
\end{align}
As a particular case, a manifold is simple if and only if $\beta_{Ter} > 1$. 

\subparagraph{Numerical test for simplicity (or absence of $\beta$-conjugate points)}

For a fixed value of $\beta$, one can test numerically whether $(M,g)$ is free of $\beta$-conjugate points by testing the non-vanishing of the function $b_\beta$ over a family of geodesics sent into the domain from a fine enough discretization of $\partial_+ SM$. This will be enough to test all $\beta$-conjugate points. Indeed, let $\gamma_{\x,v}$ a geodesic with basepoint $(\x,v)\in \partial_+ SM$ be such that $b_\beta$ does not vanish over $(0,\tau(\x,v)]$. Then by virtue of Sturm's separation theorem, no other solution of $\ddot{a} + \beta \kappa(\gamma_{\x,v}(t))a = 0$ can have two consecutive zeros over $(0,\tau(\x,v)]$. This precisely prevents the existence of pairs of $\beta$-conjugate points along the geodesic $\gamma_{\x,v}$. Doing this for every geodesic curve cast from the boundary is thus enough to ensure the absence of $\beta$-conjugate points. Simplicity is therefore tested using the particular value $\beta=1$.

\paragraph{A brief introduction to the geometry of $SM$ and the Pestov-Uhlmann inversion formulas}
We now give a brief overview of the Fredholm reconstruction formula for functions and solenoidal one-forms derived by Pestov and Uhlmann in \cite{Pestov2004}. To this end, we must introduce some additional machinery. Although the Pestov-Uhlmann formulas were not formulated in the $SM$ formalism, this latter vocabulary is slightly easier to comprehend and implement numerically, hence the author's choice to present this version, following the presentation in \cite{Paternain2011a}. 

There exists a circle action on the unit tangent bundle $SM$, whose infinitesimal generator, also called the {\em vertical vector field}, is given by $V \equiv\frac{\partial}{\partial\theta}$. From $X,V$ one may construct a global frame of $T(SM)$ by constructing the vector field $X_\perp := [X,V]$, where $[\cdot,\cdot]$ stands for the Lie bracket, or commutator, of two vector fields. One also has the additional structure equations $[V,X_\perp] = X$ and $[X,X_\perp] = \kappa V$, with $\kappa$ the Gaussian curvature. In isothermal coordinates $(x,y,\theta)$, these vector fields read
\begin{align}
    \begin{split}
	X &= e^{-\lambda} \left( \bm\theta\cdot\nabla + (\bm\theta^\perp \cdot\nabla\lambda)\ \partial_\theta \right), \\
	X_\perp &= - e^{-\lambda} \left( \bm\theta^\perp\cdot\nabla - (\bm\theta\cdot\nabla\lambda)\ \partial_\theta \right),	
    \end{split}
    \label{eq:XXperpiso}  
\end{align}
where we have defined $\nabla \equiv (\partial_1,\partial_2)$ (this notation will not conflict with the language of connections, as the latter will not be used here) as well as $\bm\theta:= \binom{\cos\theta}{\sin\theta}$. We can then define a Riemannian metric on $SM$ by declaring $(X,X_\perp,V)$ to be an orthonormal basis and the volume form of this metric will be denoted by $d\Sigma^3$ (in $(x,y,\theta)$ coordinates, this form becomes $e^{2\lambda} dx\ dy\ d\theta$). The fact that $(X,X_\perp,V)$ are orthonormal together with the structure equations implies that the Lie derivative of $d\Sigma^3$ along the three vector fields vanishes, therefore these vector fields are volume preserving. Introducing the inner product 
\begin{align*}
    (u,v) = \int_{SM} u\bar{v}\ d\Sigma^3, \qquad u,v:SM\to \Cm,    
\end{align*}
with the bar denoting conjugation, the space $L^2(SM,\Cm)$ decomposes orthogonally as a direct sum
\begin{align*}
    L^2(SM, \Cm) = \bigoplus_{k\in \Zm} H_k,
\end{align*}
where $H_k$ is the eigenspace of $-iV$ corresponding to the eigenvalue $k$. A smooth function $u:SM\to \Cm$ has a Fourier series expansion
\begin{align*}
    u = \sum_{k=-\infty}^{\infty} u_k(\x,\theta), \where \quad u_k(\x,\theta) = e^{ik\theta} \tu_k (\x), \quad \tu_k(\x) = \frac{1}{2\pi} \int_{\Sm^1} u(\x,\theta) e^{-ik\theta}\ d\theta.
\end{align*}
We also define the even/odd decomposition of such functions as
\begin{align}
    u = u_+ + u_-, \where\quad  u_+ := \sum_{k \text{ even}} u_k \qandq u_- := \sum_{k \text{ odd}} u_k.
    \label{eq:oddeven}
\end{align}
In this decomposition, a diagonal operator of particular interest is the so-called {\em Hilbert transform} $H$, whose action is best described on the Fourier components of a given function $u$ as 
\begin{align}
    (Hu)_k = -i\ \text{sgn}(k) u_k, \qquad k\in \Zm, \qquad(\text{with the convention } \text{sgn}(0) = 0).
    \label{eq:Hilbaction}
\end{align}

A crucial identity for the sequel is the following commutator formula, first derived in \cite{Pestov2005}.
\begin{lemma} \label{lem:commutator}
    The following identity holds for every $u\in \C^\infty(SM)$:
    \begin{align}
	[H,X] u = X_\perp u_0 + (X_\perp u)_0.
	\label{eq:commutator}
    \end{align}    
\end{lemma}

With the above tools at hand, let us just mention that further concepts (such as holomorphicity with respect to the fiber) can be introduced, which allowed the authors of \cite{Paternain2013a,Paternain2012,Paternain2012a,Paternain2011a,Salo2011} to prove, among other results, s-injectivity of the ray transform over tensors of any order for simple metrics, characterization of the range of the ray transform and a reconstruction procedure for the attenuated ray transform over functions.

\paragraph{The Pestov-Uhlmann reconstruction formulas}
We are now ready to present the reconstruction formulas. For $f\in L^2(SM)$, let us define $u^f(\x,\theta)$ to be the solution to the transport problem 
\begin{align}
    X u = -f, \qquad u|_{\partial_+ SM} = 0,
    \label{eq:uf}
\end{align}
so that $u|_{\partial_- SM} = If$. For $w$ defined on $\partial_+ SM$, we also denote by $w_\psi = w\circ\alpha\circ\psi$ the unique solution to the transport problem
\begin{align*}
    X u = 0, \qquad u|_{\partial_+ SM} = w,
\end{align*}
where $\alpha$ is the scattering relation and $\psi(\x,v) := (\gamma_{\x,v} (\tau(\x,v)), \dot\gamma_{\x,v} (\tau(\x,v)))\in \partial_- SM$ for any $(\x,v)\in SM$. In other words, $\alpha\circ\psi(\x,v)$ is a ``base point'' map that finds the unique point of $\partial_+ SM$ belonging to the same geodesic curve as $(\x,v)$.

When $f\in \C_0^\infty (M)$ does not depend on $\theta$, we define the operator $W f := (X_\perp u^f)_0$. It is shown in \cite{Pestov2004} that when the metric is simple, the operator $W$ can be extended as a smoothing operator $W:L^2(M)\to \C^\infty(M)$. Moreover, this operator vanishes identically if the scalar curvature is constant. The $L^2(M)$-adjoint operator $W^\star$ is given by 
\begin{align*}
  W^\star h := \left( u^{X_\perp h} \right)_0.
\end{align*}

Recall the following theorem due to Pestov and Uhlmann.
\begin{theorem}[Theorem 5.4 in \cite{Pestov2004}] \label{thm:PU}
    Consider $f\in L^2(M)$ and $h\in C_0^1(M)$ giving rise to the solenoidal vector field $X_\perp h$ and denote by $I_0 f$ and $I_1 X_\perp h$ their respective X-ray transforms. Then one has the following two formulas
    \begin{align}
	f + W^2 f &= - (X_\perp w^{(f)}_\psi)_0, \quad\text{where}\quad  w^{(f)} := (H I_0 f)_-, \label{eq:frc} \\
	h + (W^\star)^2 h &= - ( w^{(h)}_\psi )_0, \quad\text{where}\quad  w^{(h)} := (H I_1 X_\perp h)_+ . \label{eq:hrc} 
    \end{align}
\end{theorem}

The proof of Theorem \ref{thm:PU} as presented in the notation above (slightly different from that of the original paper) may be found in \cite{Monard2013a} in a slightly more general context (inversion of the ray transform over symmetric differentials). 

Although the original theorem is stated for a simple manifold, note that formulas \eqref{eq:frc}-\eqref{eq:hrc} do not require simplicity, only that the transport equation \eqref{eq:uf} be well-defined, which requires $(M,g)$ to be non-trapping. Simplicity enters the picture when proving that $W,W^\star$ are both compact, so that \eqref{eq:frc}-\eqref{eq:hrc} both satisfy Fredholm alternatives, in particular $f$ and $h$ can be reconstructed up to the finite-dimensional spaces $\ker(I + W^2)$ and $\ker (I+ (W^\star)^2)$ of smooth ghosts. 

In the case where the manifold is not simple, the operators $W,W^\star$ may no longer be smoothing operators, as predicted in \cite{Stefanov2012a} using microlocal analysis on the normal operator $I_0^\star I_0$. 

\section{Numerical Implementation} \label{sec:numerics}

All simulations and codes below are done using {\tt MatLab}.

In our simulations, the domain's parameterization is star-shaped with respect to $(0,0)$ so that the boundary of the domain may be described as $\{ x (\beta) = r(\beta)\cos\beta, y(\beta)= r(\beta)\sin\beta: 0\le \beta\le 2\pi \}$, with $r$ a smooth positive function bounded away from zero. The domain is thus described as 
\begin{align}
    M = \{(x,y)\in \Rm^2: 0\le x^2 + y^2 \le r(\arg(x,y))^2 \}.
    \label{eq:domain}
\end{align}
where $\arg$ denotes the argument function. Depending on the metric, the strict convexity may or may not be satisfied for a given domain. The grid where the function is represented is of size $n\times n$ and is an equispaced discretization of the square $[-r_{max}, r_{max}]^2$ with $r_{max} = \max(r)$. On the gridpoints that lie outside the domain \eqref{eq:domain}, the function is assumed to be zero, and the values there are never used. 

With the definition \eqref{eq:domain} of $M$, the influx boundary $\partial_+ SM$ can be thus viewed (and will be parameterized as)
\begin{align*}
    \Sm^1 \times \left(-\frac{\pi}{2},\frac{\pi}{2}\right) \ni (\beta,\alpha) \mapsto \left( x(\beta) , y(\beta), \nu(\beta) + \alpha \right) \in \partial_+ SM,
\end{align*}
where $\nu(\beta)$ denotes the angle between $\bfe_1$ and the unit inner normal at $(x(\beta),y(\beta))$. 

The parameterization above generalizes the {\em fan-beam} geometry, which has been widely studied in the case where $r(\beta) = R$ is constant (i.e. $M$ is a disk), see e.g. \cite{Natterer2001}.

\begin{figure}[htpb]
    \centering 
    \subfigure[$f$ non-smooth, circle]{
		\includegraphics[width=0.22\textwidth]{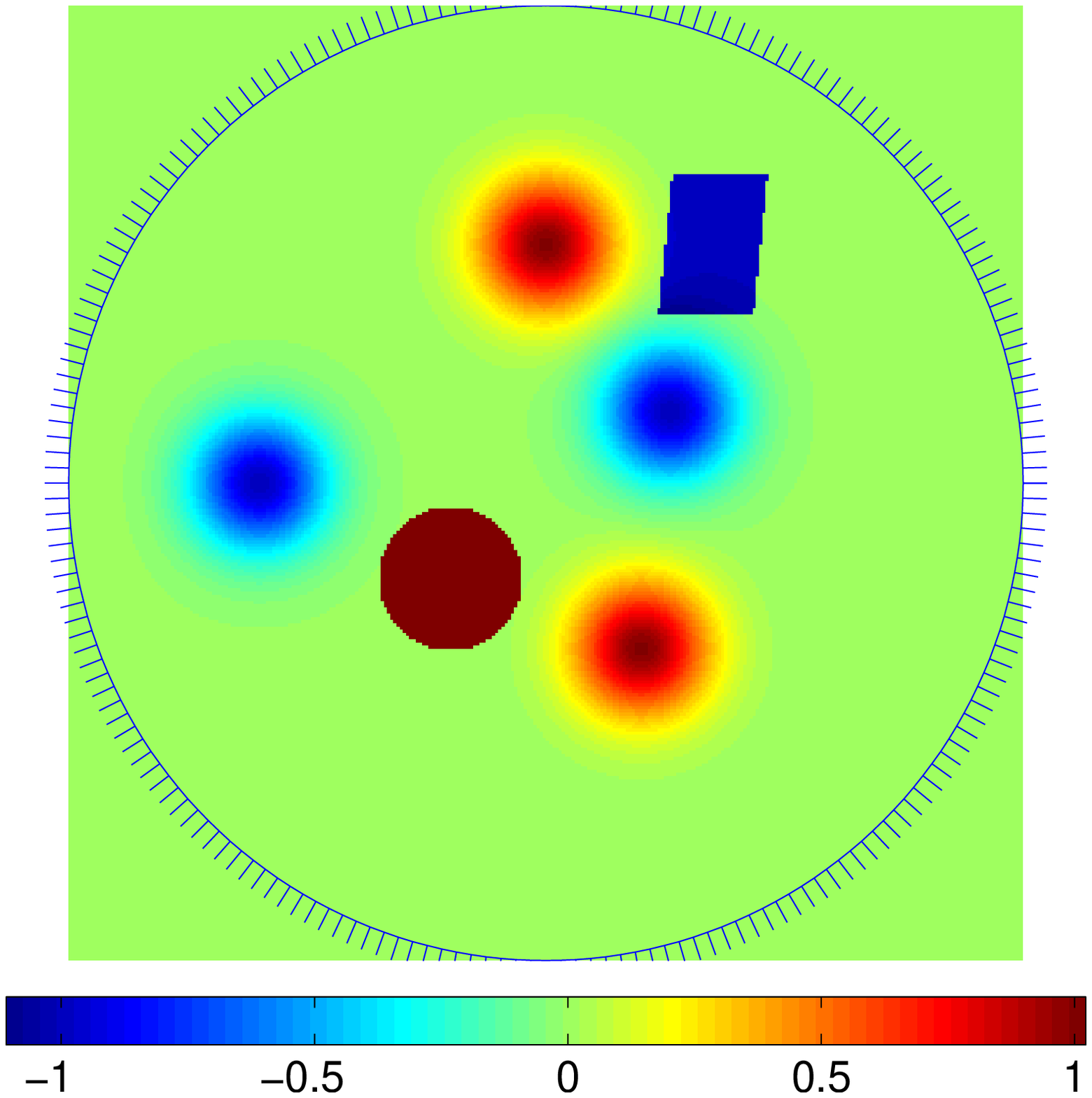}
    \label{fig:phantom1}
    }
    \subfigure[$f$ non-smooth, ellipse]{
    \includegraphics[width=0.22\textwidth]{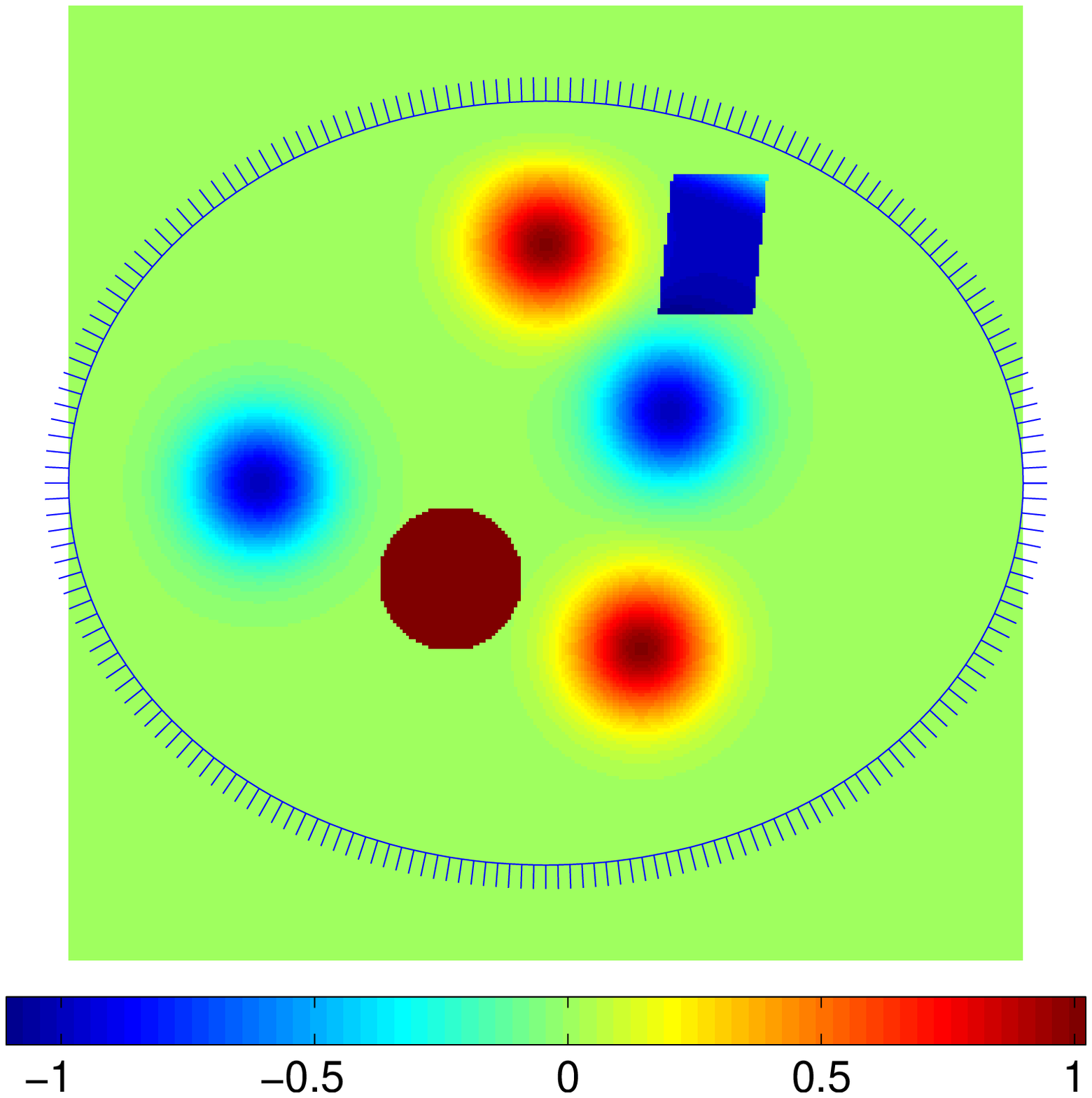}
    \label{fig:phantom2}
    }
    \subfigure[$f$ non-smooth, $\qquad$ perturbed circle]{
    \includegraphics[width=0.22\textwidth]{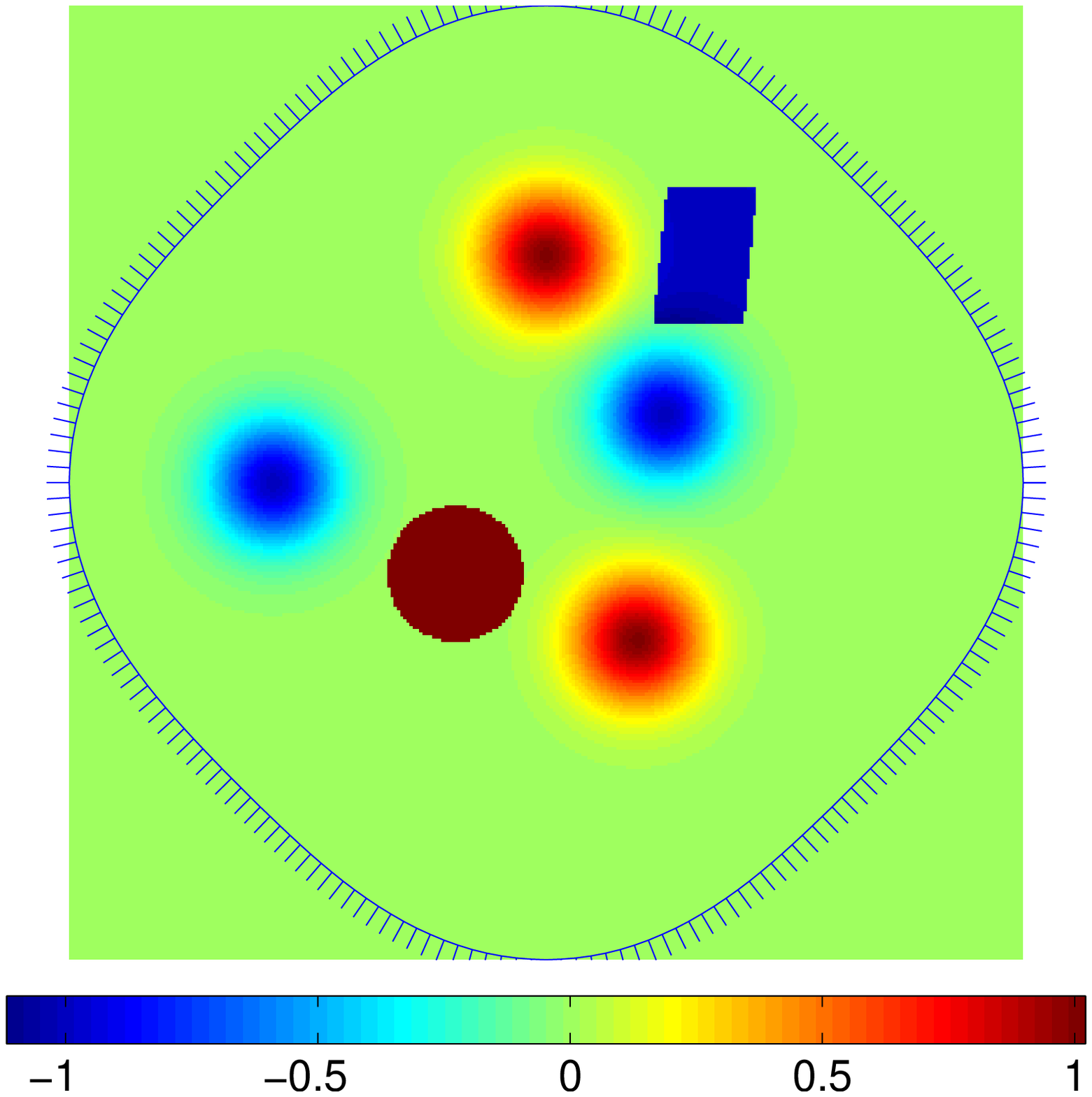}
    \label{fig:phantom3}
    }
    \subfigure[$f$ smooth, circle]{
    \includegraphics[width=0.22\textwidth]{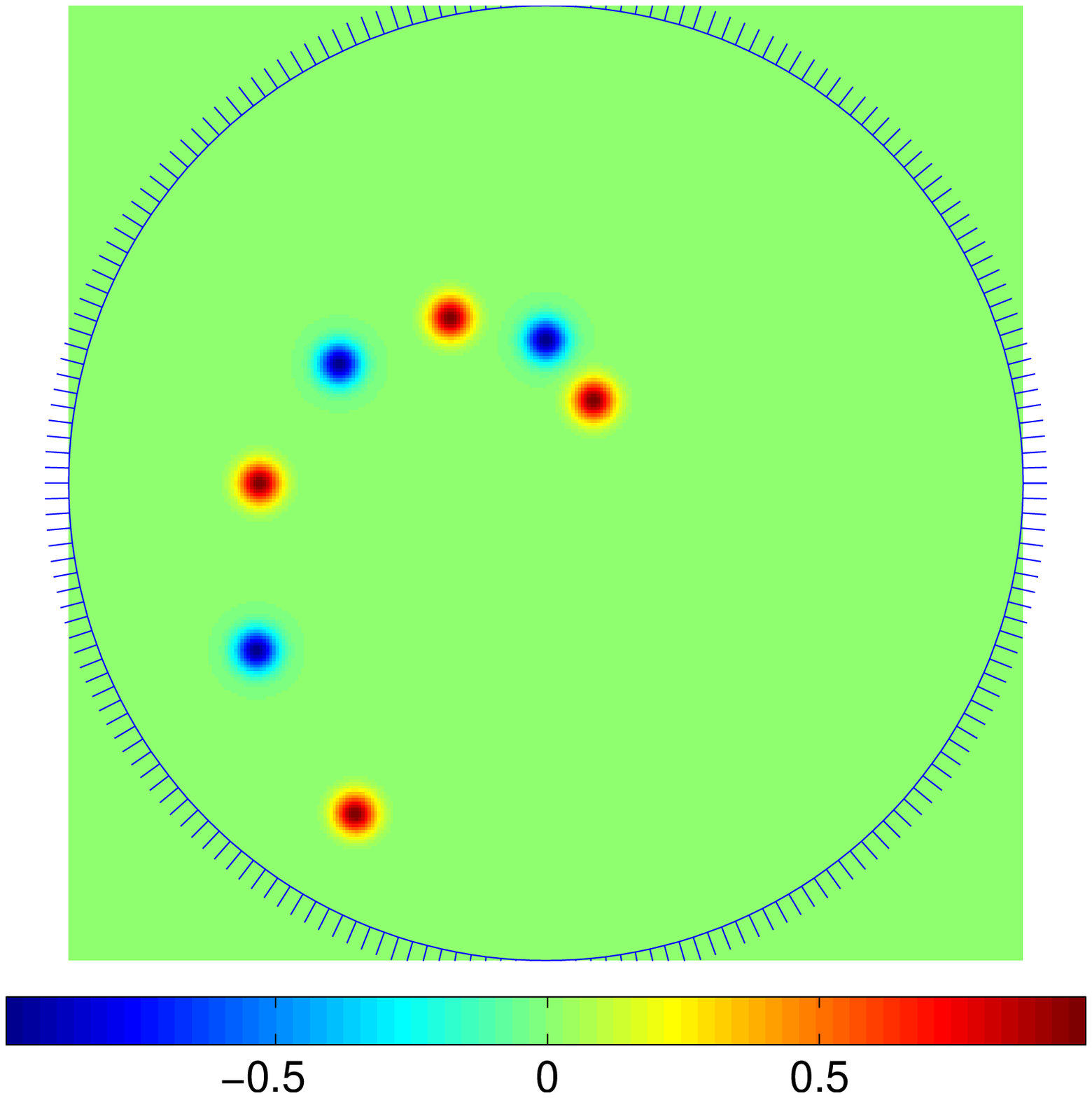}
    \label{fig:phantom4}
    }
    \caption{Phantoms and domain boundaries used throughout this paper.}
    \label{fig:unknowns}
\end{figure}

Examples of phantoms and domains used in this paper are given Fig. \ref{fig:unknowns}.



\subsection{Building blocks} \label{sec:buildingblocks}

\subsubsection{The forward operators $I_0 f$ and $I_1 [X_\perp h]$}

Computing the forward operators $I_0 f$ of a given function $f$ and $I_1 [X_\perp h]$ of a given solenoidal vector field $X_\perp h$ consists of the following steps. 
\begin{remunerate}
    \item Discretizing $\partial_+ SM$ appropriately. Here we will choose $2n\times n$ equispaced points in $[0,2\pi]\times [-\frac{\pi}{2},\frac{\pi}{2}]$, where $n$ is the sidelength of the reconstruction grid.
    \item For each boundary point $(\beta,\alpha)$ in this discretization, compute the corresponding geodesic $\{ (\gamma_{\beta,\alpha},\dot{\gamma}_{\beta,\alpha})(t),\ 0\le t\le \tau(\beta,\alpha) \}$ by solving system \eqref{eq:geodesicSM} numerically with initial conditions
	\begin{align*}
	    x(0) = r(\beta) \cos \beta, \quad y(0) = r(\beta) \sin\beta, \quad \theta(0) = \nu(\beta) + \alpha. 
	\end{align*}
	This is done by marching forward in time with stepsize $\Delta t$, until the geodesic exits the domain. The exit test is given by the condition defining $M$ in \eqref{eq:domain}. The outcome of such a procedure is a collection of points of the form $(x_{\beta,\alpha}^p, y_{\beta,\alpha}^p, \theta_{\beta,\alpha}^p)_{p=1}^N$. It is sufficient to take $N$ as an integer larger than $\diam(M)/\Delta t$. The metric and its partial derivatives are defined by analytic expressions so that there is no particular underlying Eulerian grid in the forward problem. 
    \item Using the computed geodesics, compute $I_0 f(\beta,\alpha)$ and $I_1[X_\perp h](\beta,\alpha)$ by the following quadrature rules
	\begin{align*}
	    I_0 f(\beta,\alpha) &\approx \Delta t\sum_{p=1}^N f(x_{\beta,\alpha}^p, y_{\beta,\alpha}^p), \\ 
	    I_1[X_\perp h](\beta,\alpha) &\approx \Delta t\sum_{p=1}^N e^{-\lambda(x_{\beta,\alpha}^p,y_{\beta,\alpha}^p)} \frac{h(x_{\beta,\alpha}^{p,+}, y_{\beta,\alpha}^{p,+}) - h(x_{\beta,\alpha}^{p,-}, y_{\beta,\alpha}^{p,-})}{2\Delta t},
	\end{align*}	
	where we have defined $x_{\beta,\alpha}^{p,\pm} := x_{\beta,\alpha}^p \pm \Delta t\sin \theta_{\beta,\alpha}^p$ and $y_{\beta,\alpha}^{p,\pm} := y_{\beta,\alpha}^p \mp \Delta t\cos \theta_{\beta,\alpha}^p$. When computing either of the expressions above, accessing the values of $f$ and $h$ can be handled in two different ways:
	\begin{romannum}
	    \item If the input function is described by analytic expressions (i.e. function handle), its values are computed exactly.
	    \item If the input function is given over a cartesian grid, its values are computed by bilinear interpolation of its values at the gridpoints. 
	\end{romannum}
\end{remunerate}

For a fixed boundary point characterized by $\beta$, the operations above are vectorized with respect to $\alpha$ so that the only {\tt for} loop is in $\beta$.

\subsubsection{Inversion}

\paragraph{Right-hand-side of \eqref{eq:frc}}
We first rewrite the right-hand-side of \eqref{eq:frc} in such a way that differentiation only occurs on the final cartesian grid (rewrite $w^{(f)} = w$ in the calculation below):
\begin{align*}
    -\frac{1}{2\pi} \int_{\Sm^1} X_\perp w_\psi(\x,\theta)\ d\theta &= \frac{e^{-\lambda(\x)}}{2\pi} \int_{\Sm^1} \bm\theta^\perp\cdot\nabla w_\psi - (\bm\theta\cdot\nabla\lambda) \partial_\theta w_\psi\ d\theta \\
    &= \frac{e^{-\lambda(\x)}}{2\pi} \int_{\Sm^1} \bm\theta^\perp\cdot\nabla w_\psi + (\bm\theta^\perp\cdot\nabla\lambda) w_\psi\ d\theta \\
    &= \frac{e^{-2\lambda(\x)}}{2\pi} \int_{\Sm^1} \bm\theta^\perp\cdot\nabla (e^{\lambda} w_\psi)\ d\theta \\
    &= \frac{e^{-2\lambda(\x)}}{2\pi} \nabla\cdot \left( e^{\lambda} \int_{\Sm^1} \bm\theta^\perp w_\psi(\x,\theta)\ d\theta \right).
\end{align*}
At each point of the domain, the computation of the right-hand-side of \eqref{eq:frc} consists of the following steps:
\begin{remunerate}
    \item Compute $w^{(f)} = (H I_0 f)_- = H (I_0 f)_-$. First extend the data in an odd manner in the $\alpha$ variable. Then compute the fiberwise Hilbert transform: this is done slice-by-slice via Fast Fourier Transform. Finally, restrict it back to the influx boundary.
    \item For each gridpoint, compute $$u(\x) = \int_{\Sm^1} w^{(f)}_\psi(\x,\theta) \cos\theta\ d\theta \quad \text{and} \quad v(\x) = \int_{\Sm^1} w^{(f)}_\psi(\x,\theta)\sin\theta \ d\theta, $$ where each access $w_\psi(\x,\theta)$ requires computing the basepoint $\alpha\circ\psi(\x,\theta)$ by following the geodesic with initial conditions $(\x,\theta)$ backwards. 
    \item Compute $\frac{e^{-2\lambda}}{2\pi} (- \partial_x (e^\lambda v) + \partial_y (e^\lambda u))$ at each point of the reconstruction grid using centered finite differences and pointwise multiplications. 
\end{remunerate}

\paragraph{Right-hand-side of \eqref{eq:hrc}}

Computing the right-hand-side of \eqref{eq:hrc} requires fewer steps than the previous one:
\begin{remunerate}
    \item Compute $w^{(h)} = (H I_1 X_\perp h)_- = H (I_1 X_\perp h)_-$. First extend the data in an even manner in the $\alpha$ variable. Then compute the fiberwise Hilbert transform, done as above. Restrict back to the influx boundary.
    \item For each gridpoint, compute $h(\x) = \frac{1}{2\pi} \int_{\Sm^1} w^{(h)}_\psi(\x,\theta)\ d\theta$. Again, each access $w^{(h)}_\psi(\x,\theta)$ requires computing the basepoint $\alpha\ \circ\ \psi(\x,\theta)$ by following the geodesic with initial conditions $(\x,\theta)$ backwards. 
\end{remunerate}

In both inversions above, the computational bottleneck comes from computing the basepoint of every gridpoint and for every direction. The codes above are vectorized with respect to the gridpoint $\x$ so that the only {\tt for} loop is the computation of the integrals in $\theta$. 

\begin{figure}[htpb]
    \centering 
    \subfigure[$I_0 f$]{
		\epsfig{figure=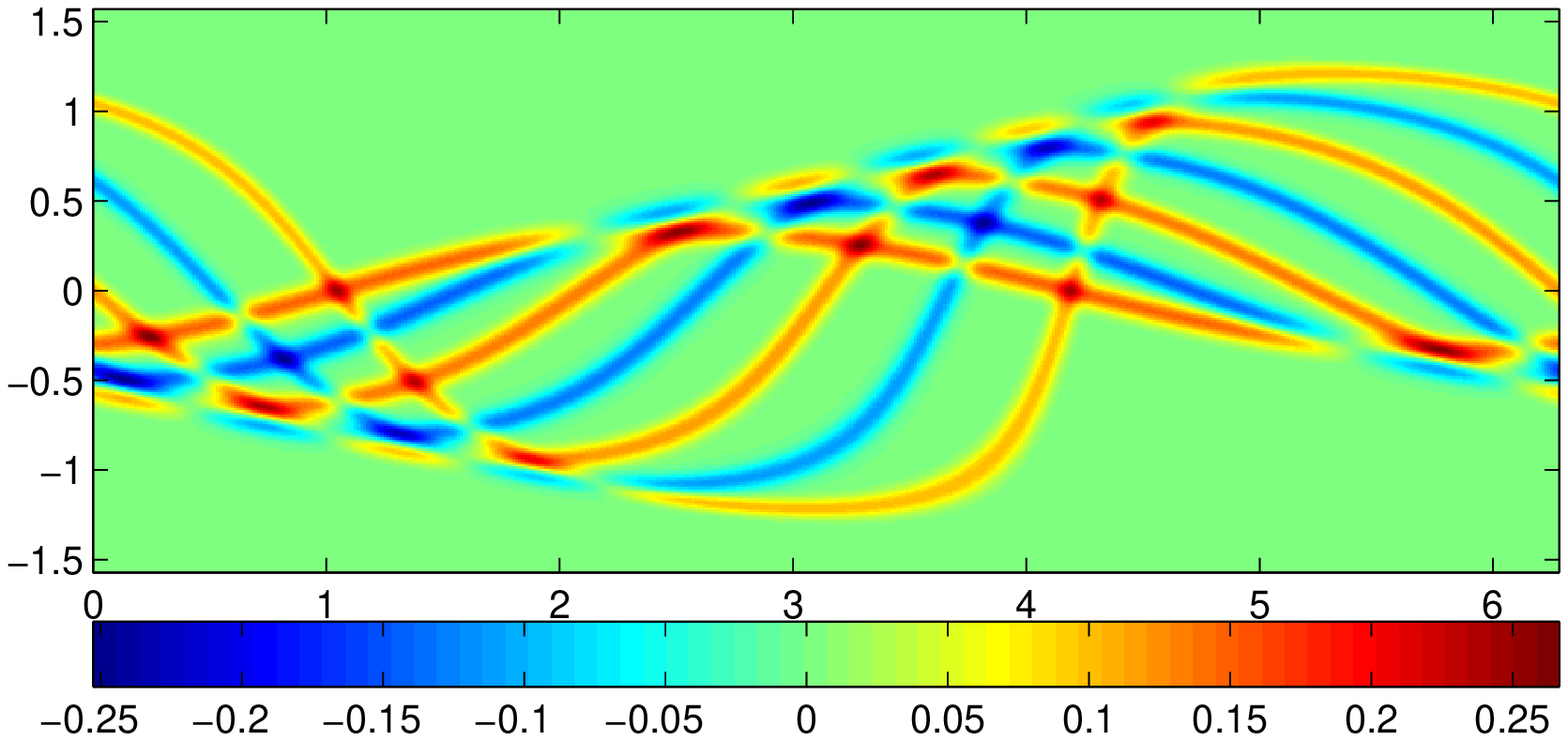,width=0.47\textwidth}    
    \label{fig:fwd1}
    }
    \subfigure[$I_1 X_\perp f$]{
    \epsfig{figure=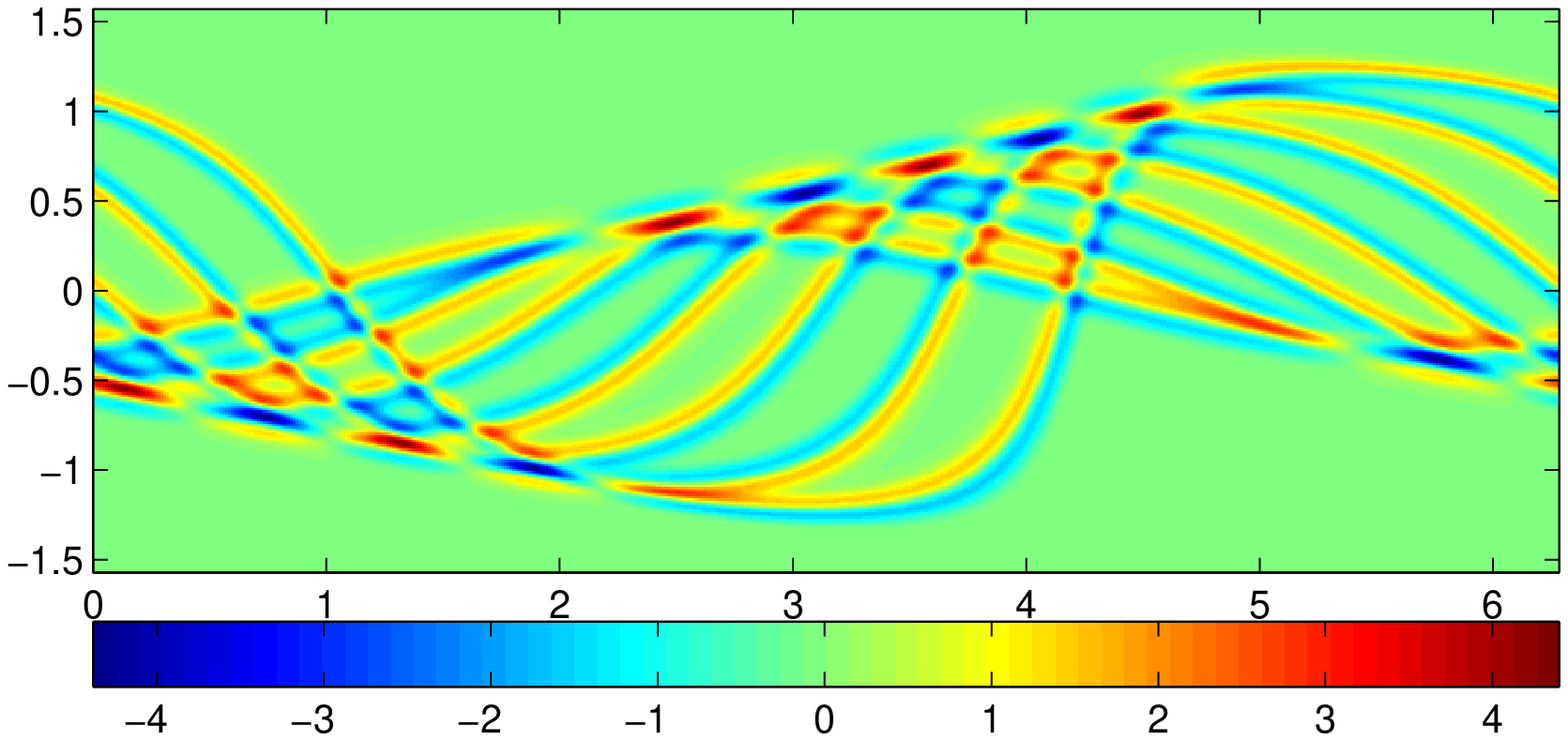,width=0.47\textwidth}    
    \label{fig:fwd2}
    }
    \subfigure[$H (I_0 f)_-$]{
    \epsfig{figure=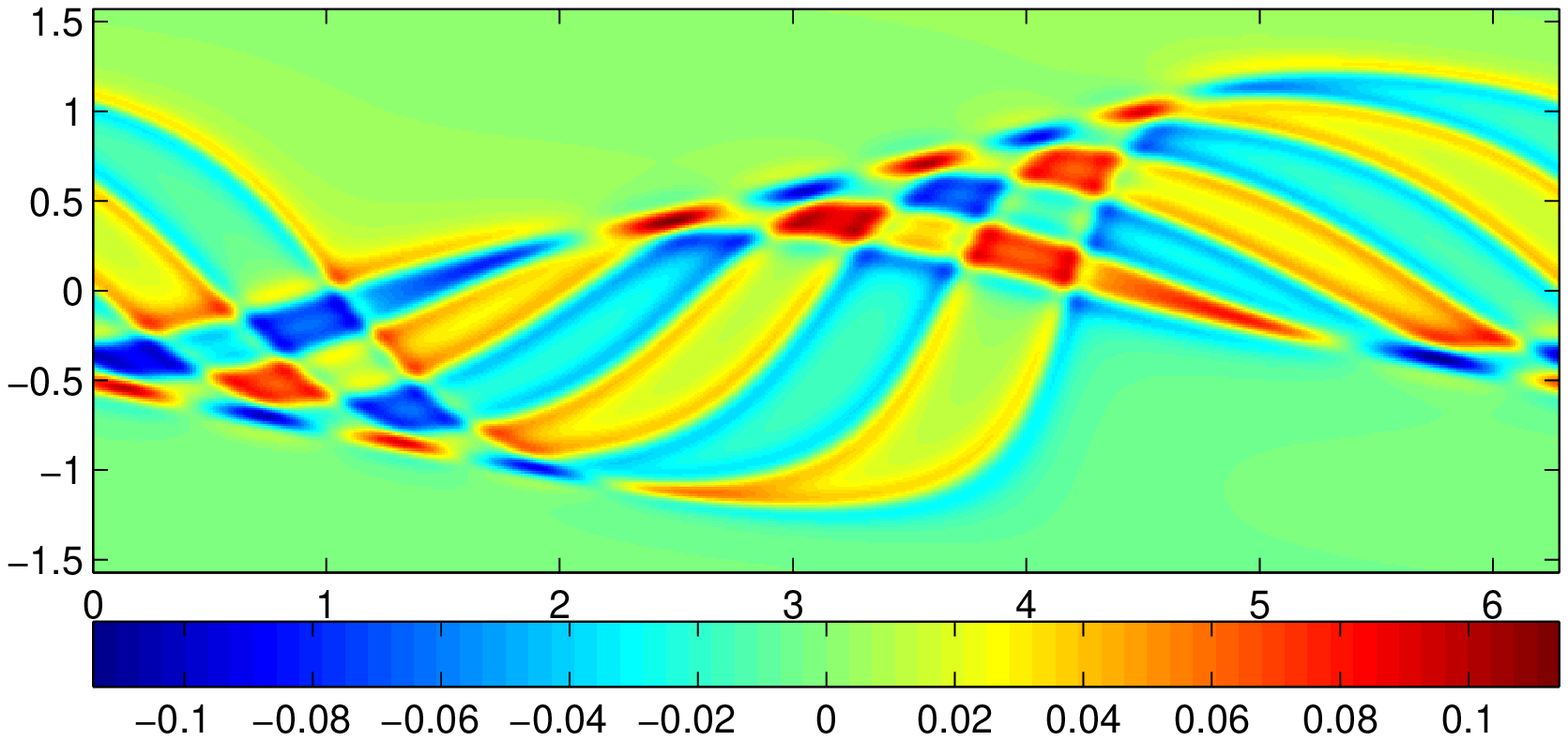,width=0.47\textwidth}    
    \label{fig:fwd3}
    }
    \subfigure[$H (I_1 X_\perp f)_+$]{
    \epsfig{figure=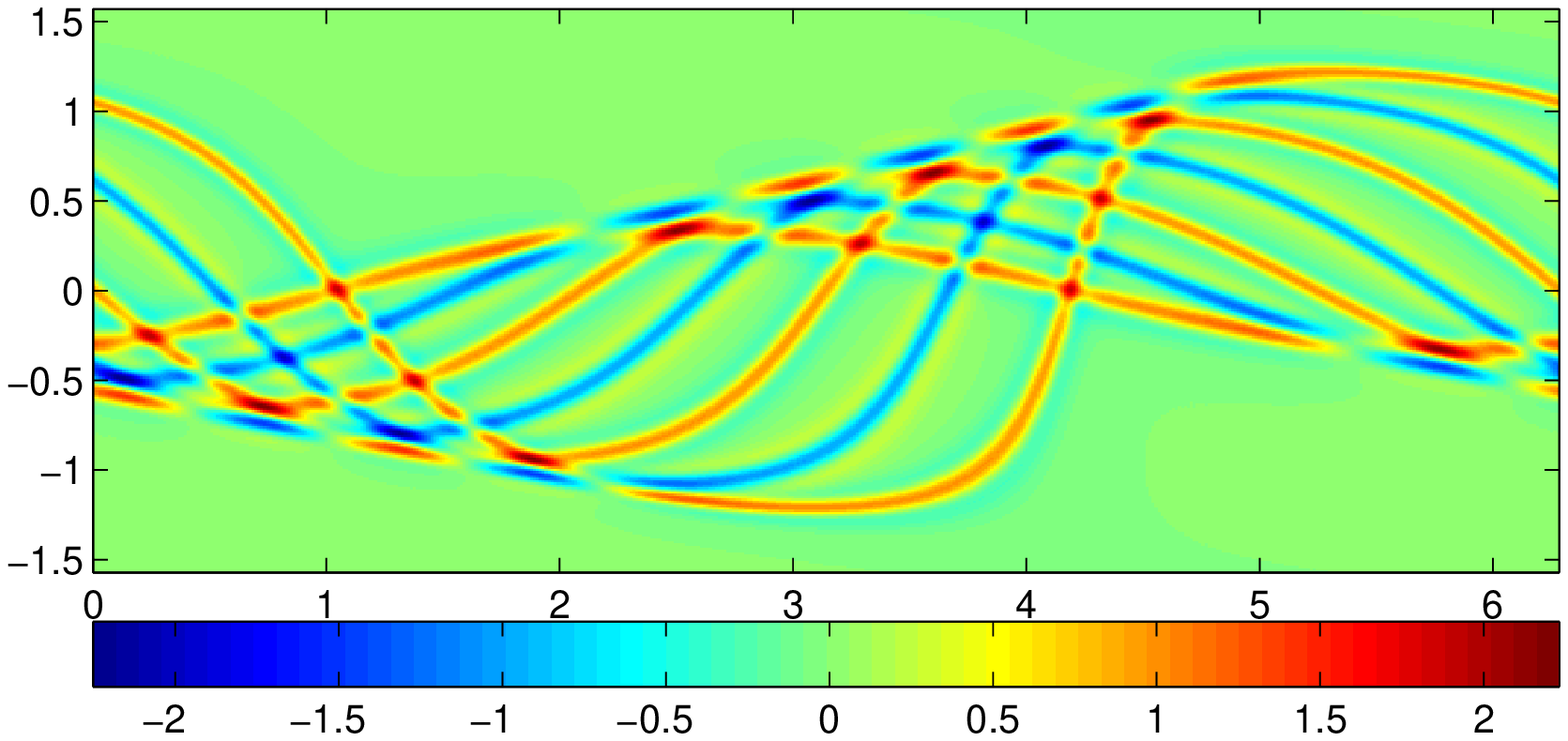,width=0.47\textwidth}    
    \label{fig:fwd4}
    }
    \caption{Examples of ray transforms and their Hilbert transforms. The function and domain's boundary are those of Fig. \ref{fig:phantom4}. The metric is given by \eqref{eq:poscurvg} with $R=1.2$. The $(x,y)$ axes are the variables $(\beta,\alpha)$ ranging in $[0,2\pi]\times \left[ -\frac{\pi}{2}, \frac{\pi}{2} \right]$.}
    \label{fig:fwd}
\end{figure}

Examples of forward transforms $I_0$ and $I_1$, as well their preprocessing before backprojection (odd or even extension, then Hilbert transform, then restriction to influx boundary), are presented Fig. \ref{fig:fwd}. 

\subsection{Constant curvature manifolds and one-shot inversions} \label{sec:constcurv}

We start with the case of manifolds with constant curvature, in which case the operators $W$ and $W^\star$ vanish identically, so that \eqref{eq:frc} and \eqref{eq:hrc} are exact reconstruction formulas. Such problems were studied and solved early on in the case of Euclidean space (Radon transform, see \cite{Radon1917}) and symmetric functions on the two-sphere (Funk transform, see \cite{Funk1916}), then generalized to the context of pairs of homogeneous spaces of the same group, see \cite{Helgason1999} for a thorough account. Such accounts do not necessarily consider manifolds-with-boundary, although questions of injectivity and reconstruction formulas from the latter to the former could be done by trivially extending the unknown function to the whole space (by zero in the Euclidean case, or into a symmetric function in the two-sphere case and considering that the simplicity condition forces the initial manifold to be stricly included in a hemisphere). 

Nonetheless, when the manifold is no longer a subset of a homogeneous space, the family of geodesics can only be parameterized from the influx boundary, generalizing the {\em fan-beam} coordinates (see e.g. \cite{Natterer2001}). This is what we implement in the present paper. The formulas implemented also allow one to reconstruct solenoidal vector fields, which are not considered in the literature mentioned above.

\paragraph{Positive curvature} 

On $\Rm^2$, and given $R>0$ fixed, the isotropic metric
\begin{align}
    g_{R,+}(x,y) := \frac{4R^4}{(x^2+y^2+R^2)^2},
    \label{eq:poscurvg}
\end{align}
has constant positive Gaussian curvature $\kappa = \frac{1}{R^2}$. A way to obtain it is by considering the centered sphere of radius $R$ with the metric induced by Euclidean $\Rm^3$, and pulling back this metric to $\Rm^2$ using the inverse of the stereographic projection map. The circle of center $(0,0)$ and radius $R$ is a notable closed geodesic, and any domain strictly enclosed in it corresponds to a subdomain of the $2$-sphere strictly included in a hemisphere. It does not contain antipodal points and is therefore free of conjugate points. 

\paragraph{Negative curvature} 

On the centered open disk of radius $R$, the isotropic metric 
\begin{align}
    g_{R,-}(x,y) := \frac{4R^4}{(x^2+y^2-R^2)^2},
    \label{eq:negcurvg}    
\end{align}
has constant negative Gaussian curvature $\kappa = -\frac{1}{R^2}$. When $R=1$, this is the model of the Poincar\'e disk. Because this model does not have a boundary (every geodesic has infinite length), one must choose a computational domain that is included in some disk of radius $R'<R$.

\paragraph{Experiments with simple domains}

Note that in both models above, when sending $R$ to $\infty$ while keeping a fixed bounded domain, the geometry becomes the Euclidean one. 

We now present one-shot inversions in cases of constant positive and negative curvature on a $n\times n$ grid with $n=300$. We will implement formula \eqref{eq:frc} here. The phantom is the non-smooth one appearing in e.g. Fig. \ref{fig:phantom1}, and the domains considered are
\begin{remunerate}
  \item The unit disk, of boundary equation $r(\beta) = 1$, see Fig. \ref{fig:phantom1}. 
  \item An ellipse, of boundary equation $r(\beta) = \frac{ab}{\sqrt{ (b\cos\beta)^2 + (a\sin\beta)^2}}$ with $(a,b) = (1,0.8)$, see Fig. \ref{fig:phantom2}. 
  \item A perturbation of the disk, of boundary equation $r(\beta) = a + b\cos(4\beta)$, with $(a,b) = (1,0.05)$, see Fig. \ref{fig:phantom3}. 
\end{remunerate}

\begin{figure}[htpb]
    \centering
    \subfigure[$R = 1.2$. Phantom/domain from Fig. \ref{fig:phantom1}. Relative $L^2$ error: $11.7\%$]{
	\includegraphics[width=0.21\textwidth]{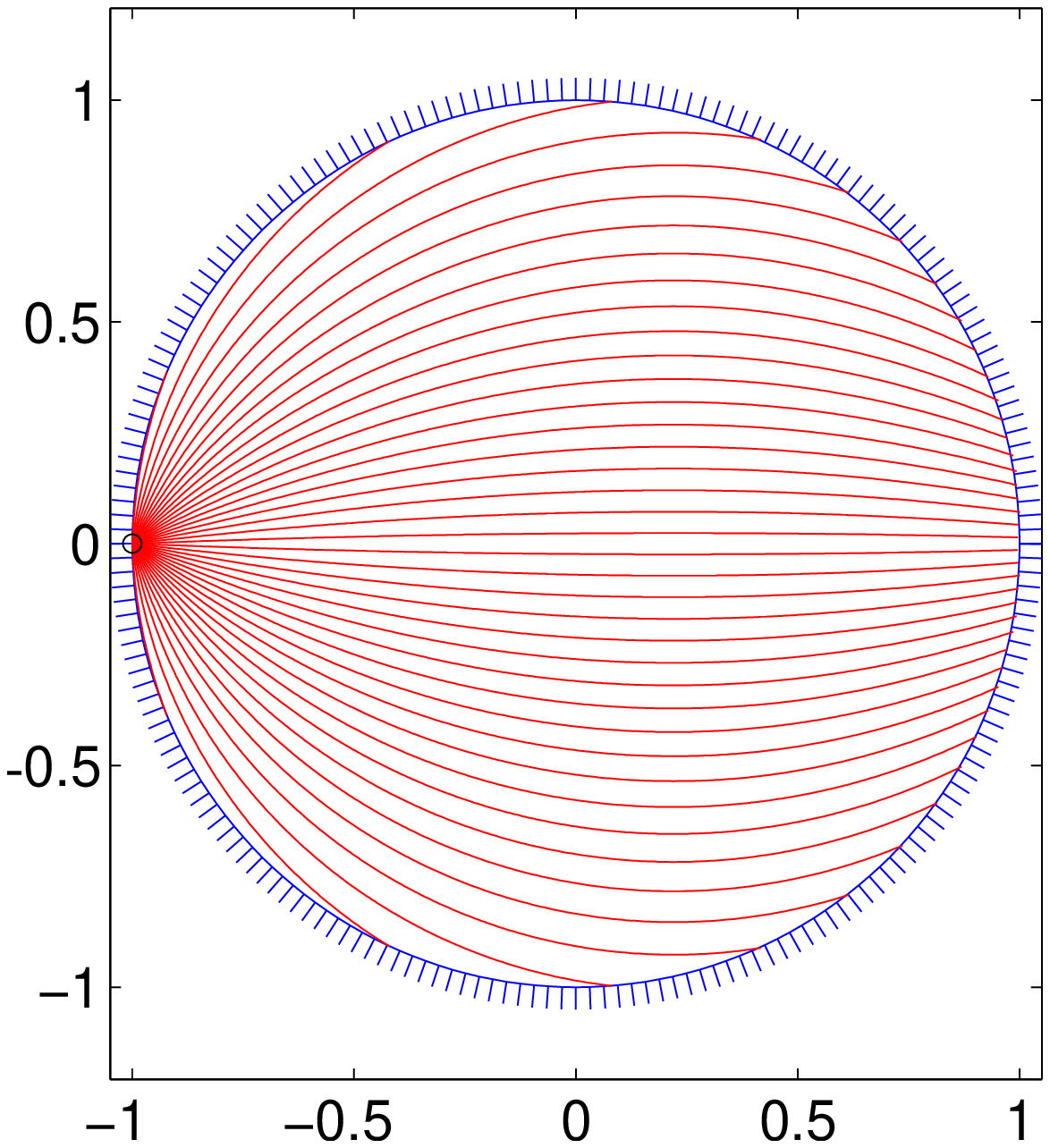}
	\includegraphics[width=0.48\textwidth]{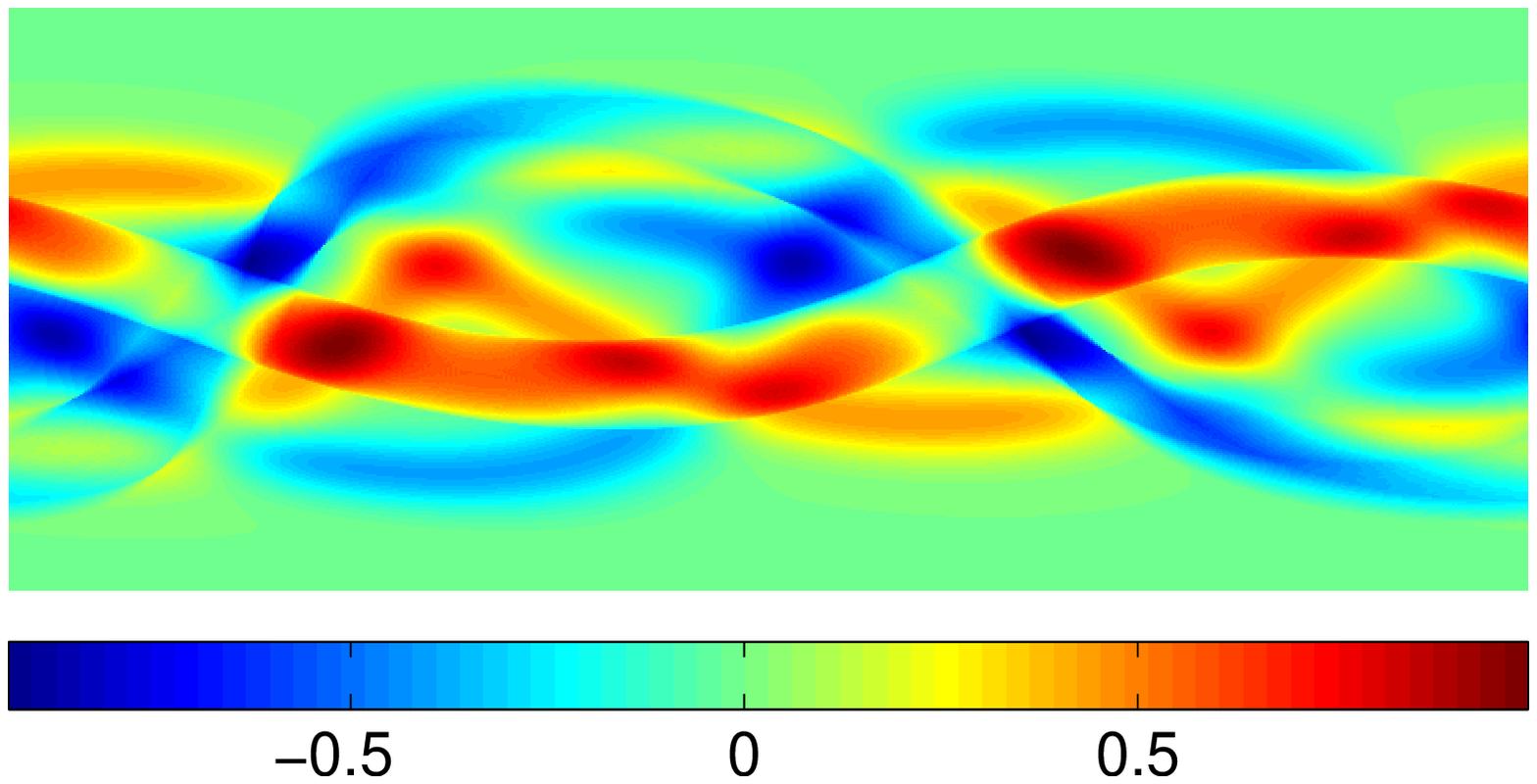}
	\includegraphics[width=0.25\textwidth]{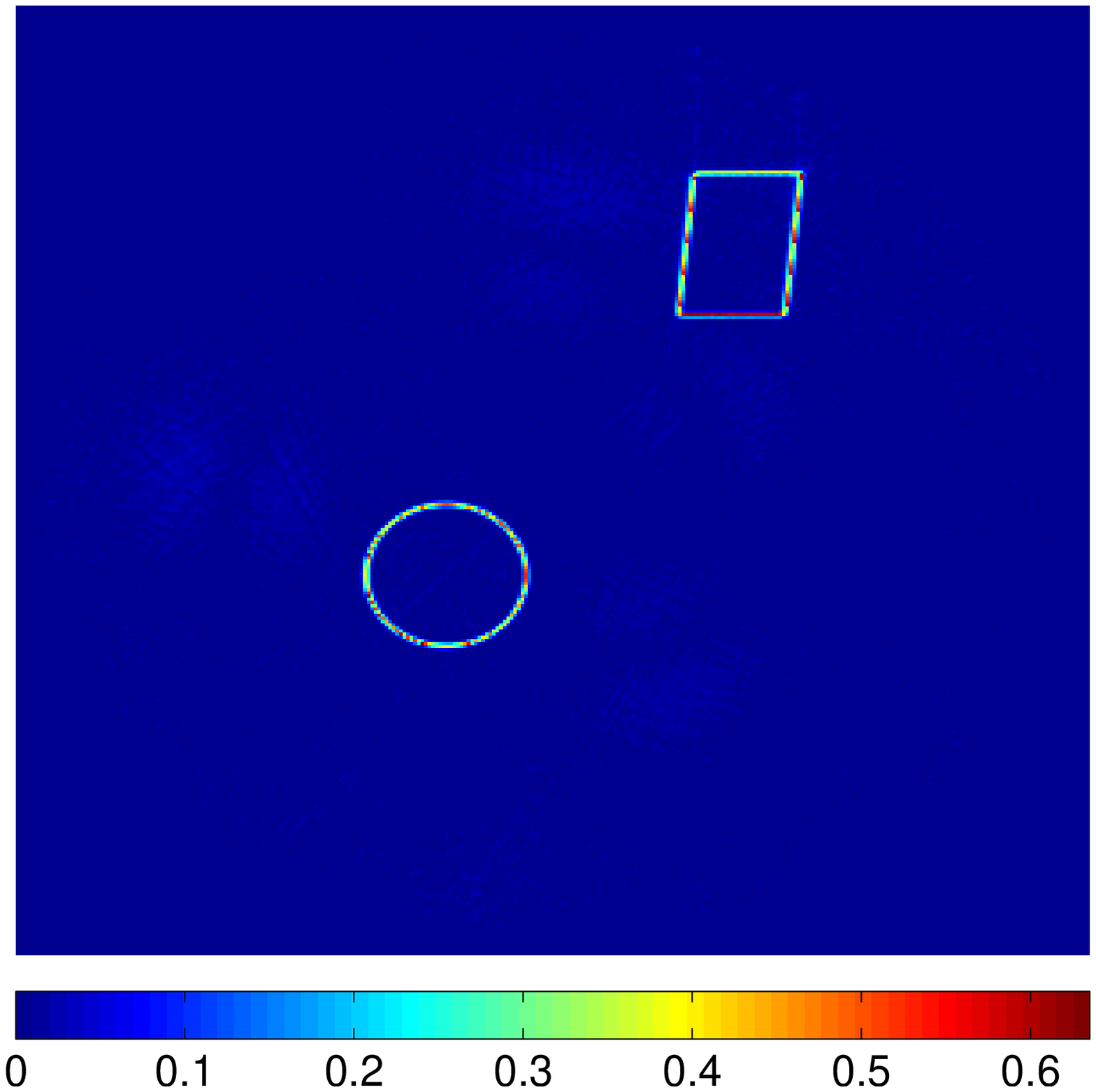}
	\label{fig:CPC1}
    }
    \subfigure[$R = 2$. Phantom/domain from Fig. \ref{fig:phantom1}. Relative $L^2$ error: $12.3\%$]{
	\includegraphics[width=0.21\textwidth]{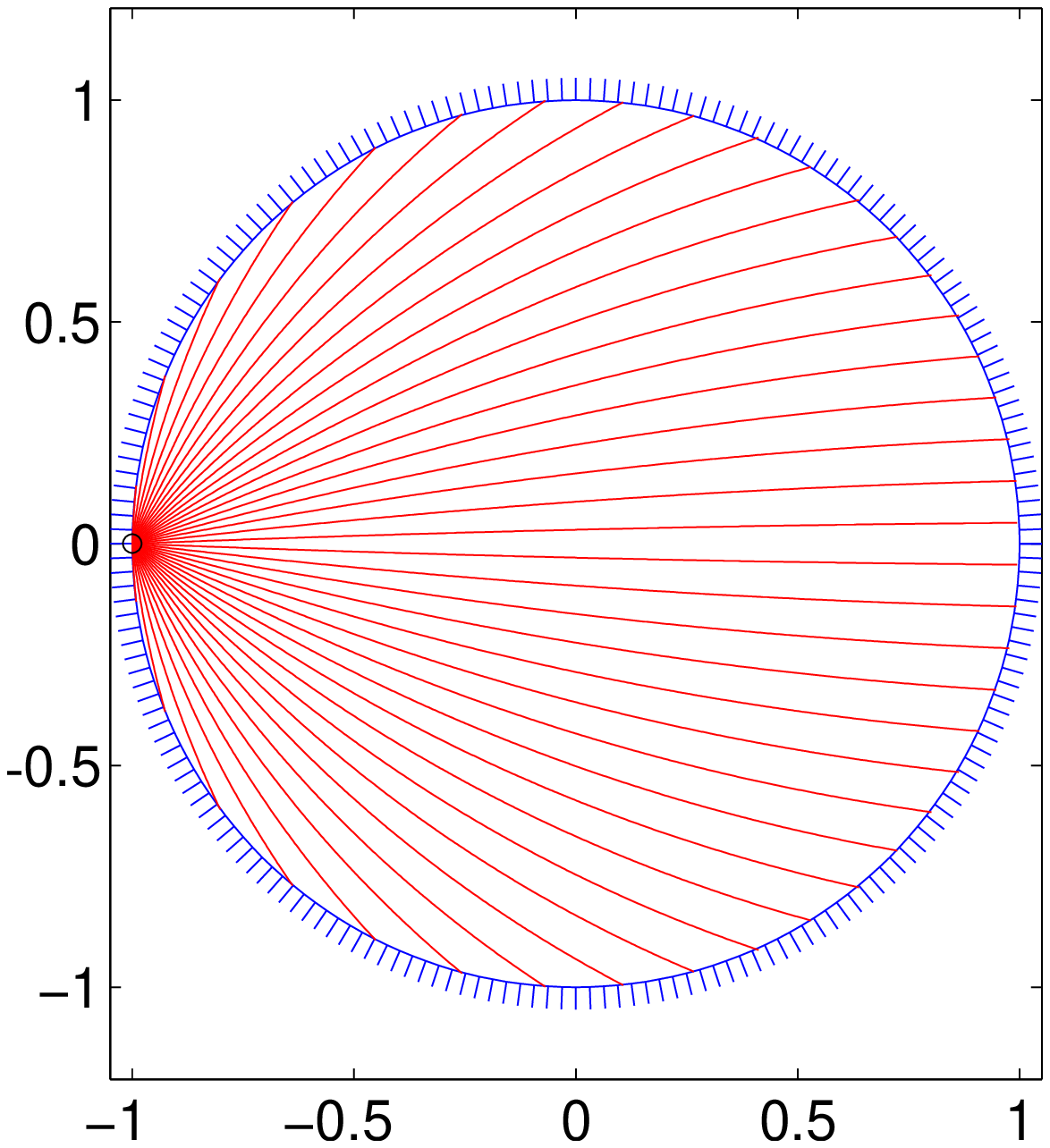}
	\includegraphics[width=0.48\textwidth]{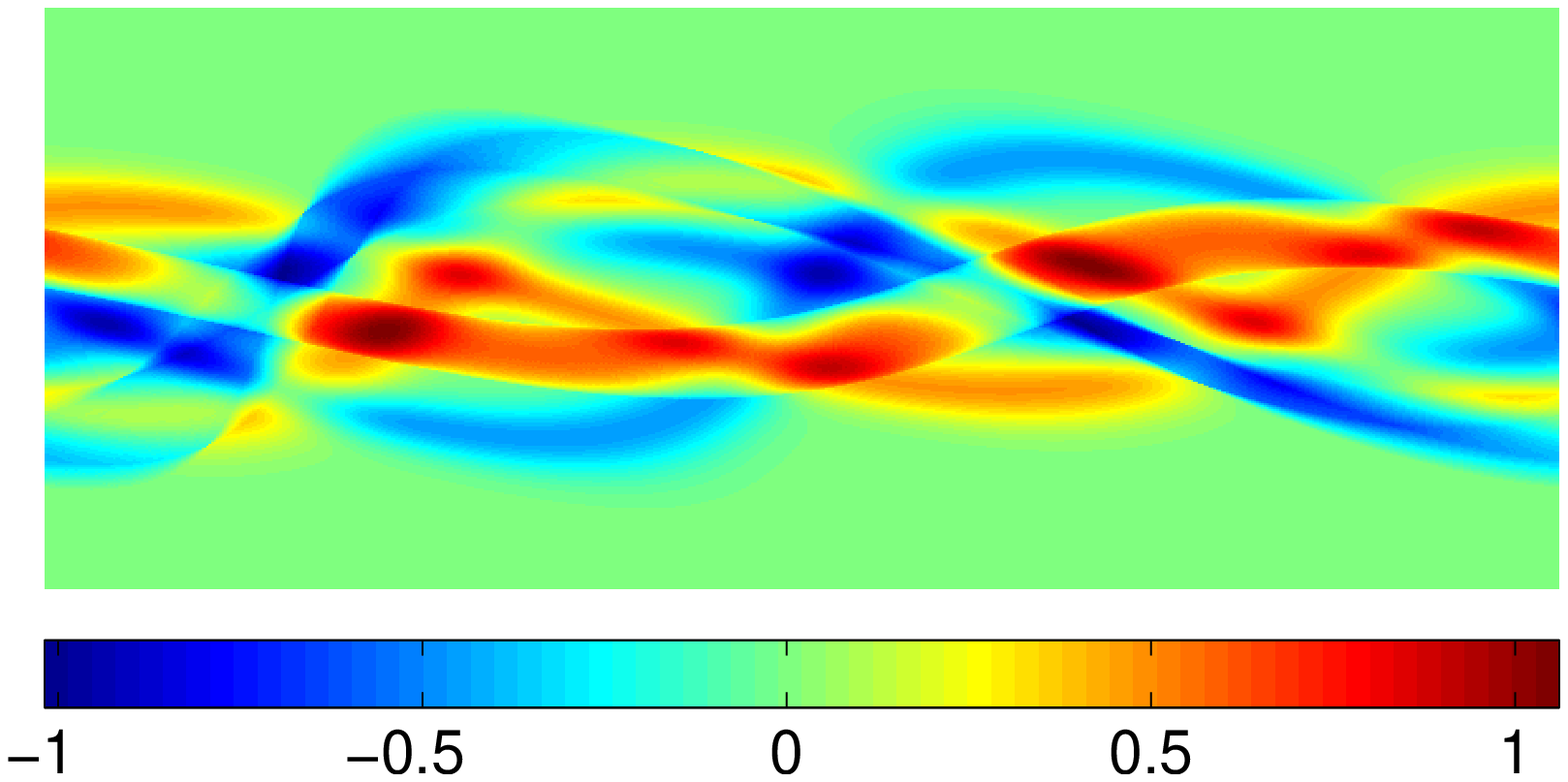}
	\includegraphics[width=0.25\textwidth]{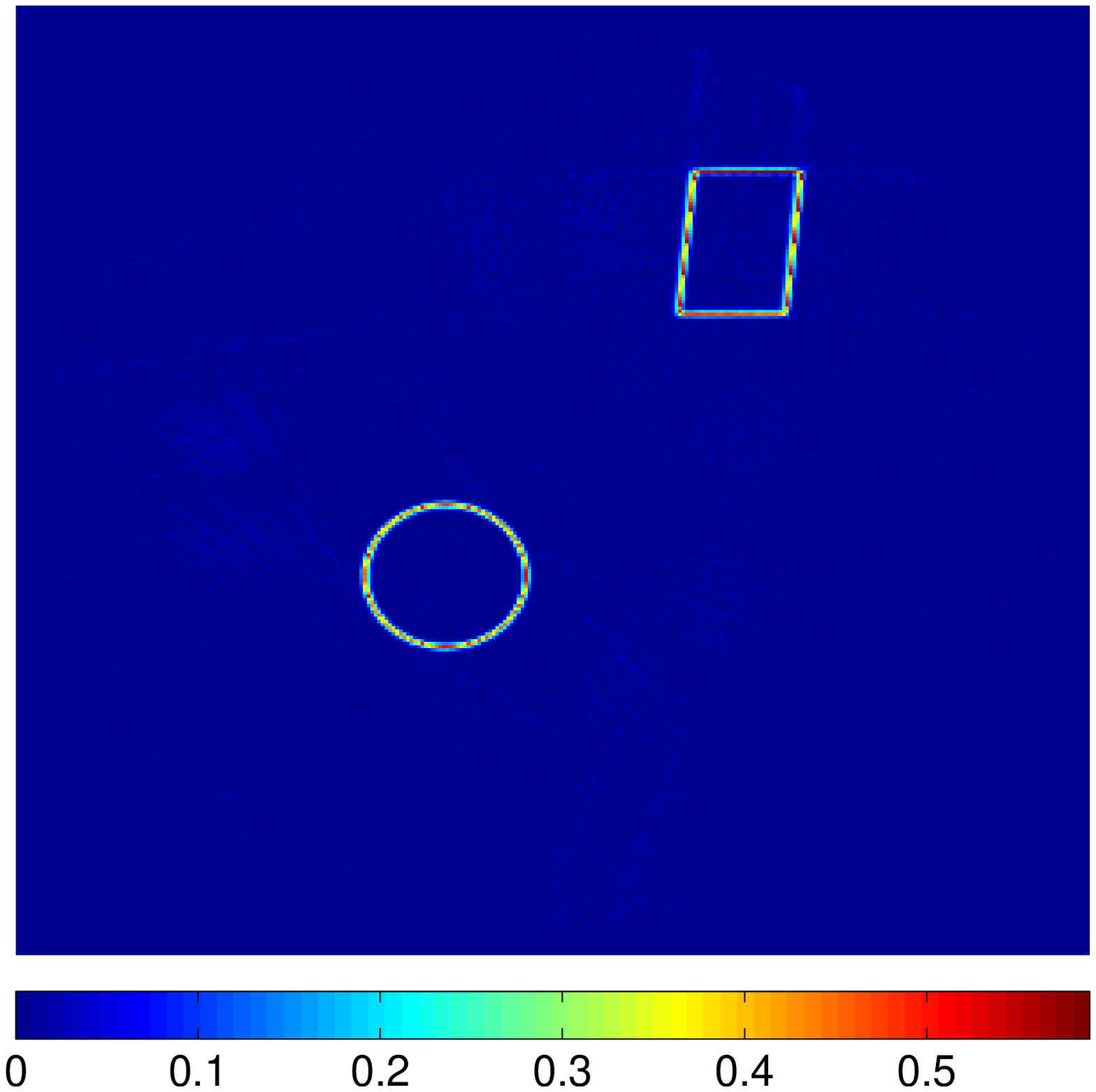}
    \label{fig:CPC2}
    }
    \subfigure[$R = 1.2$. Phantom/domain from Fig. \ref{fig:phantom2}. Relative $L^2$ error: $11.3\%$]{
	\includegraphics[width=0.21\textwidth]{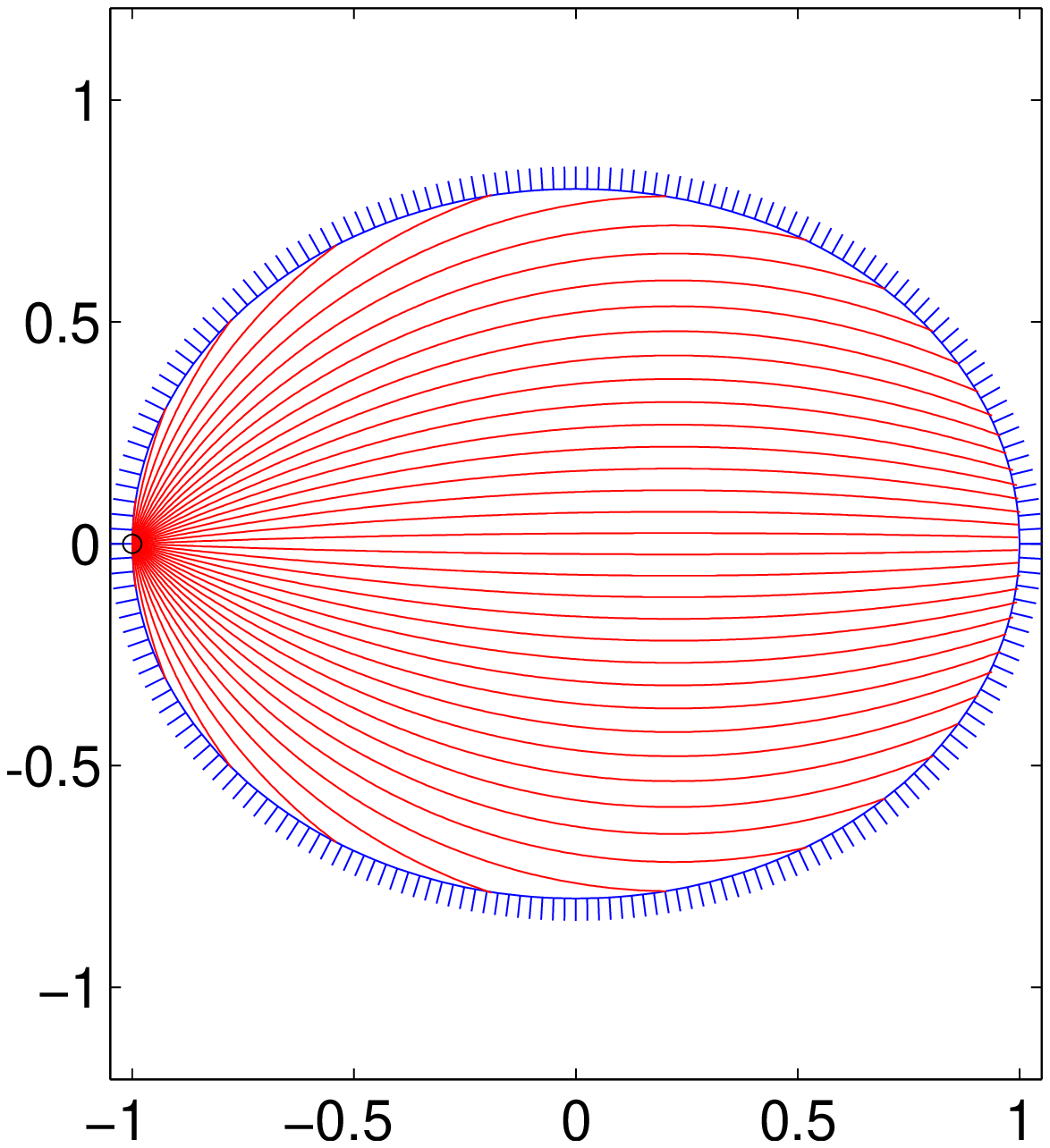}
	\includegraphics[width=0.48\textwidth]{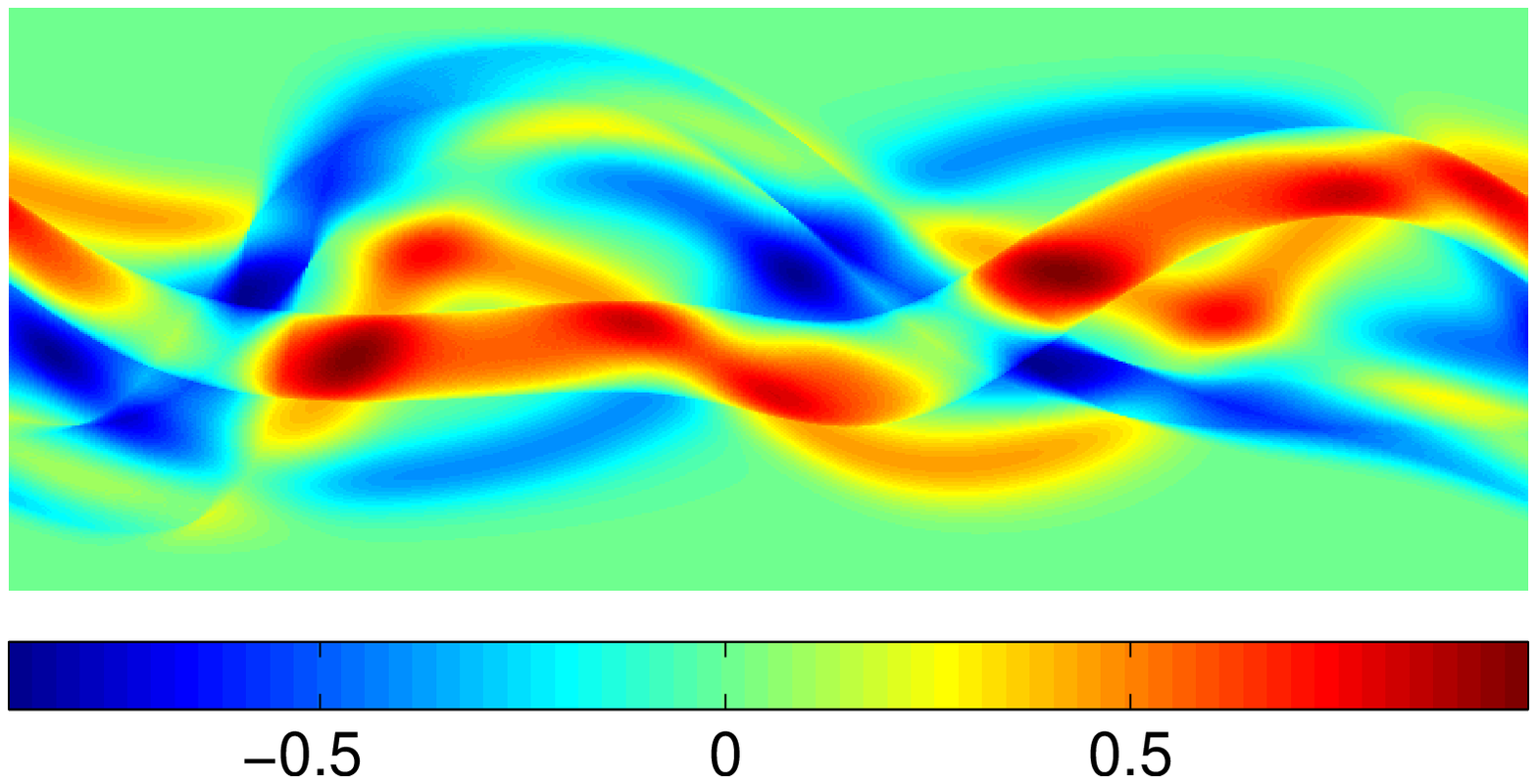}
	\includegraphics[width=0.25\textwidth]{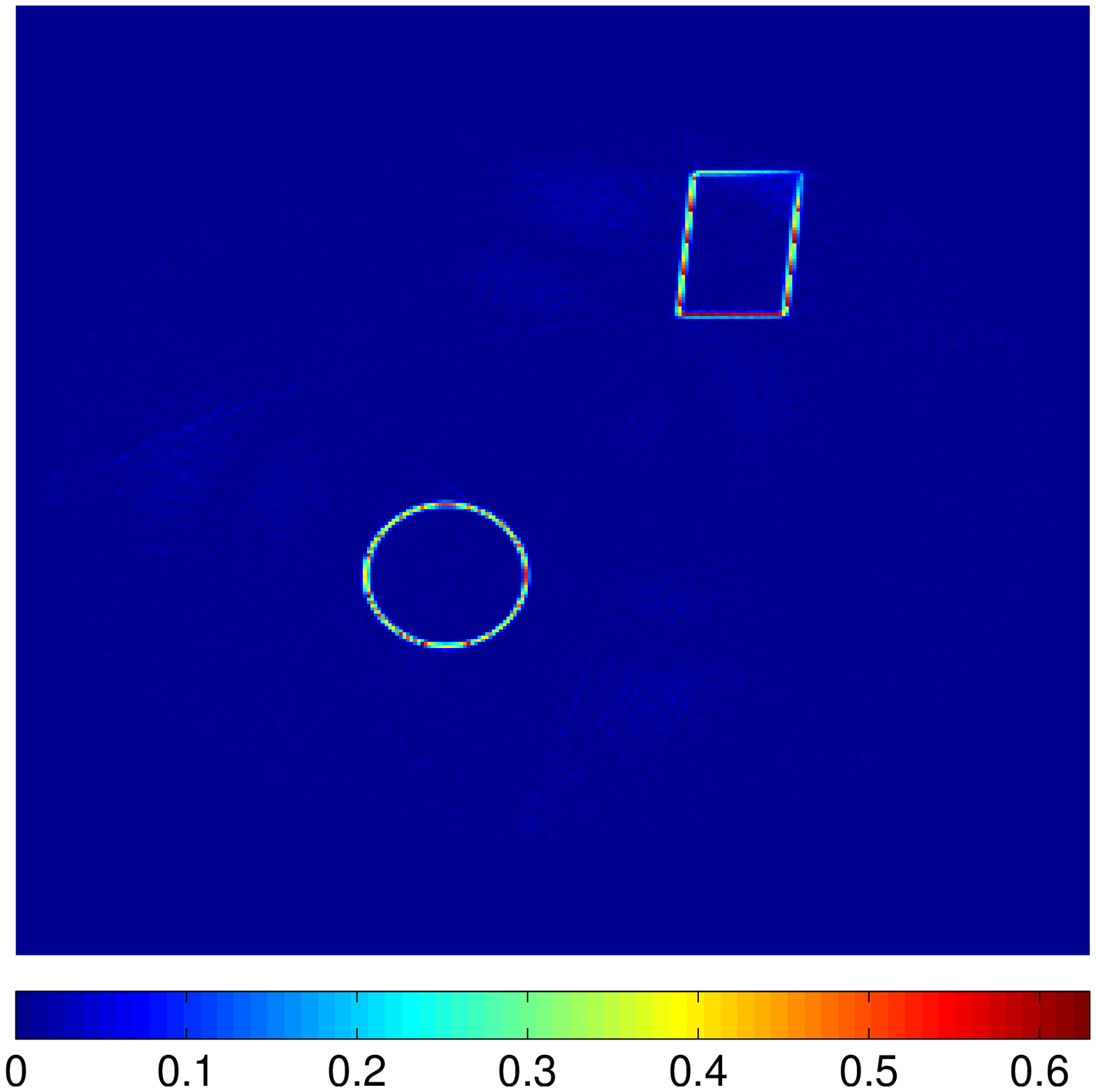}
	\label{fig:CPC3}
    }    
    \subfigure[$R = 2$. Phantom/domain from Fig. \ref{fig:phantom2}. Relative $L^2$ error: $11.7\%$]{
	\includegraphics[width=0.21\textwidth]{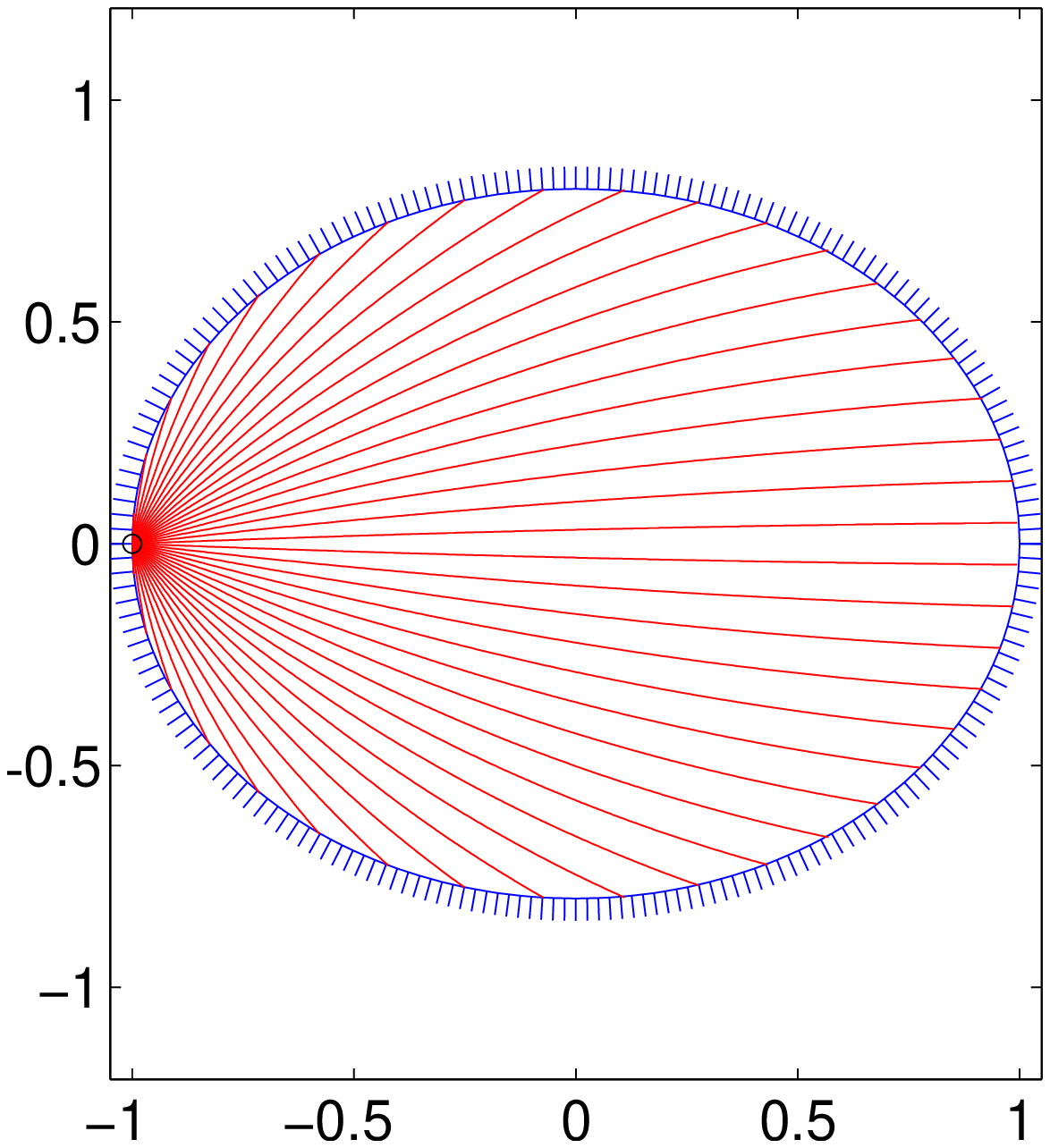}
	\includegraphics[width=0.48\textwidth]{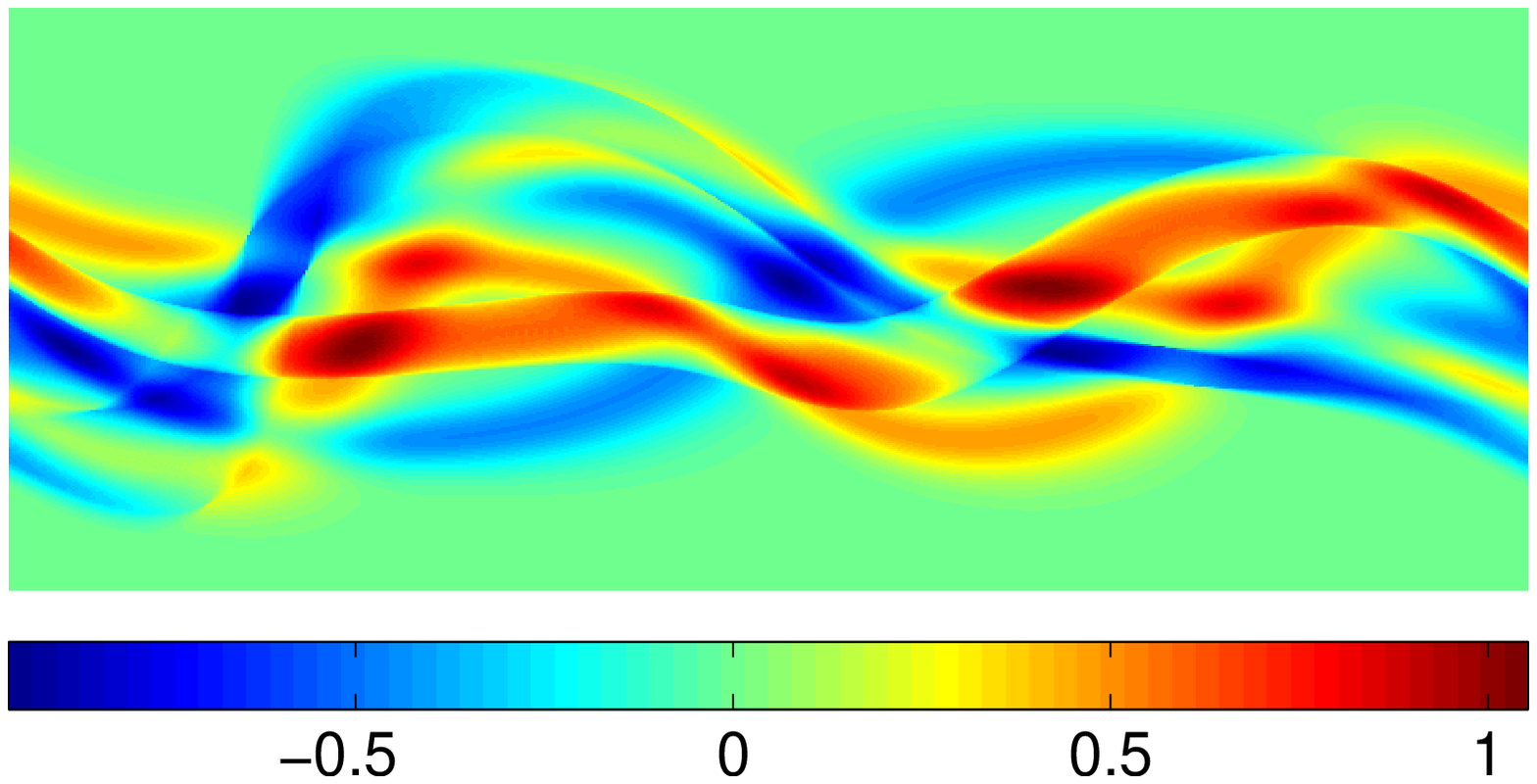}
	\includegraphics[width=0.25\textwidth]{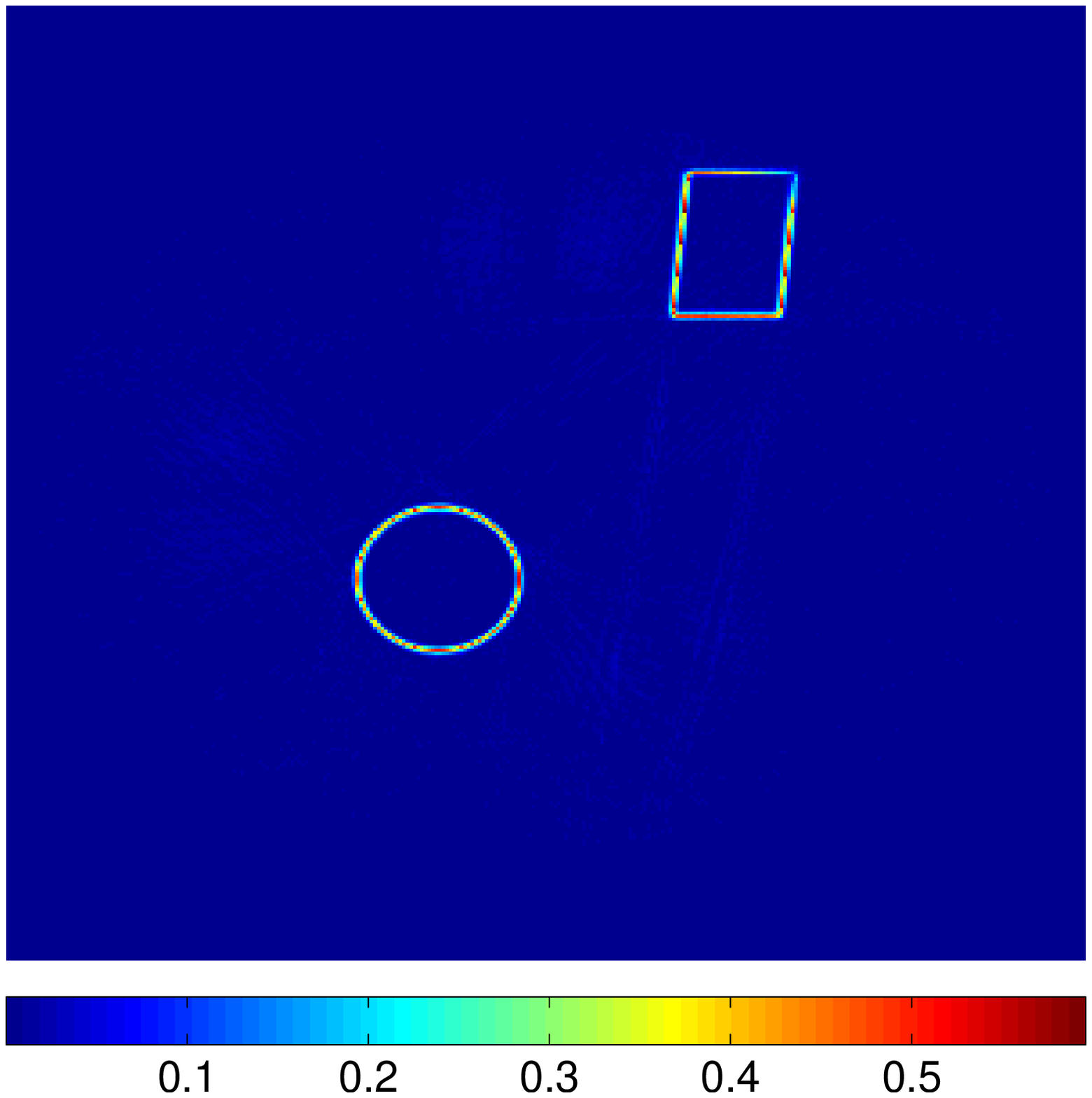}
	\label{fig:CPC4}
    }
    \caption{One-shot inversions in the constant positive curvature case. From left to right: some geodesics inside the domain, data $I_0 f$, pointwise error $|f-f_{rc}|$ after one-shot inversion.}
    \label{fig:CPC}
\end{figure}

\begin{figure}[htpb]
    \centering
    \subfigure[$R = 1.2$. Phantom/domain from Fig. \ref{fig:phantom1}. Relative $L^2$ error: $20.3\%$]{
	\includegraphics[width=0.21\textwidth]{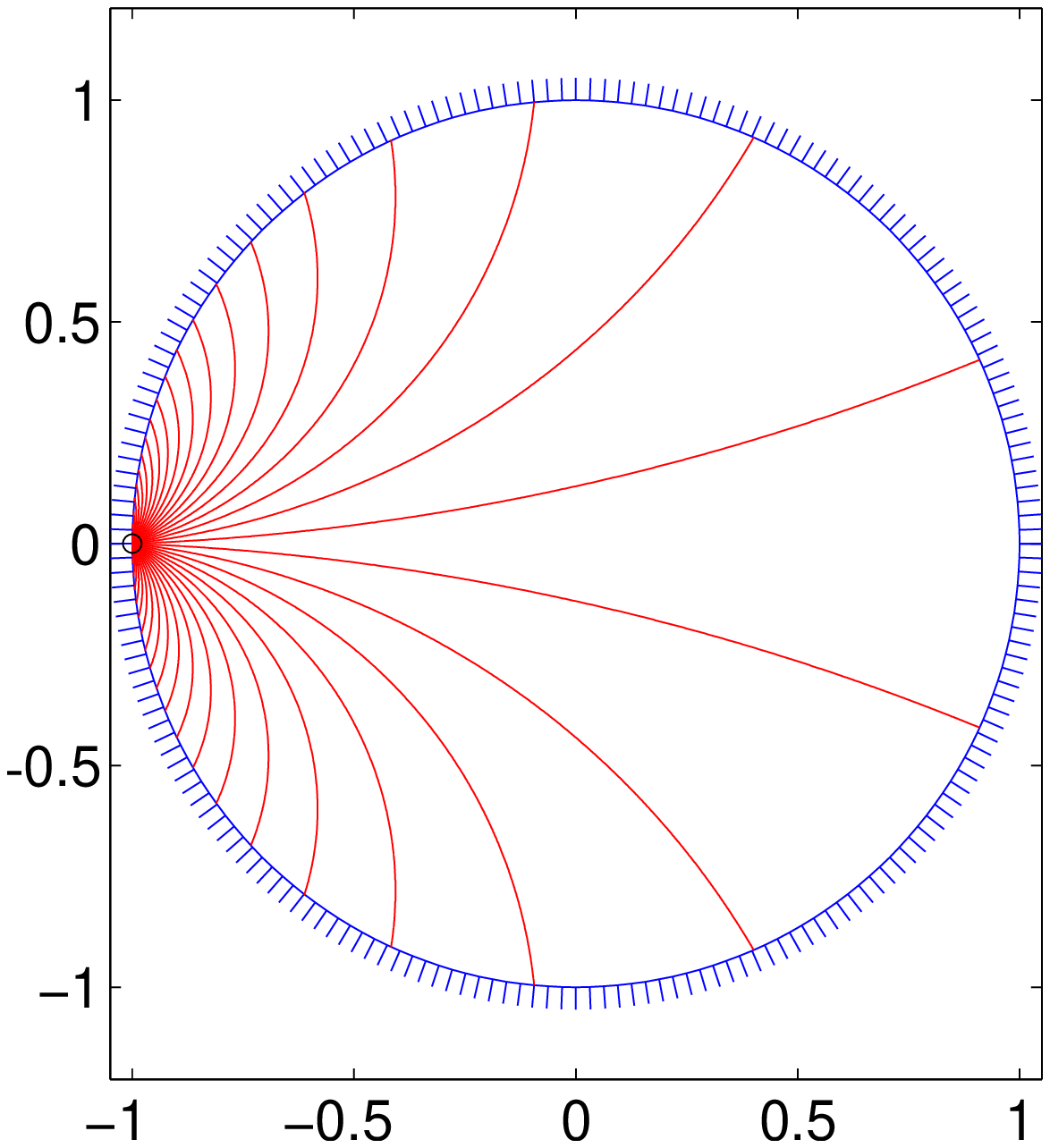}
	\includegraphics[width=0.48\textwidth]{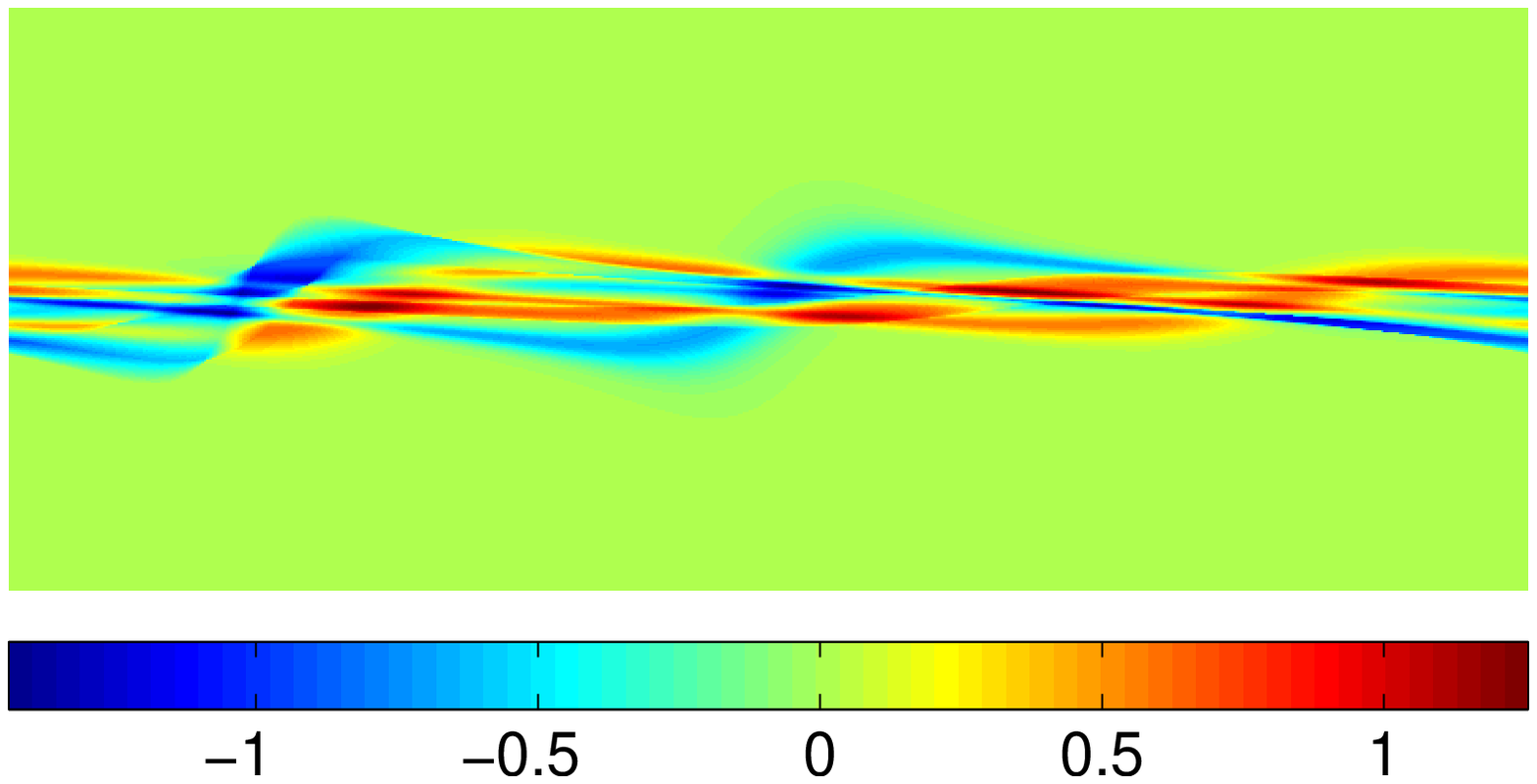}
	\includegraphics[width=0.25\textwidth]{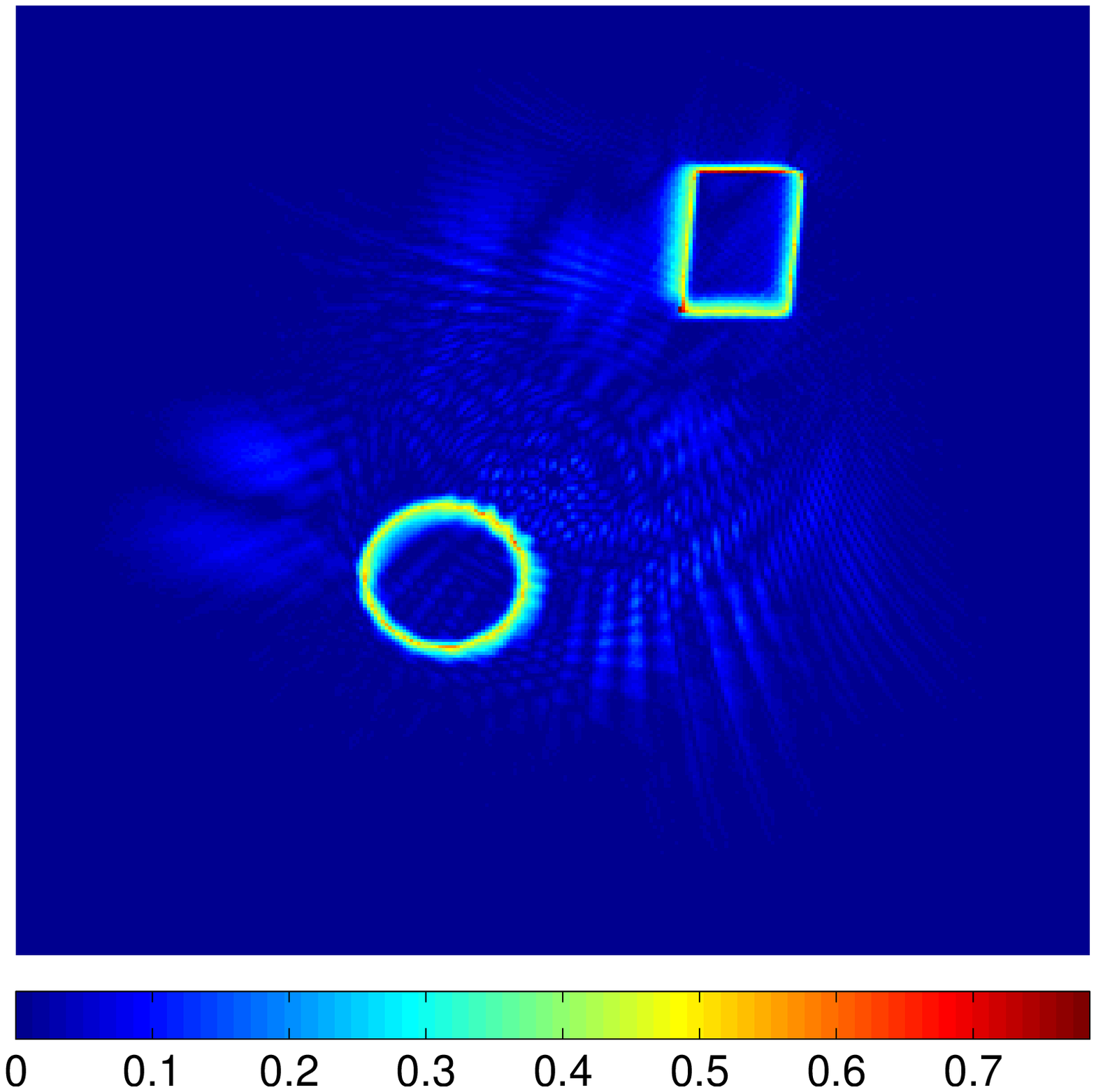}
	\label{fig:CNC1}
    }
    \subfigure[$R = 2$. Phantom/domain from Fig. \ref{fig:phantom1}. Relative $L^2$ error: $15.2\%$]{
	\includegraphics[width=0.21\textwidth]{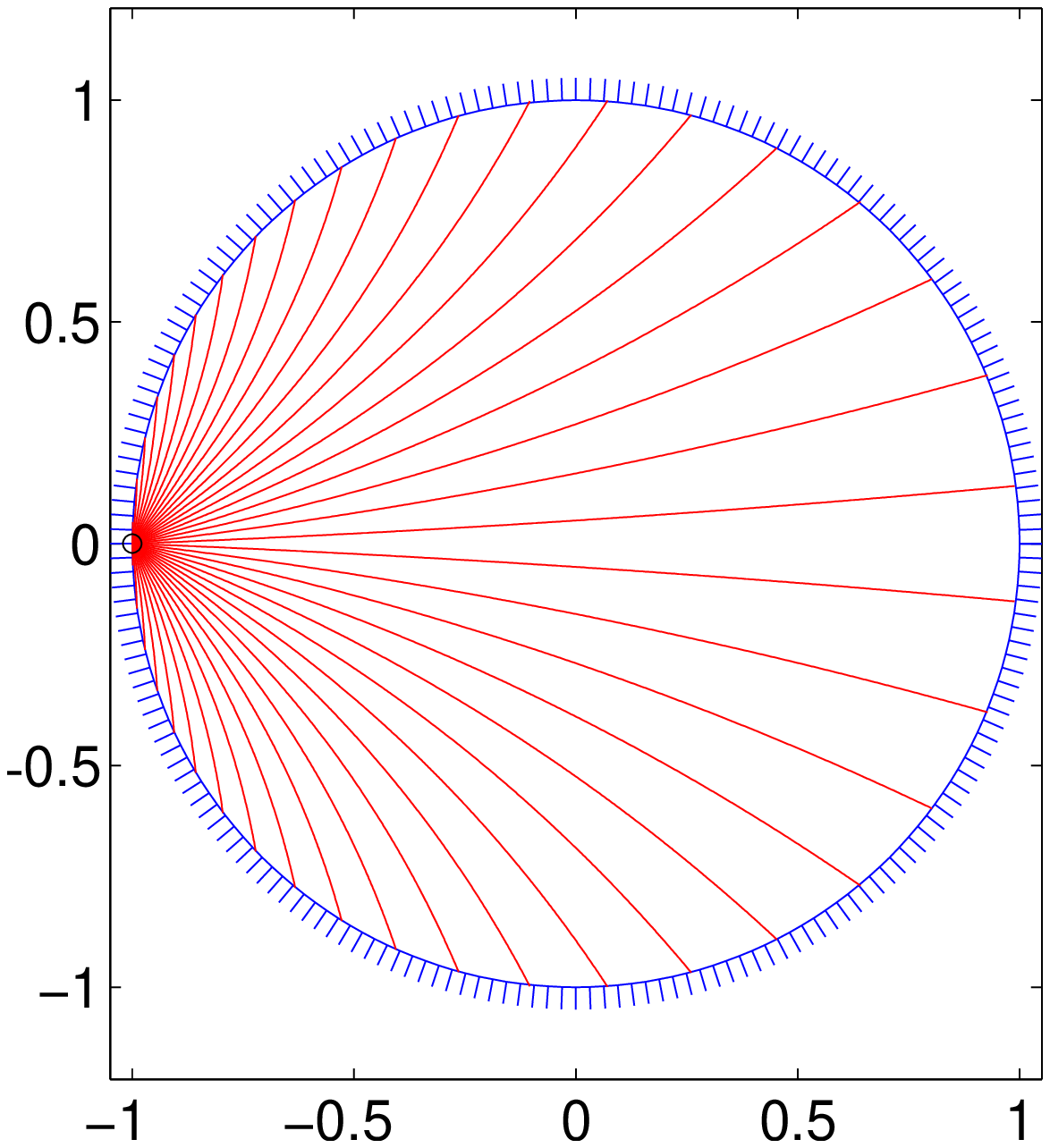}
	\includegraphics[width=0.48\textwidth]{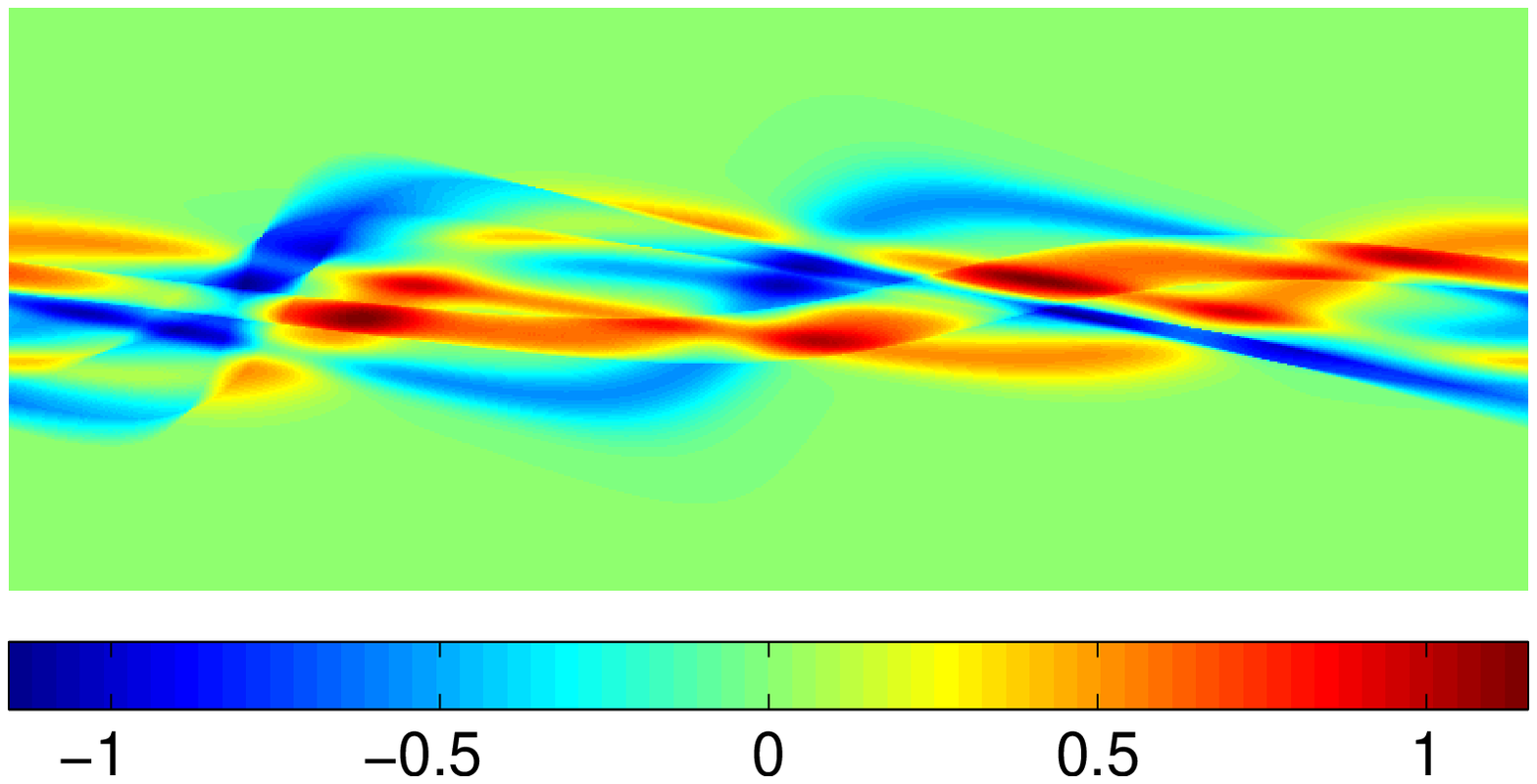}
	\includegraphics[width=0.25\textwidth]{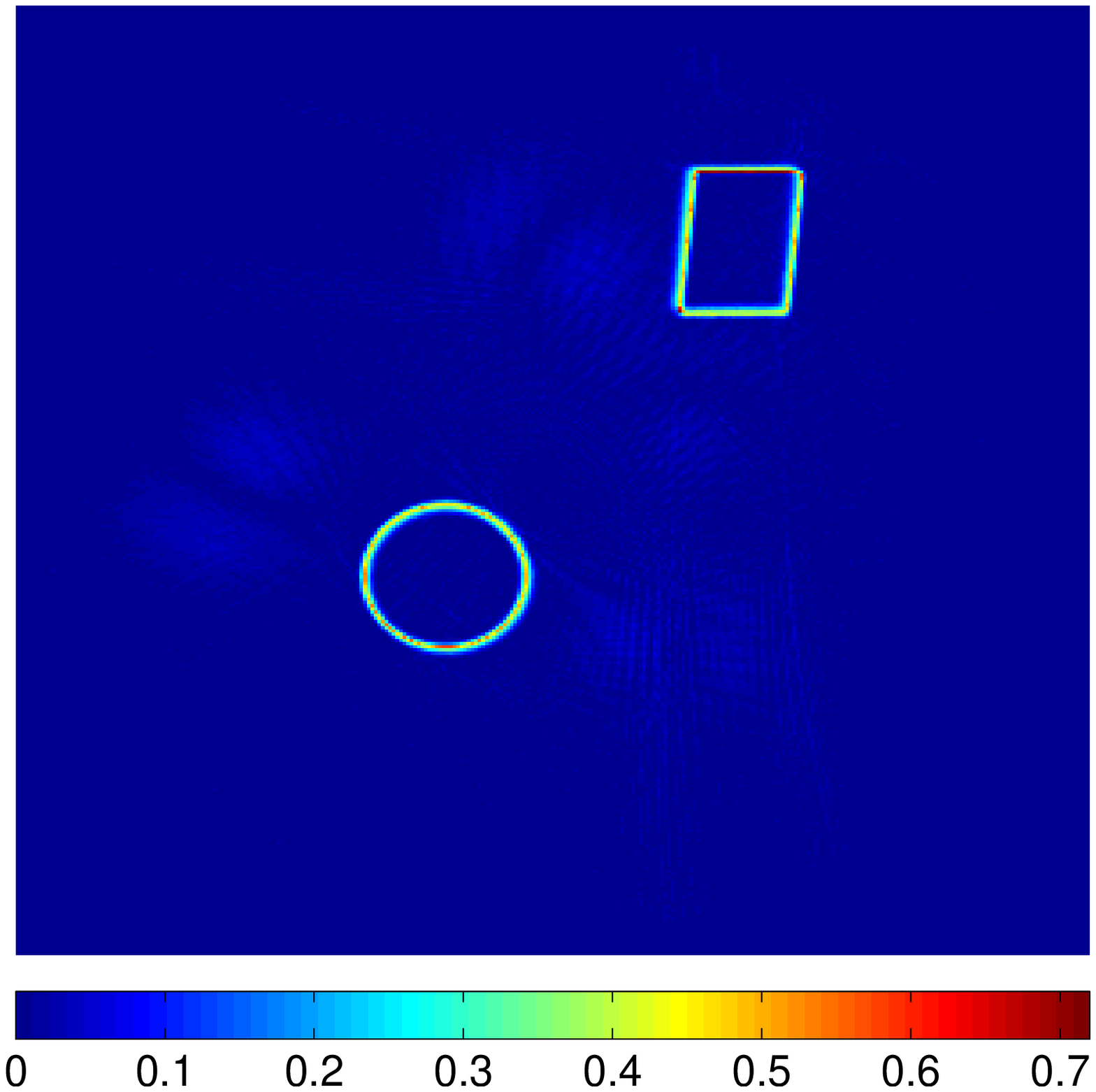}
	\label{fig:CNC2}
    }
    \subfigure[$R = 1.2$. Phantom/domain from Fig. \ref{fig:phantom3}. Relative $L^2$ error: $20.5\%$]{
	\includegraphics[width=0.21\textwidth]{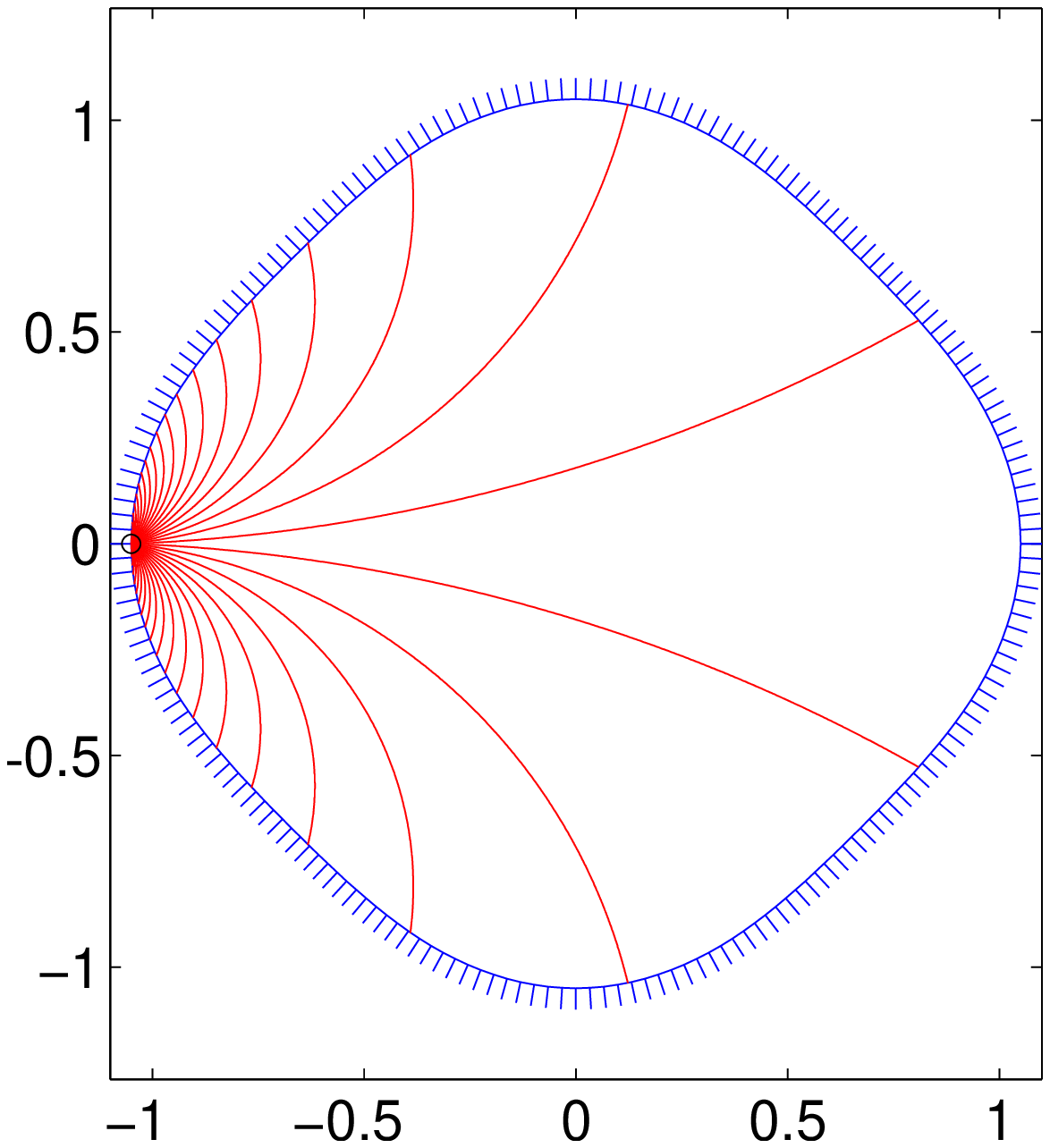}
	\includegraphics[width=0.48\textwidth]{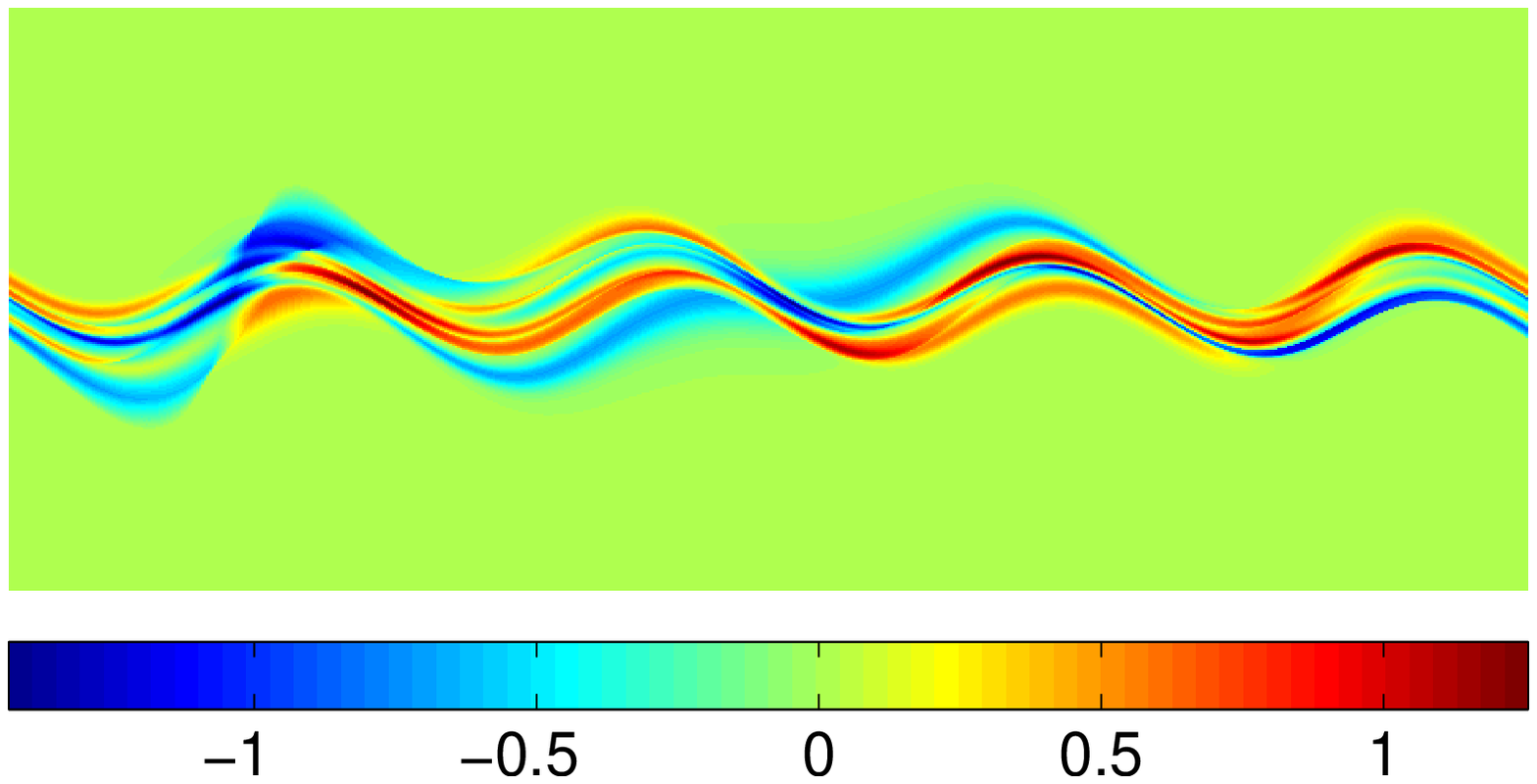}
	\includegraphics[width=0.25\textwidth]{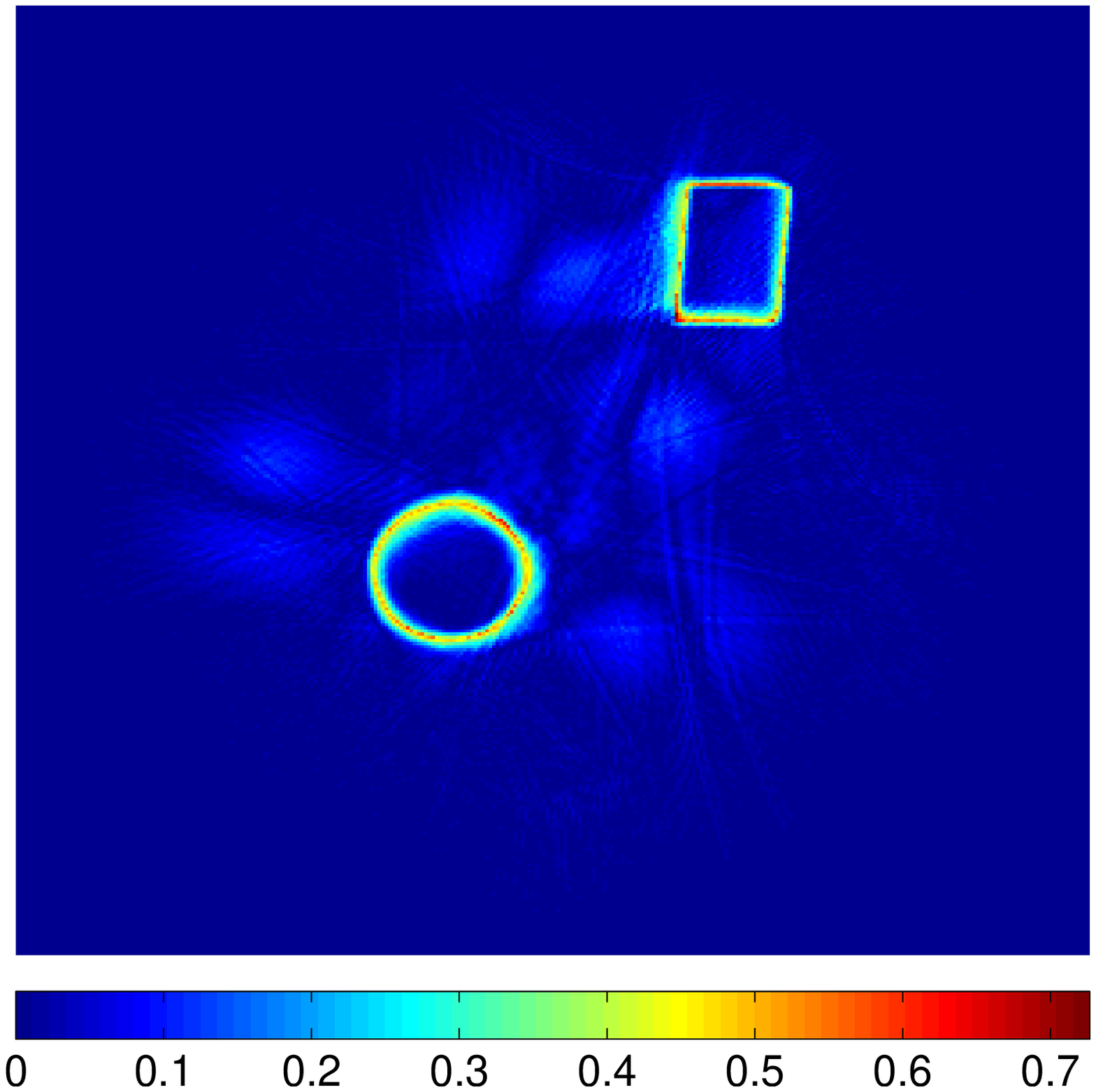}
	\label{fig:CNC3}
    }    
    \subfigure[$R = 2$. Phantom/domain from Fig. \ref{fig:phantom3}. Relative $L^2$ error: $15.6\%$]{
	\includegraphics[width=0.21\textwidth]{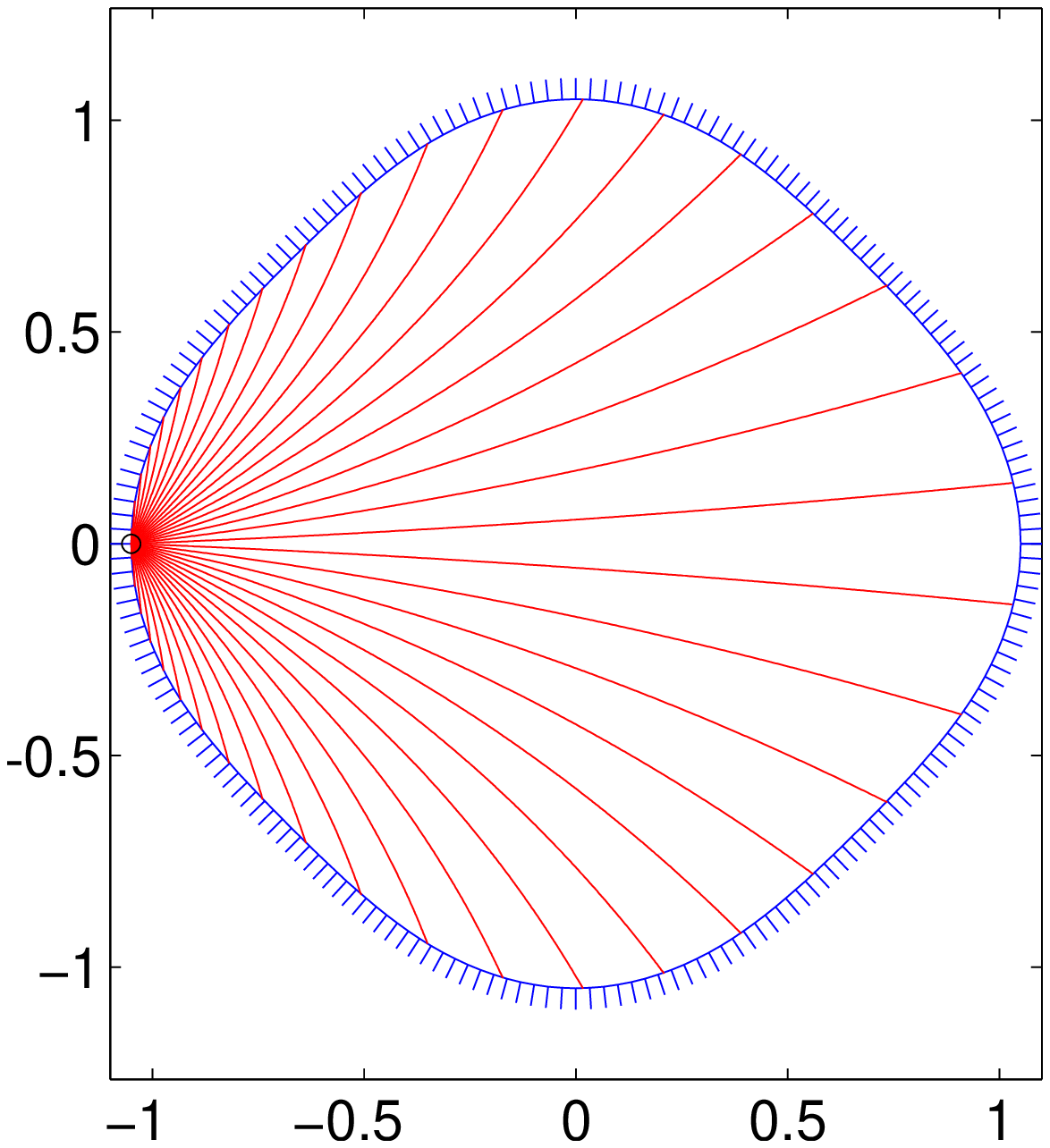}
	\includegraphics[width=0.48\textwidth]{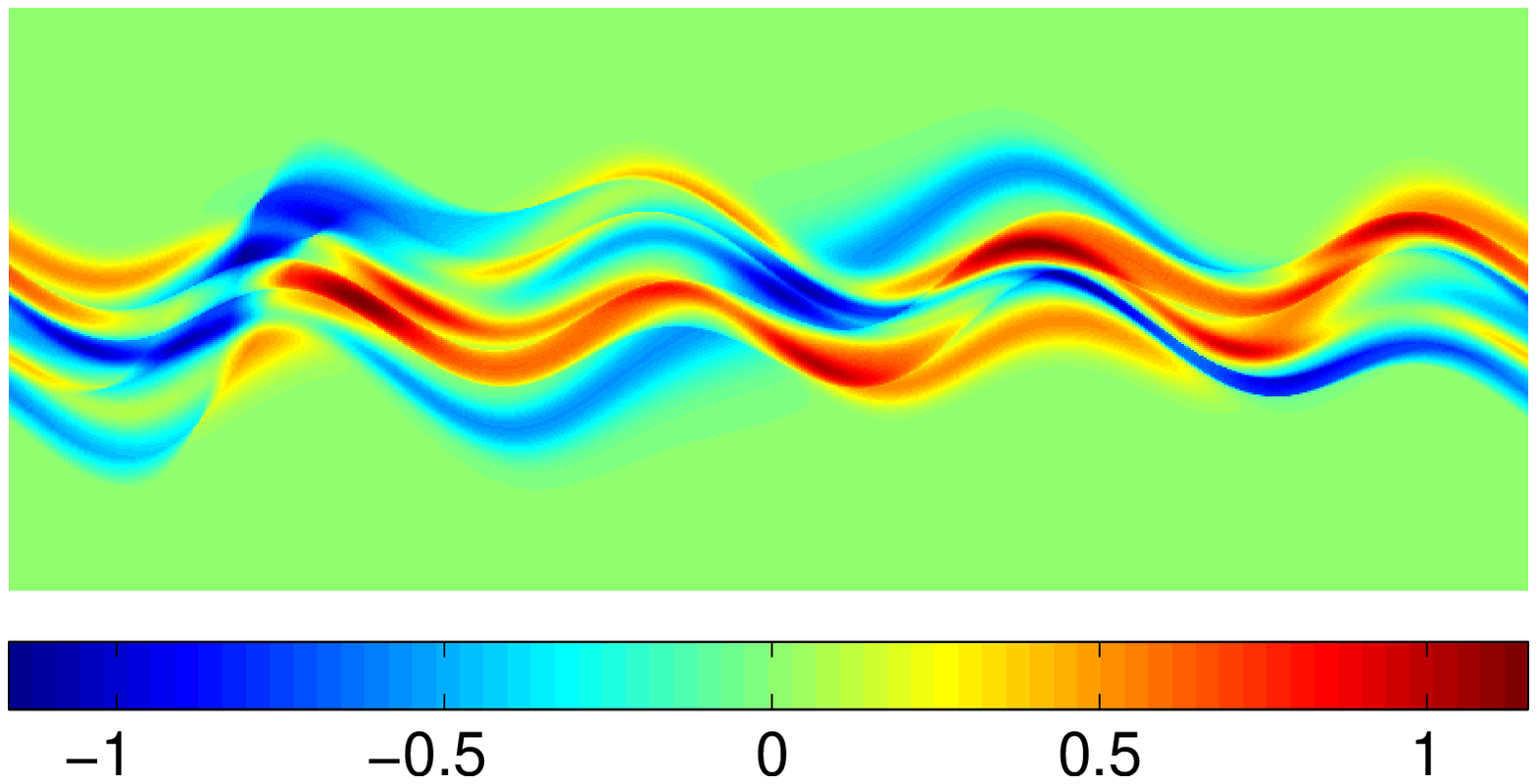}
	\includegraphics[width=0.25\textwidth]{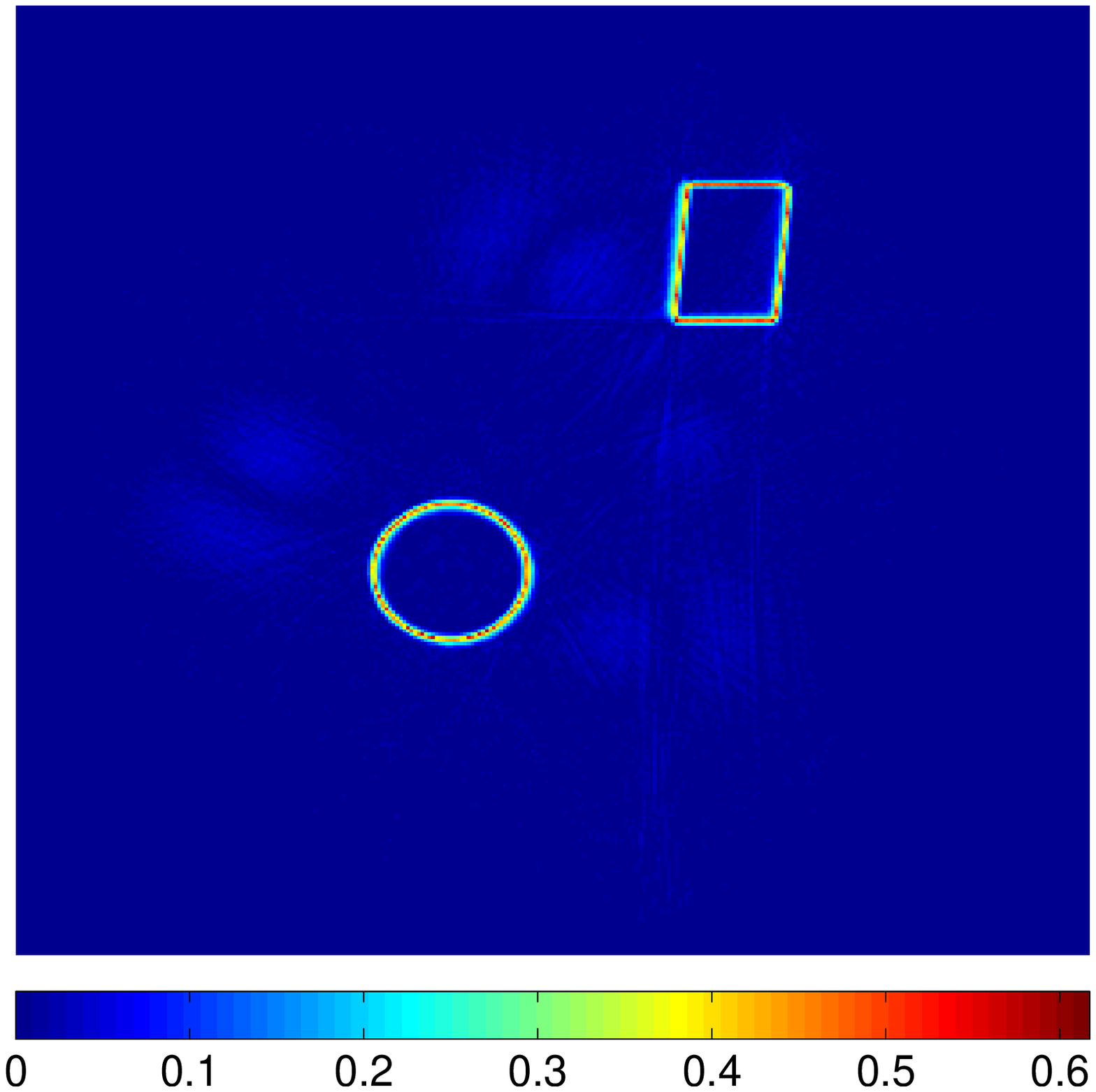}
	\label{fig:CNC4}
    }
    \caption{One-shot inversions in the constant negative curvature case. From left to right: some geodesics inside the domain, data $I_0 f$, pointwise error $|f-f_{rc}|$ after one-shot inversion.}
    \label{fig:CNC}
\end{figure}

Some one-shot reconstructions are presented on Figure \ref{fig:CPC} for the case of positive curvature and on Figure \ref{fig:CNC} for the case of negative curvature. In each case, the smooth part is accurately recovered while the error is concentrated at the sharp edges, as they cannot be resolved with perfect accuracy. In each group of pictures, the left plot has 40 curves shot from the leftmost point of the domain, with equispacing in direction. Looking at the negative curvature plots with $R = 1.2$ (Figs. \ref{fig:CNC1} and \ref{fig:CNC3}), we see that in comparison to the other cases, many fewer curves actually sample the central part of the domain. This is responsible for the data being supported close to the axis $\alpha=0$ and the undersampling artifacts on the corresponding pointwise error images. 

Based on this observation, it becomes clear that the uniform sampling in $\alpha$ is not adapted to every case of metric. This is an issue that will be adressed in future work. 

\begin{figure}[htpb]
    \centering
    \subfigure[Phantom $f$ and domain]{
	\includegraphics[height=0.18\textheight]{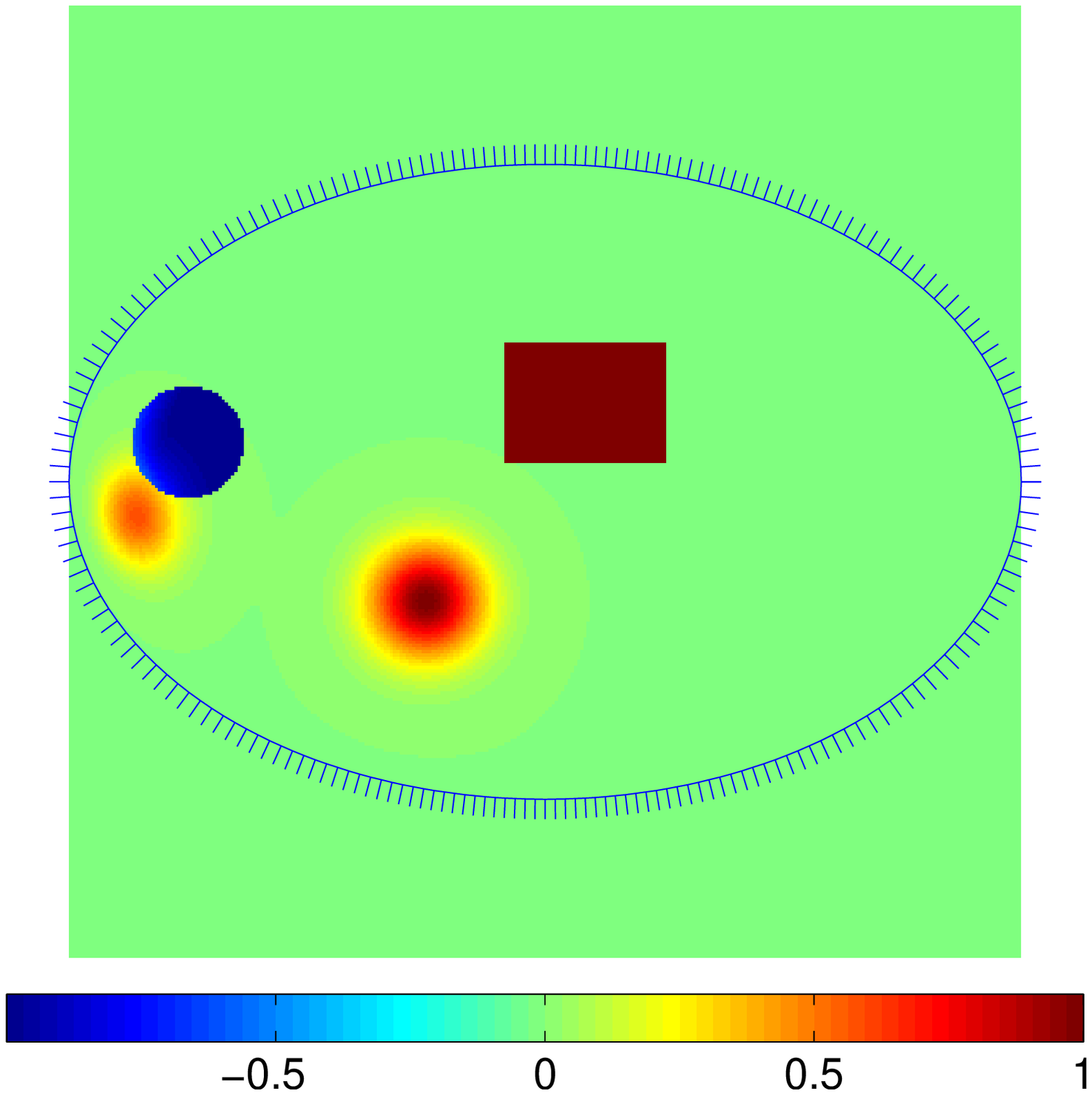}
	\label{fig:CPCns1}
    }
    \subfigure[Sample geodesics for $g_{R,+}$ with $R=1$]{
	\includegraphics[height=0.18\textheight]{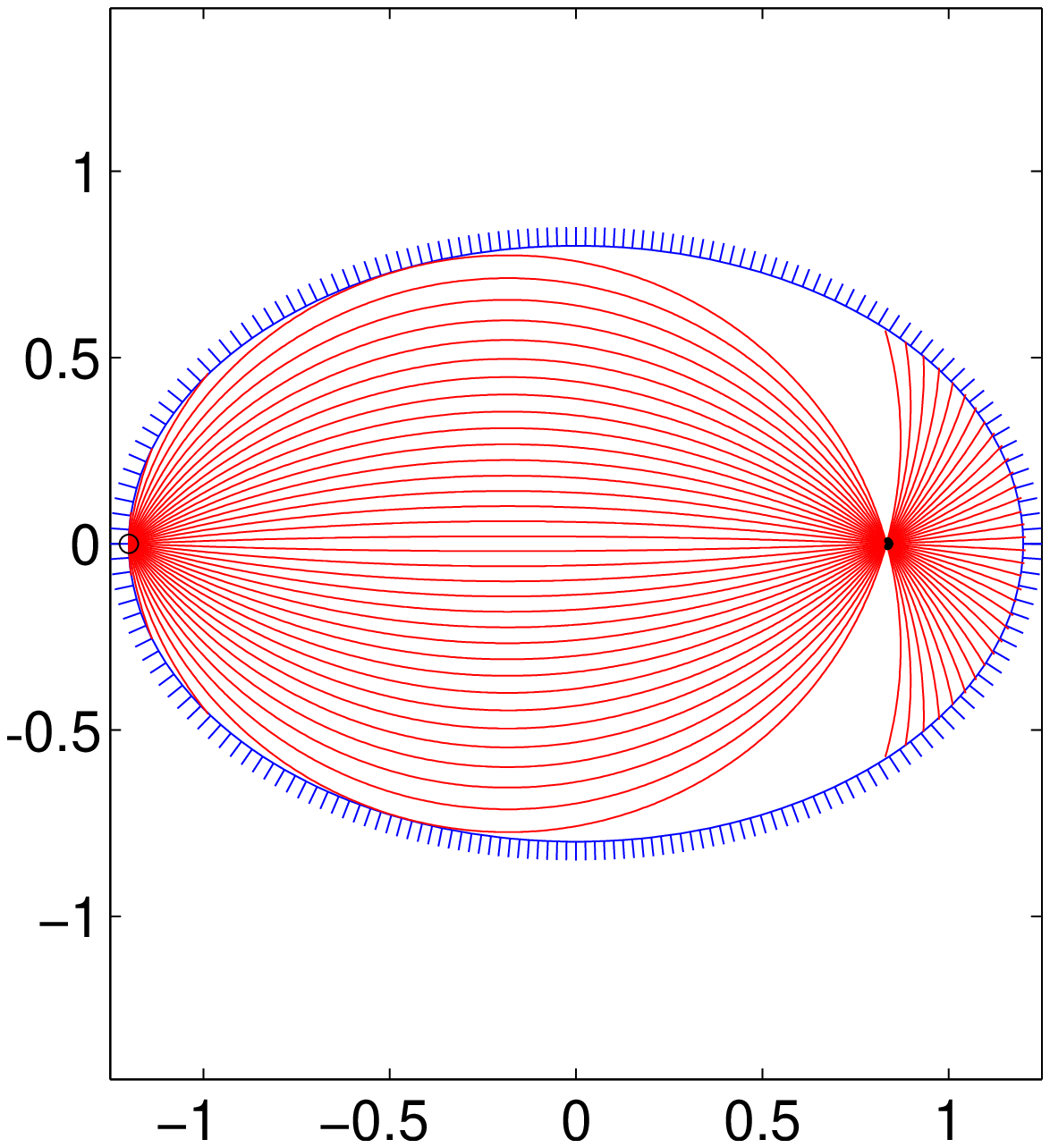} 
	\includegraphics[height=0.18\textheight]{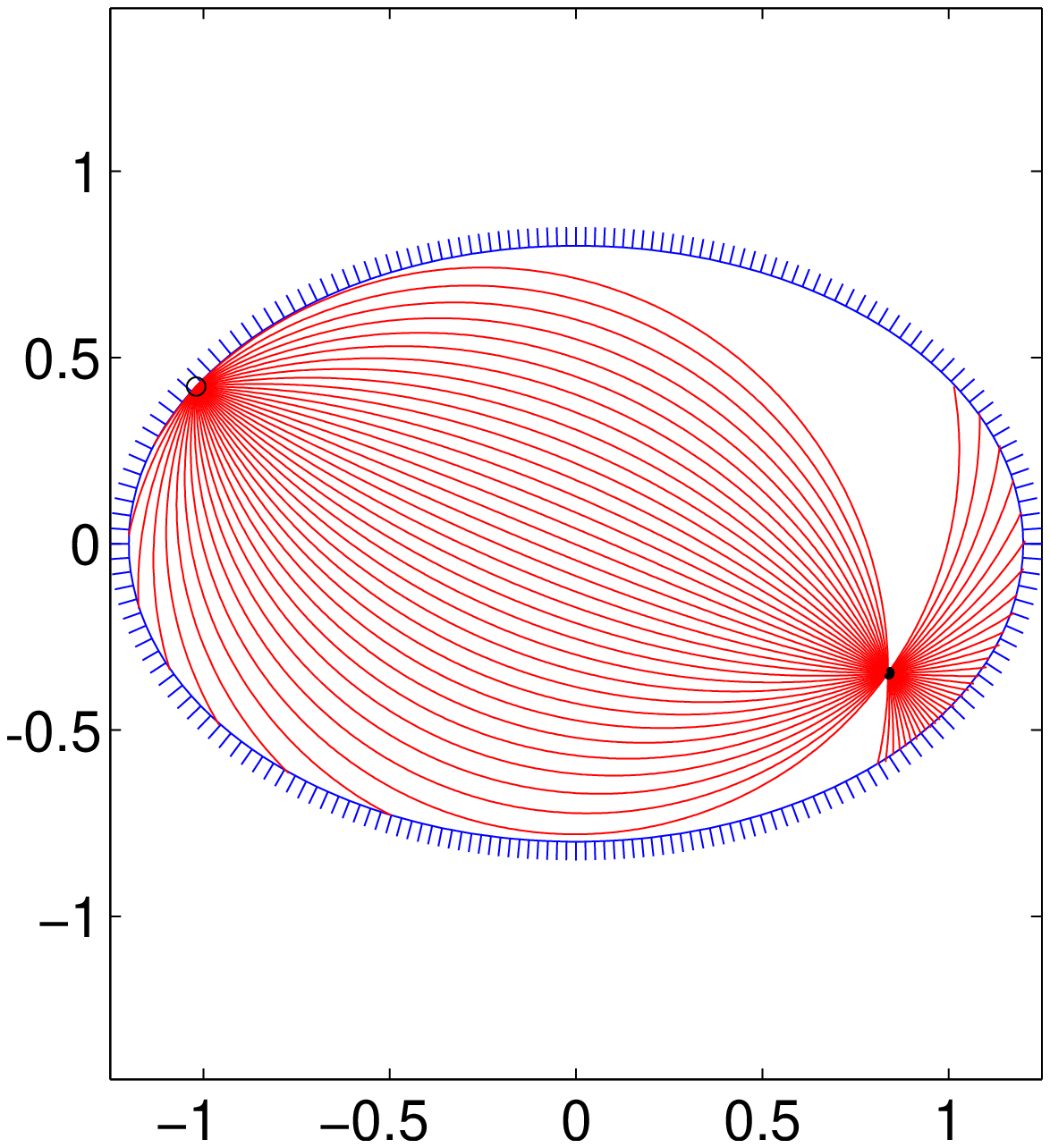}
	\includegraphics[height=0.18\textheight]{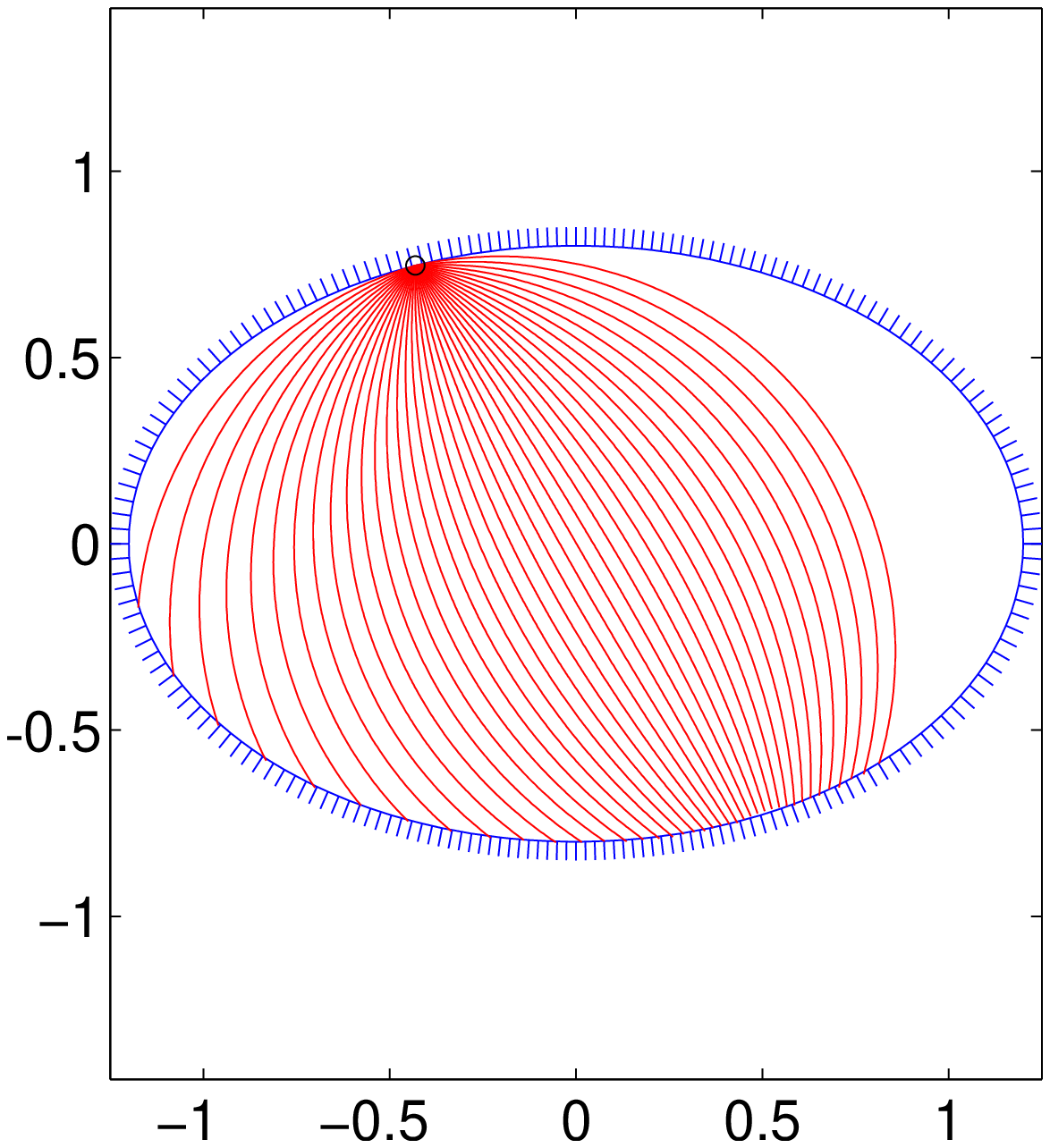}  
	\label{fig:CPCns2}
    }
    \subfigure[Forward data]{
	\includegraphics[height=0.18\textheight]{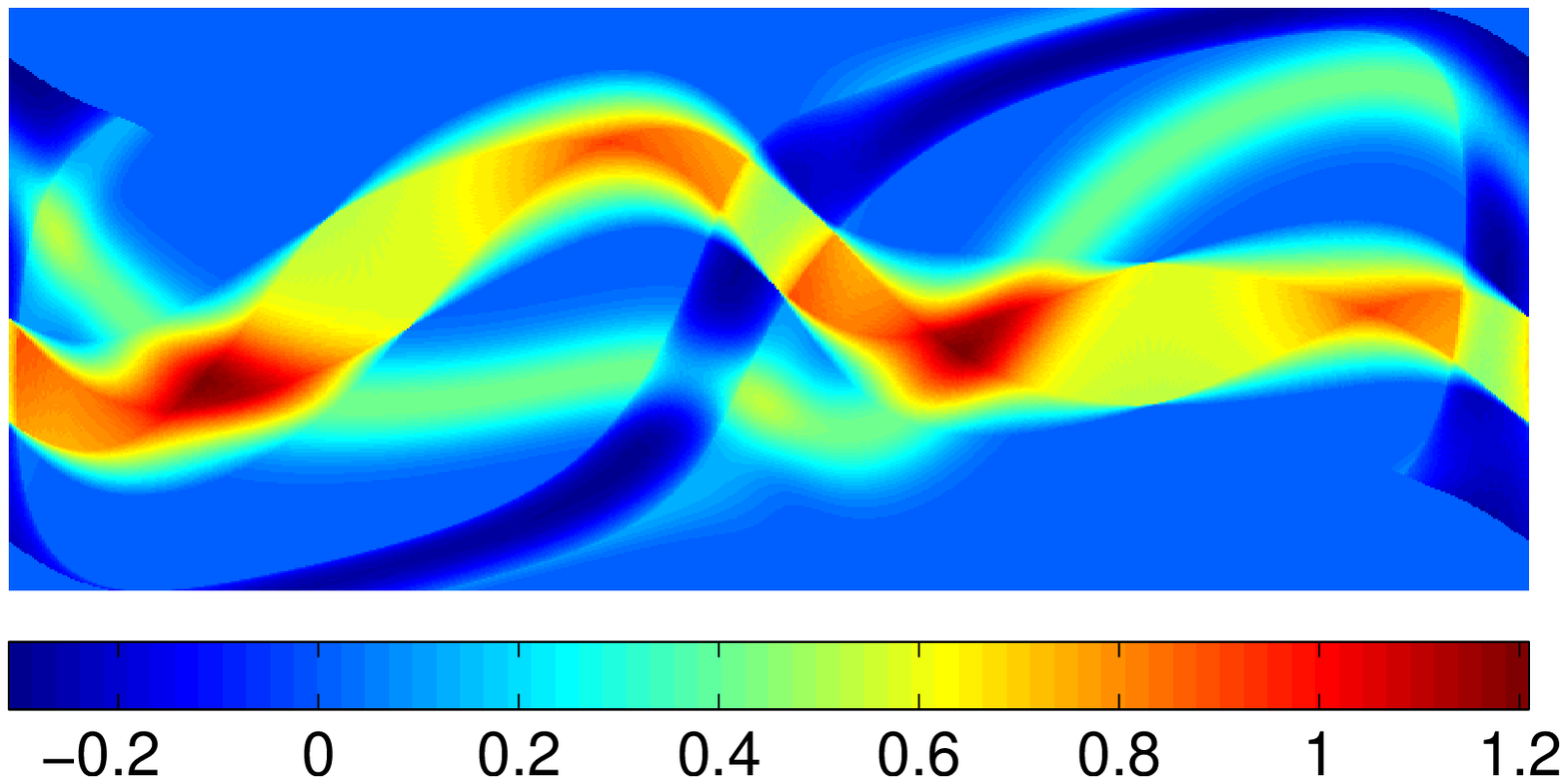}
	\label{fig:CPCns3}
    }
    \subfigure[Pointwise error $|f-f_{rc}|$]{
	\includegraphics[height=0.18\textheight]{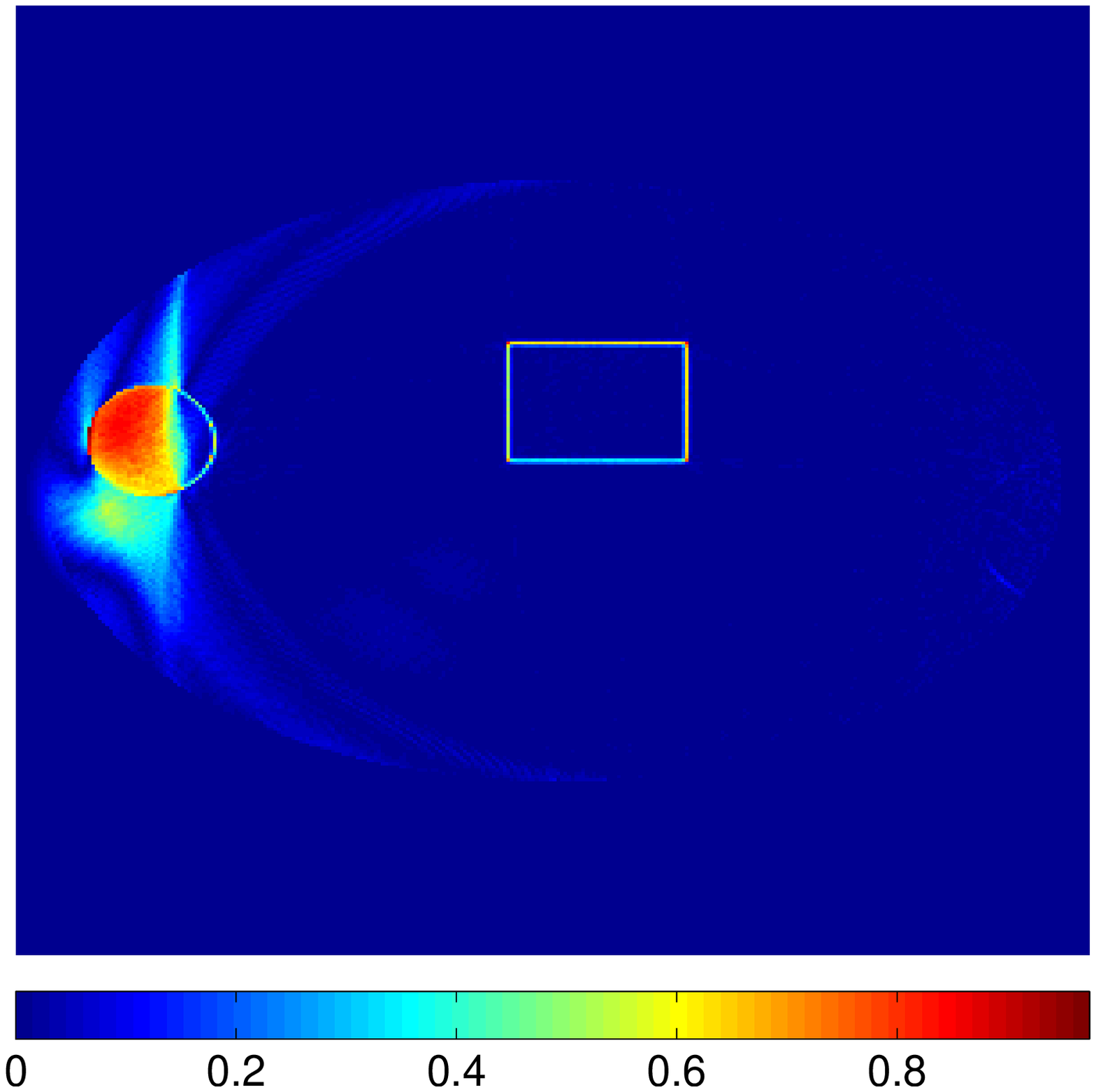}
	\label{fig:CPCns4}
    }
    \caption{Non-simple domain with constant positive curvature. }
    \label{fig:CPCnonsimple}
\end{figure}

\paragraph{Experiment with a non-simple domain} In this experiment, we now allow the domain to include conjugate points while remaining non-trapping. In order to do this, we consider the following case: pick the metric $g_{R,+}$ defined in \eqref{eq:poscurvg} with $R=1$, so that the circle $\{x^2 + y^2 = 1\}$ is the ``equator'', i.e. a trapped geodesic. Consider the elliptical domain $\left\{ \frac{x^2}{a^2} + \frac{y^2}{b^2} < 1  \right\}$ with $(a,b) = (1.2,0.8)$, so that the domain contains antipodal (i.e., in this case, conjugate) points, e.g. $(-1,0)$ and $(1,0)$, though does not contain an entire great circle, thus guaranteeing that the domain is not simple yet not trapping.
Results may be found Fig. \ref{fig:CPCnonsimple}, where the pointwise error clealy demonstrate that at points whose antipodal points are not included in the domain (e.g. the central part), the initial function is accurately reconstructed, while on the left part of the domain, where every point has a conjugate point inside the domain, the function is hardly reconstructed at all. This should be contrasted with the fact that equation \eqref{eq:frc} predicts one-shot and exact reconstruction of $f$ even in the non-simple case when curvature is constant. 

\subsection{Metrics with non-constant curvature} \label{sec:nonconstcurv}

We now move on to the case where the metric has non-constant curvature. Reconstruction formulas \eqref{eq:frc}-\eqref{eq:hrc} no longer allow for a one-shot inversion due to the presence of the error operators $W^2$ and $(W^\star)^2$, and one wishes to remove the Fredholm error from the reconstructions whenever possible. 

\subsubsection{Known results and iterative reconstruction algorithms}

A case leading again to direct inversion is when the metric has curvature close to constant. Indeed, it is established in \cite{Krishnan2010} a bound of the form
\begin{align*}
    \|W\|_{\L(L^2)} \le C \|\nabla \kappa\|_\infty, \quad (\kappa:\text{ Gaussian curvature}),
\end{align*}
for $C$ a fixed constant. This in turn ensures that if the curvature is small enough, the operators $W$ and $W^\star$ are contractions in $L^2(M)$, therefore the equations \eqref{eq:frc} and \eqref{eq:hrc} can be inverted for $f$ and $h$ via the Neumann series
\begin{align}
    f = - \sum_{k=0}^\infty (-W^2)^k (X_\perp w_\psi^{(f)})_0 \qandq h = - \sum_{k=0}^\infty (-(W^\star)^2)^k (w_\psi^{(h)})_0.
    \label{eq:neumann}
\end{align}

As one will see, implementing these Neumann series {\em does not require implementing the operators} $W^2$ and $(W^\star)^2$, as they can be directly expressed as 
\begin{align*}
    -W^2 f = f + (X_\perp w_\psi^{(f)})_0 \qandq -(W^\star)^2 h = h + (w^{(h)}_\psi)_0. 
\end{align*}
This principle has been seen in prior literature, see e.g. \cite{Bal2010a,Qian2011}: if one has an approximate reconstruction formula modulo a contractive error, the exact reconstruction can be deduced by setting up an iterative scheme, each step of which requires solving a forward problem and an approximate inversion. Put in formal terms, if the following equation holds for any $f$ in some Hilbert space
\begin{align*}
    f + \K f = AIf,
\end{align*}
with $I$ a ``forward operator", $A$ an ``approximate inverse" (parametrix) and $\K$ a contractive operator, then $f$ can be reconstructed via the following iterated sum
\begin{align}
    f = \sum_{k=0}^\infty (-\K)^k AIf = \sum_{k=0}^\infty (Id - AI)^k AIf. 
    \label{eq:formalNeumann}
\end{align}
An implementation of a partial sum of series \eqref{eq:formalNeumann} is described in Algorithm \ref{algo:neumann}.

\begin{algorithm}
    \caption{Computation of a partial sum of the Neumann series \eqref{eq:formalNeumann}.}
    \begin{algorithmic}
	\STATE Denote the data $\D = If$ 
	\STATE Compute $f = A\D$ 
	\STATE Set $g=f$
	\FOR{ $k=1$ to $N_{iter}$ }
	\STATE Compute and update: $g = g - AI g$ \COMMENT{at step $k$, $g = (Id - AI)^k A\D$}
	\STATE Update: $f = f + g$ \COMMENT{at step k, $f = \sum_{p=0}^k (Id - AI)^k A\D$}
	\ENDFOR
	\STATE Return $f$
    \end{algorithmic}
    \label{algo:neumann}
\end{algorithm}

Theoretical results predict that it makes sense to compute this series when the curvature is ``close enough'' to constant. However, how close to constant curvature the metric must be is not very well quantified, and numerics indicate good reconstructions even when implementing Algo. \ref{algo:neumann} with metrics close to non-simple. 

\subsubsection{Manifolds with a radially symmetric metric}

A first simplified case leading to an interesting toy model is that of isotropic and radially symmetric metrics, i.e. the scalar function $g$ only depends on the radial variable $r = \sqrt{x^2 + y^2}$, denoted $g = g(r) = c(r)^{-2}$ ($c(r)$ is refered to as the local speed of geodesics). This case was considered first by Herglotz in \cite{Herglotz1905}. See also the work in \cite{Sharafudtinov1997} where it is established injectivity over solenoidal tensors of any order and in any dimension $d \ge 2$, on spherically symmetric layers of the form $\{x\in \Rm^d : \rho_0 < |x|<\rho_1 \}$ for some $0<\rho_0<\rho_1$. 

When considering additionally that the domain is a centered disk of radius $R$, the geodesic flow has invariant properties under rotations about $0$. As it is useful for the study of the next toy model, let us first recall the Herglotz condition. Reparameterizing geodesics as $\gamma(t) = r(t) (\cos\alpha (t), \sin\alpha (t))$ and $\dot\gamma(t) = c(\gamma(t)) (\cos\theta(t), \sin\theta (t))$, direct calculations show that the geodesic equation is now the following system
\begin{align*}
    \dot{r} &= c(r) \cos (\theta-\alpha), \\
    \dot\alpha & = \frac{c(r)}{r} \sin(\theta-\alpha), \\
    \dot\theta & = c'(r) \sin(\theta-\alpha). 
\end{align*}
In particular, we have the ODE 
\begin{align*}
    \dot{\overline{\theta-\alpha}} = \left( c'(r) - \frac{c(r)}{r} \right) \sin(\theta-\alpha),
\end{align*}
from which we derive
\begin{align*}
    \dot{\overline{\cos(\theta-\alpha)}} = - r \frac{\partial}{\partial r} \left( \frac{c(r)}{r} \right) \sin^2(\theta-\alpha). 
\end{align*}
Looking at this system, we can make the following observations: if there exists $0<r_0<R$ such that $\frac{\partial}{\partial r} \left( \frac{c(r)}{r} \right) |_{r_0} = 0$, then the centered circle of radius $r_0$ is a trapped geodesic. Therefore, if we want the disk to be non-trapping, the function $\frac{\partial}{\partial r} \left( \frac{c(r)}{r} \right)$ cannot vanish and must therefore have constant sign. Now a geodesic will go toward the center if the quantity $\cos(\theta-\alpha)$ is decreasing. In order to avoid this case, one must ensure that 
\begin{align*}
    \frac{\partial}{\partial r} \left( \frac{c(r)}{r} \right) <0 \quad \Leftrightarrow \quad \frac{\partial}{\partial r} \left( \frac{r}{c(r)} \right) >0, 
\end{align*}
a condition first found by Herglotz in \cite{Herglotz1905}.

\paragraph{A toy model for focusing lens}

Let us consider $M$ to be the unit disk $\{x^2 + y^2 \le 1\}$, and consider the family of isotropic and radially symmetric metrics with parameter $k\in (0,\infty)$
\begin{align*}
    g_k(r) = \exp \left( k\exp \left( -\frac{r^2}{2\sigma^2} \right)\right), \quad 0\le r\le 1,  \quad \sigma\text{ fixed}
\end{align*}
with corresponding velocity $c_k(r) = g_k^{-\frac{1}{2}}(r) = \exp\left( -\frac{k}{2} \exp\left( -\frac{r^2}{2\sigma^2} \right) \right)$. Considering the Herglotz condition for a non-trapping metric, we compute
\begin{align*}
    \frac{\partial}{\partial r} (r c_k^{-1}(r)) = \exp\left( \frac{k}{2} \exp\left( -\frac{r^2}{2\sigma^2} \right) \right) \left( 1 - k f\left( \frac{r^2}{2\sigma^2} \right) \right),
\end{align*}
where we have defined $f(x) := xe^{-x}$. The function $f$ satisfies $0\le f(x) \le e^{-1}$ for $0\le x\le \infty$ and reaches its maximum at $x=1$ where $f(1) = e^{-1}$. This means that for $k<e$, the Herglotz condition is satisfied and the manifold is non-trapping. For $k\ge e$, the manifold becomes trapping. For instance, at $k=e$, the circle $r = \sqrt{2}\sigma$ is a trapped geodesic (pick $\sigma$ small enough so that this circle lies inside the unit disk). 

We shall draw the heuristic conclusions:
\begin{romannum}
    \item As $k=0$, the metric is flat and thus $\beta_{Ter} = \infty$.
    \item As $k$ increases, the metric goes from simple to non-simple non-trapping. One would expect that in that range, $\beta_{Ter}$ decreases from $\infty$ to $1$ and even less whenever simplicity no longer holds. Numerics indicate that the transition to non-simple occurs approximately at $k\approx 0.47$ when the domain is the unit disk and $\sigma = 0.25$.
    \item As $k$ reaches $e$, the manifold becomes trapping. 
\end{romannum}

\begin{figure}[htpb]
    \centering 
    \subfigure[$k=0.20$ (simple)]{
	\includegraphics[width=.22\textwidth]{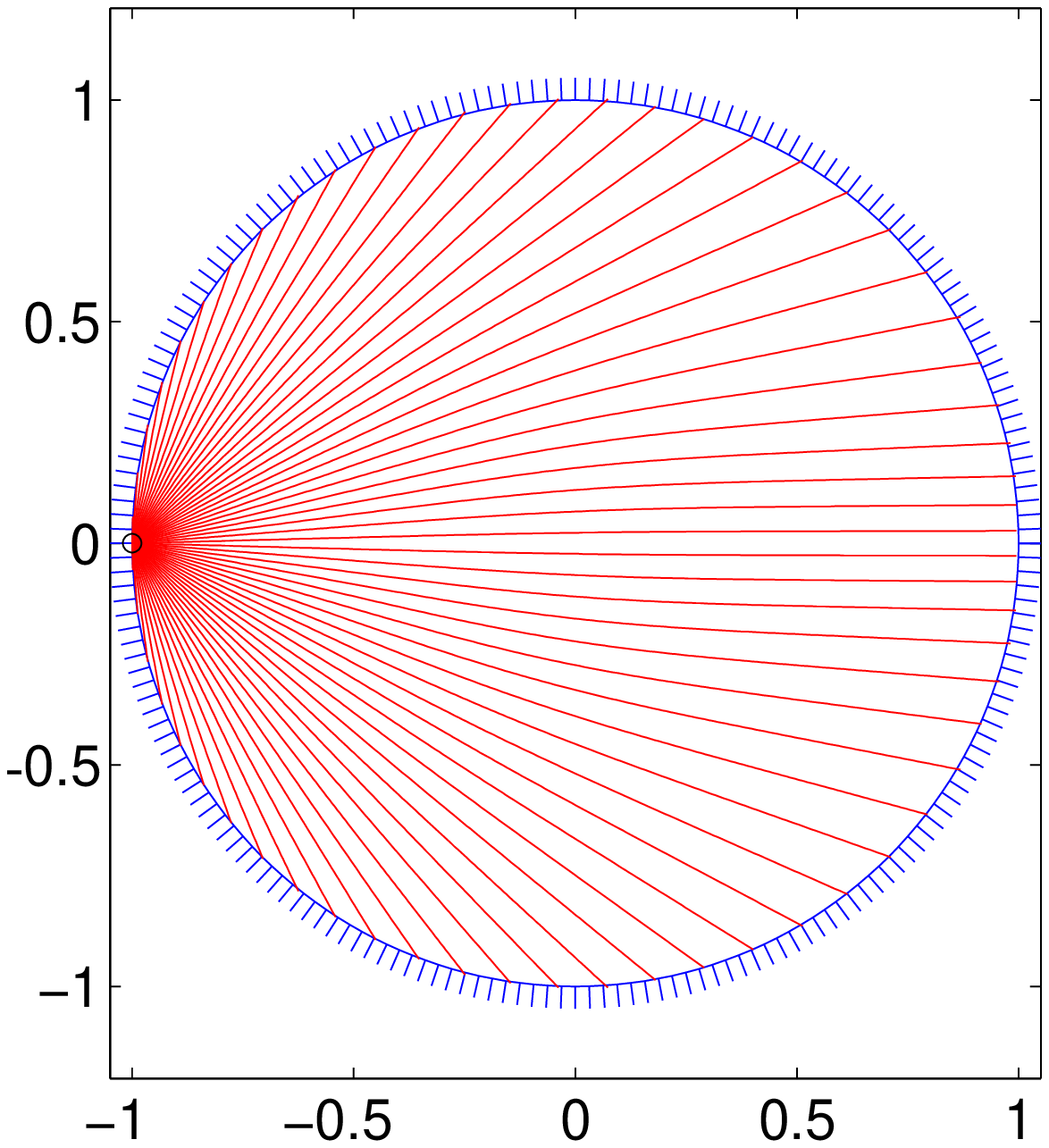}
	\label{fig:1lens1}
    }
    \subfigure[$k=0.49$ (non-simple)]{
	\includegraphics[width=.22\textwidth]{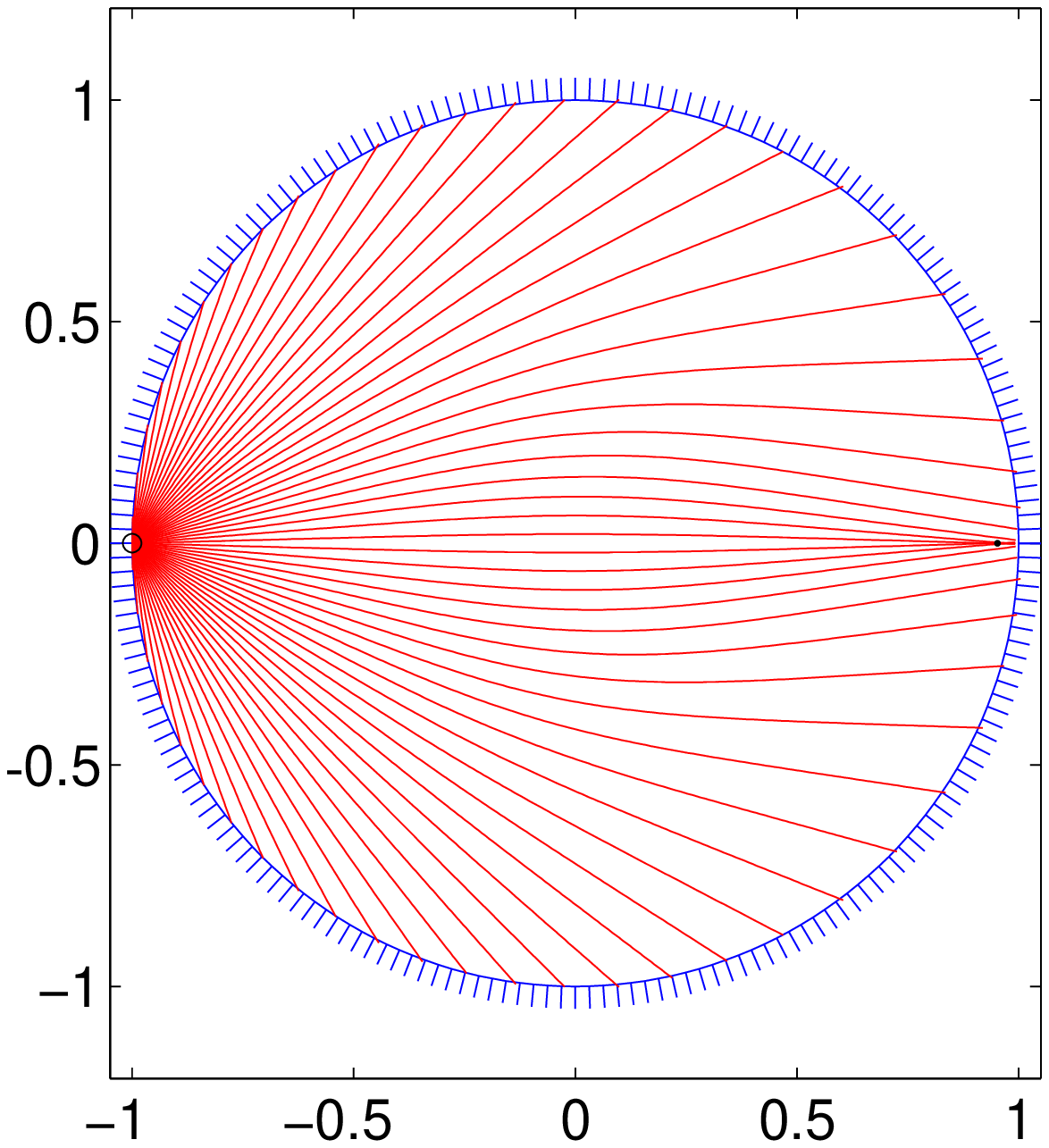}
	\label{fig:1lens2}
    }
    \subfigure[$k=1.23$ (non-simple)]{
	\includegraphics[width=.22\textwidth]{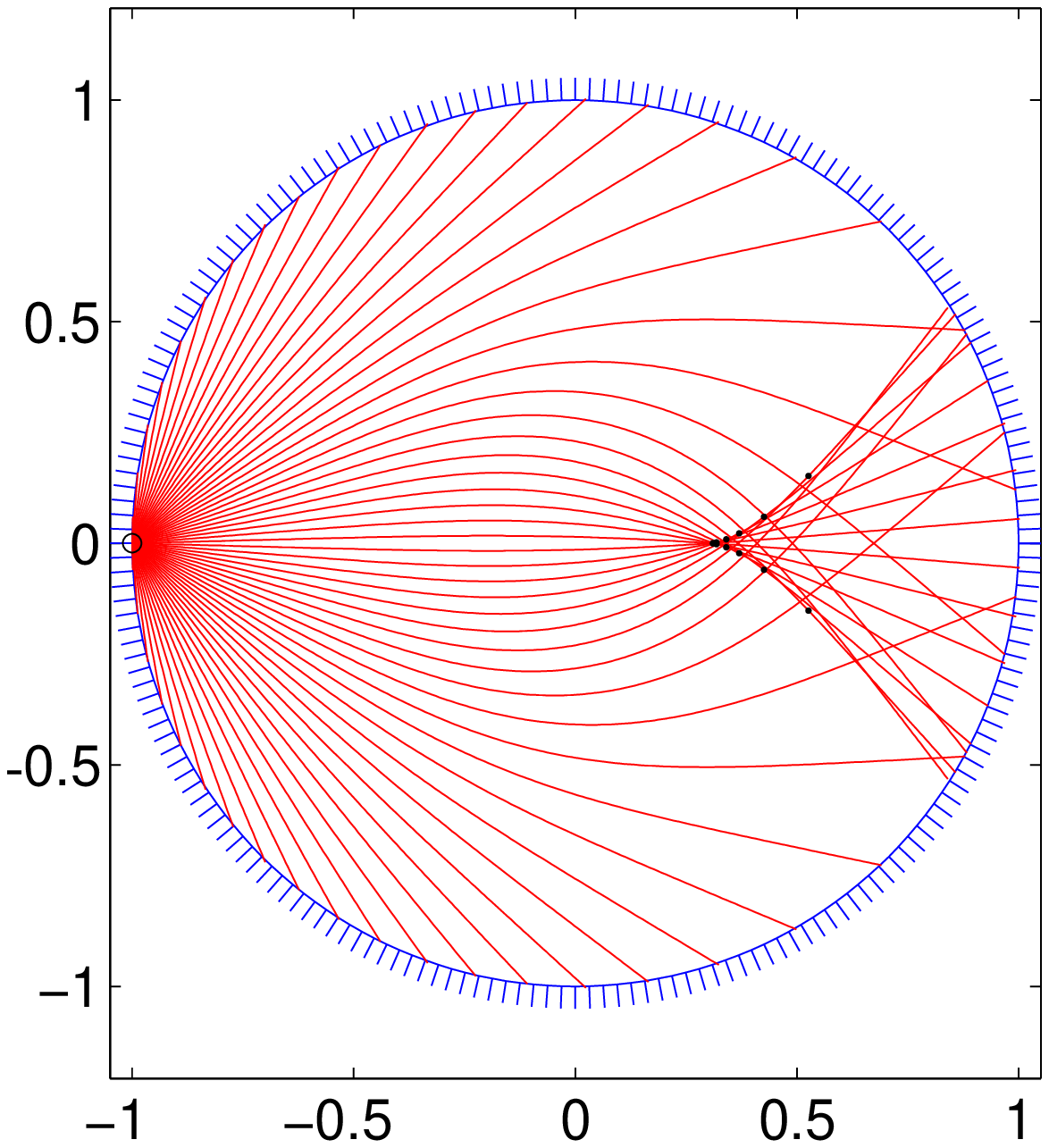}
	\label{fig:1lens3}
    }
    \subfigure[plot of $\beta_{Ter}$ versus $k$]{
	\includegraphics[width=.22\textwidth]{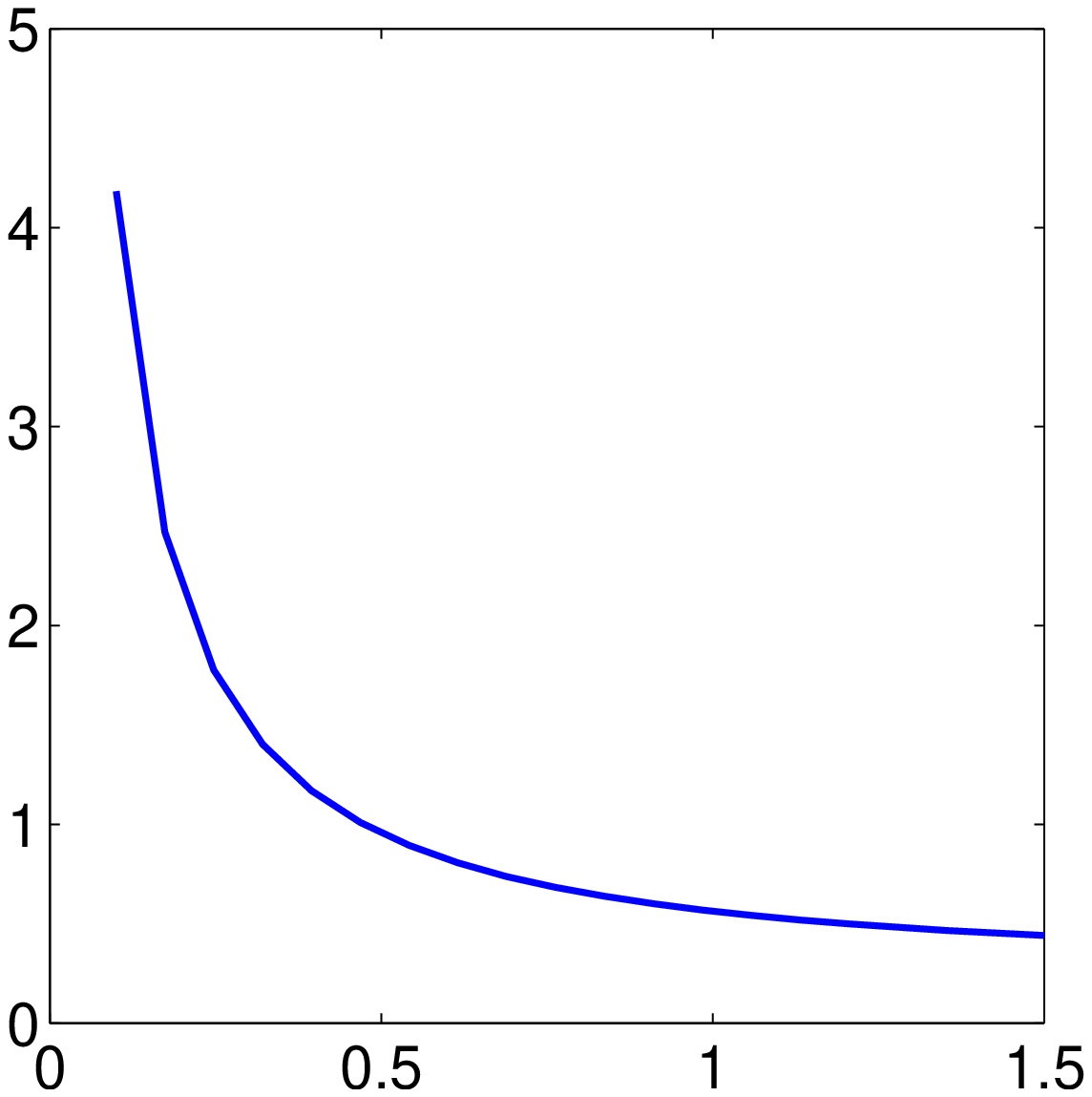}
	\label{fig:1betater}
    }
    \caption{Geodesics emanating from $(-1,0)$ for various values of $k$ for some metrics of lens type. Black dots indicate conjugate points. The conjugate locus consists of a cusp from which emanate two fold branches. \subref{fig:1betater}: plot of $\beta_{Ter}$ versus $k$ within $10^{-3}$ computed using Algo. \ref{algo:beta}.}
    \label{fig:1}
\end{figure}

The computation of $\beta_{Ter}$ may be obtained by dichotomy as described in Algorithm \ref{algo:beta} and Fig. \ref{fig:1betater} displays a computed plot of $\beta_{Ter}$ versus $k$ using this algorithm.

\begin{algorithm}
    \caption{Computation of $\beta_{Ter}$ by dichotomy.}
    \begin{algorithmic}
	\STATE Set threshold $\varepsilon$, pick initial $\beta_m << 1$ and $\beta_M >>1$. 
	\WHILE{ $\beta_M - \beta_n\ge \varepsilon$ }
	\STATE Set $\beta = \frac{1}{2} (\beta_M + \beta_m)$.
	\STATE Find the $\beta$-conjugate points of $M$ by extracting the points where $b_\beta$ vanishes. 
	\IF{$M$ is free of $\beta$-conjugate points}
	\STATE $\beta_m = \beta$ 
	\ELSE 
	\STATE $\beta_M = \beta$
	\ENDIF
	\ENDWHILE
	\STATE Return $\frac{1}{2} (\beta_M + \beta_m)$.	
    \end{algorithmic}
    \label{algo:beta}
\end{algorithm}

\subsubsection{Examples} \label{sec:examples}

\begin{figure}[htpb]
    \centering 
    \subfigure[$k=0.3$ (simple)]{
    \includegraphics[width=0.25\textwidth]{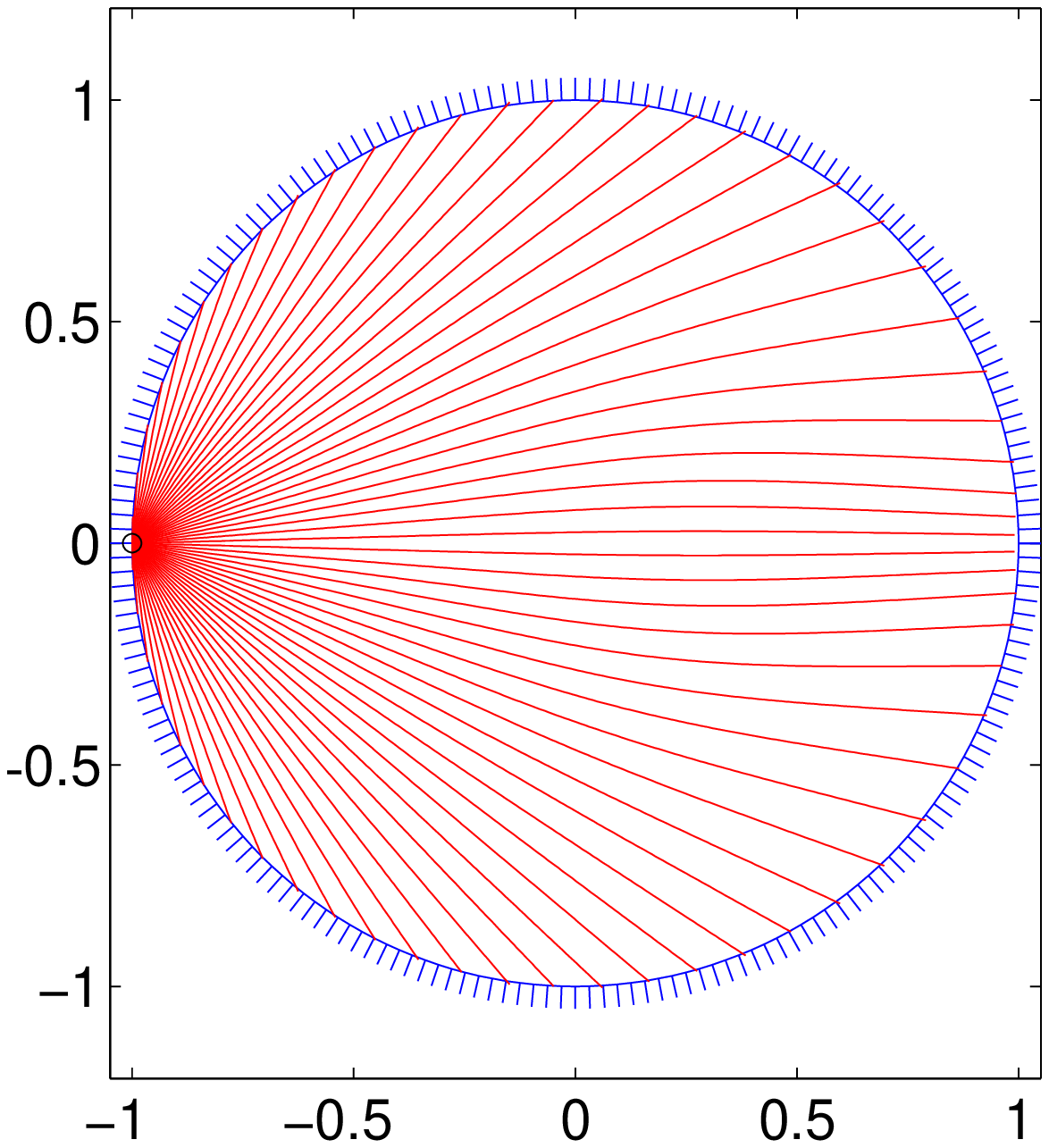}
    \label{fig:lensgeo1}
    }
    \subfigure[$k=0.6$ (non-simple)]{
    \includegraphics[width=0.25\textwidth]{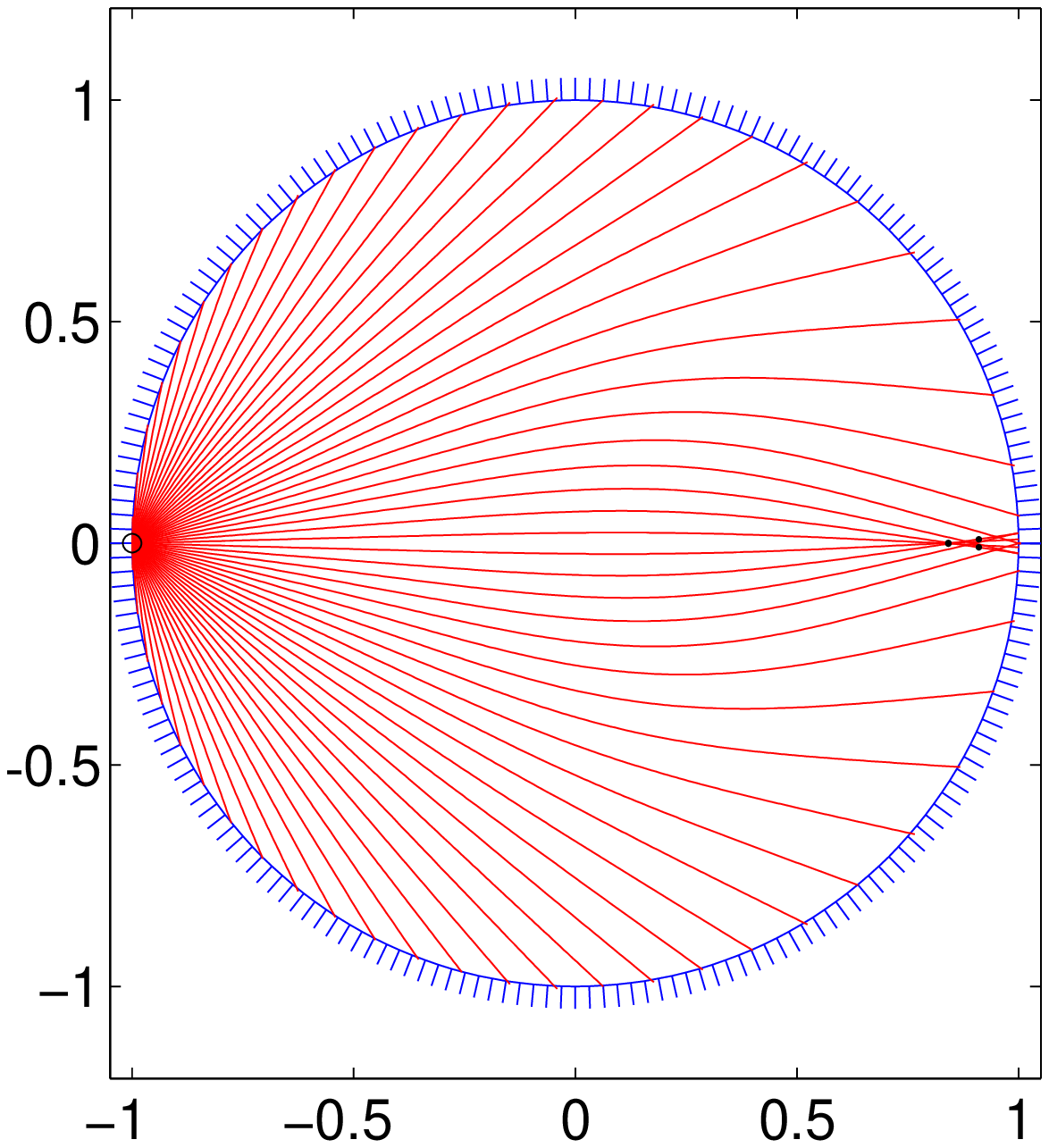}
    \label{fig:lensgeo2}
    }
    \subfigure[$k=1.2$ (non-simple)]{
    \includegraphics[width=0.25\textwidth]{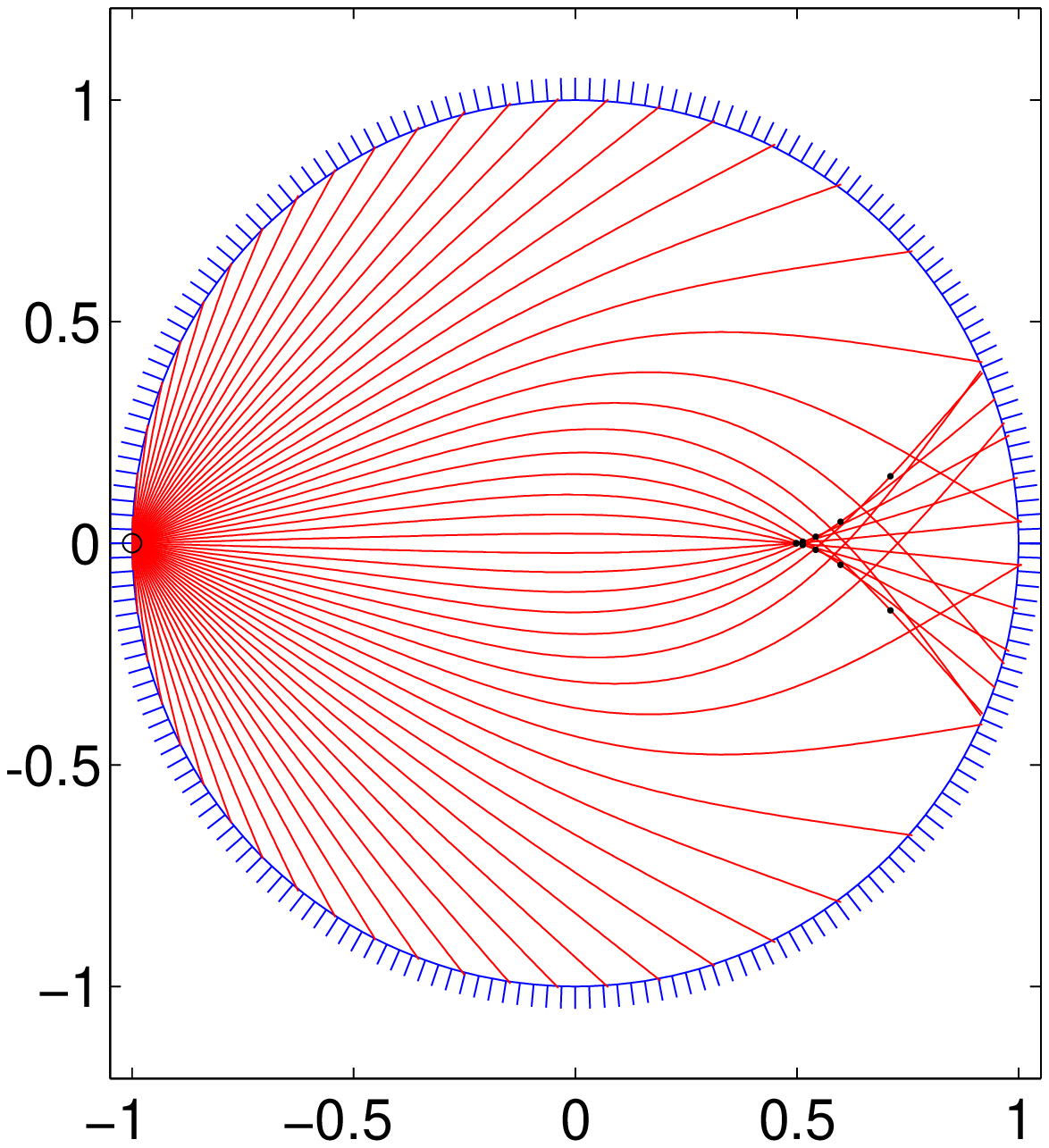}
    \label{fig:lensgeo3}
    }
    \caption{Some geodesics for the metric of lens type, described by equation \eqref{eq:lens}.}
    \label{fig:lens_init}
\end{figure}

We now present numerical implementations of Algo. \ref{algo:neumann} for both reconstruction formulas \eqref{eq:neumann}, using resolution $n=250$. The domain is the unit disk, with both smooth (Fig. \ref{fig:phantom4}) and non-smooth (Fig. \ref{fig:phantom1}) phantoms, and the metrics are of the same type as before (``focusing lens'') shifted away from the center to break symmetry, of scalar expression
\begin{align}   
    g(x,y) = \exp \left( -k  \exp \left( - (( x-0.2 )^2 + y^2 )/ (2\sigma^2) \right) \right),
    \label{eq:lens}
\end{align}
with $\sigma=0.25$ and lens parameter $k$ taking values in $\{0.3, 0.6, 1.2\}$. The case $k=0.3$ is simple while the remaining two are not, though the case $k=0.6$ is ``closer to simple'' than $k=1.2$ in the sense that $\beta_{Ter}(1.2)<\beta_{Ter}(0.6)<1$. Some geodesics are being displayed Fig. \ref{fig:lens_init}.

We perform the following simulations: \\
\noindent{\bf Experiment 1:} Smooth phantom (Fig. \ref{fig:phantom4}), iterative inversion from data $I_0 f$. \\
\noindent{\bf Experiment 2:} Smooth phantom (Fig. \ref{fig:phantom4}), iterative inversion from data $I_1 X_\perp f$. \\
\noindent{\bf Experiment 3:} Non-smooth phantom (Fig. \ref{fig:phantom1}), iterative inversion from data $I_1 X_\perp f$. 

\begin{figure}[htpb]
    \centering
    \subfigure[Metric defined in \eqref{eq:lens} with $k=0.3$ (simple)]{
	\includegraphics[height=0.15\textheight]{fig43.eps}    
	\includegraphics[height=0.16\textheight]{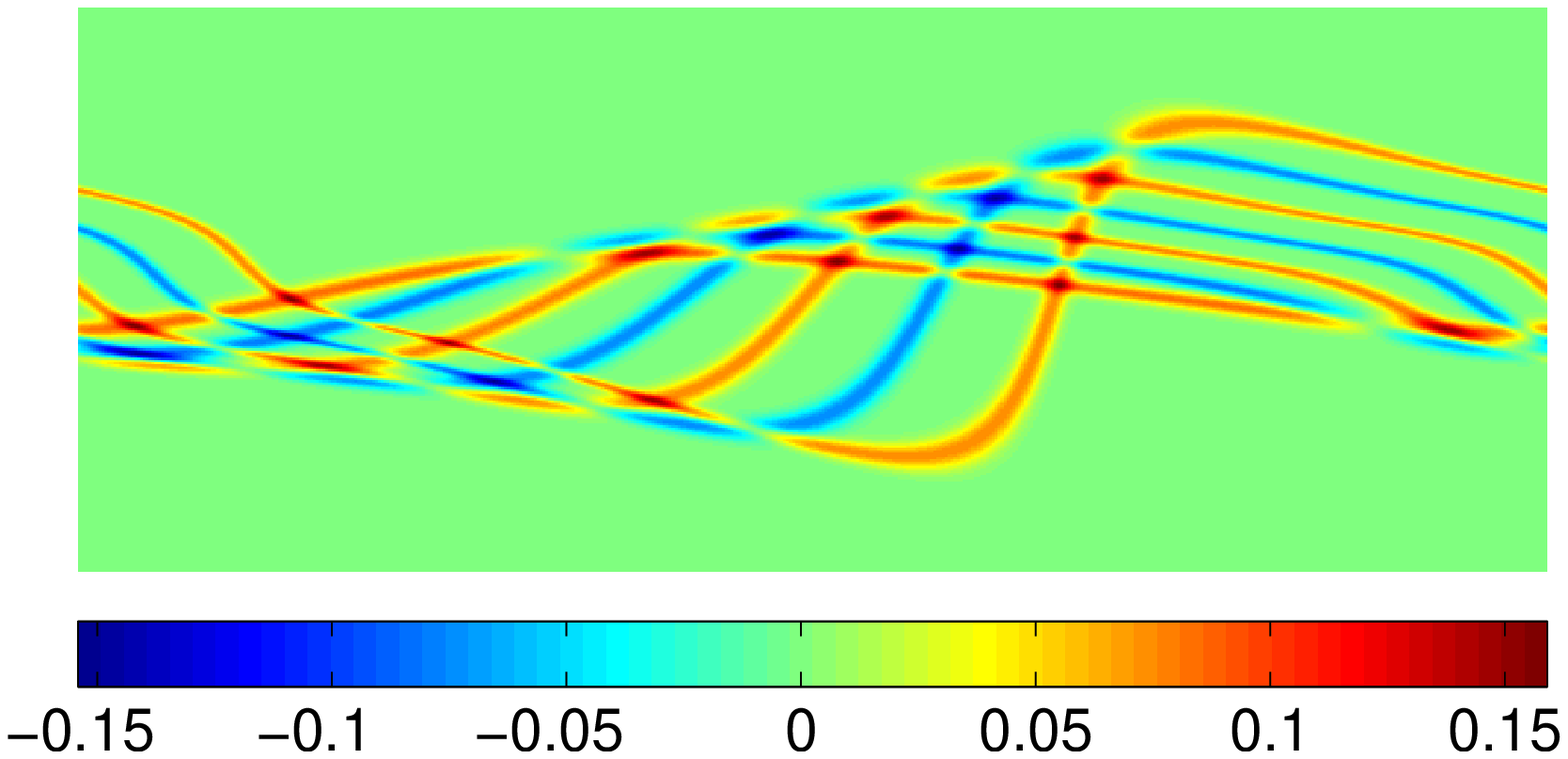}
	\includegraphics[height=0.16\textheight]{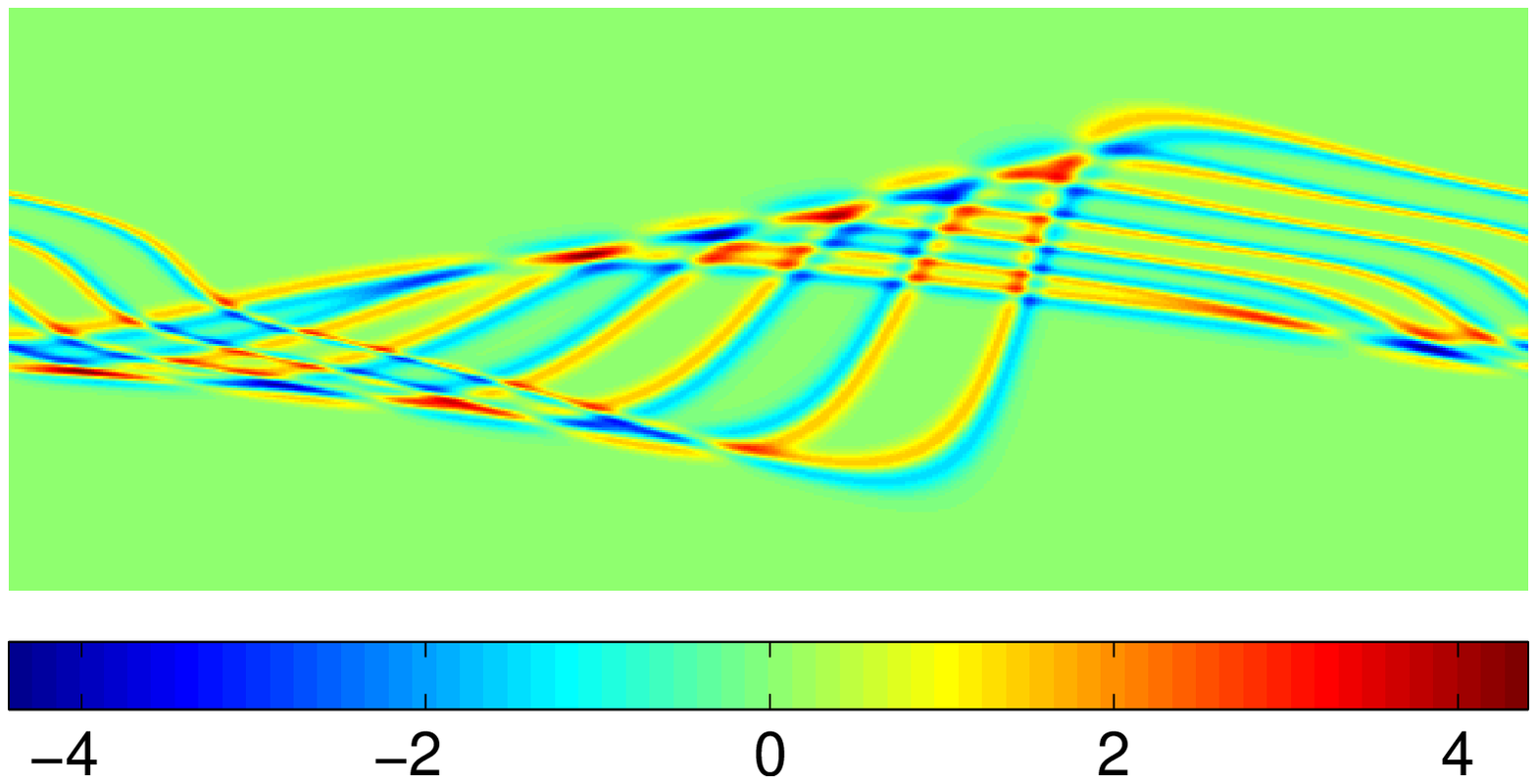}  
	\label{fig:lens1}
    }
    \subfigure[Metric defined in \eqref{eq:lens} with $k=0.6$ (not simple)]{
	\includegraphics[height=0.15\textheight]{fig44.eps}    
	\includegraphics[height=0.16\textheight]{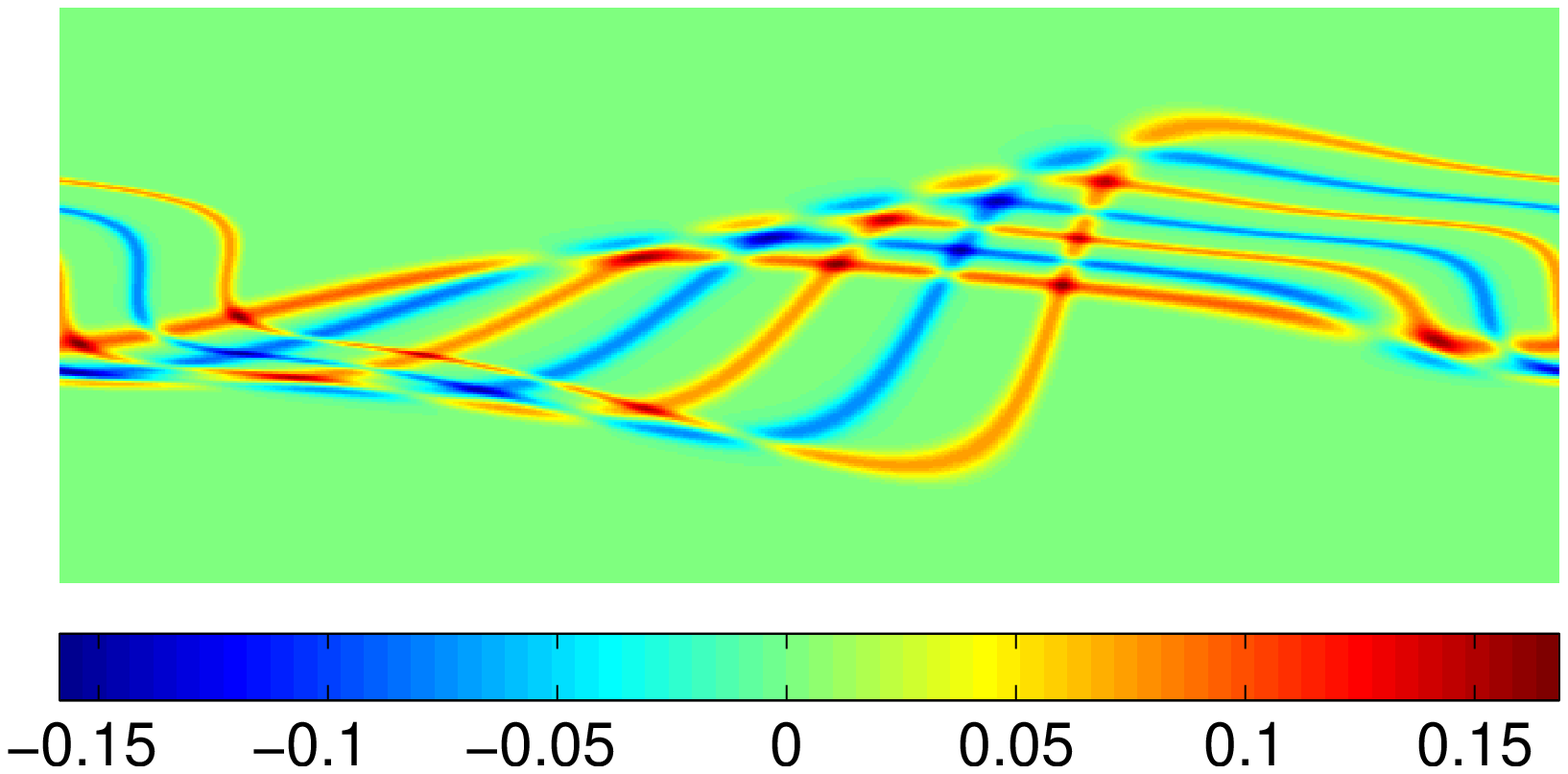}
	\includegraphics[height=0.16\textheight]{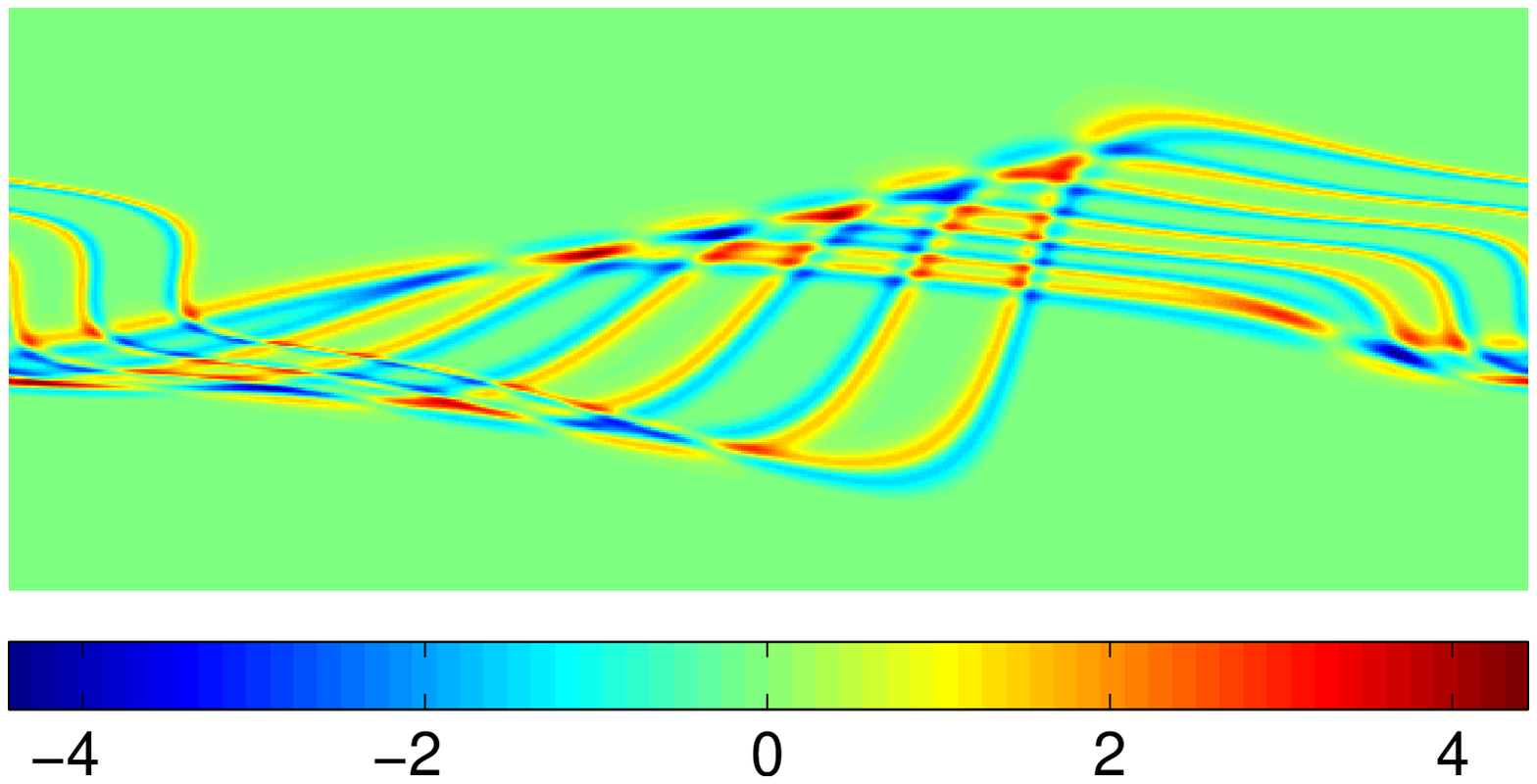}
	\label{fig:lens2}
    }
    \subfigure[Metric defined in \eqref{eq:lens} with $k=1.2$ (not simple)]{
	\includegraphics[height=0.15\textheight]{fig45.eps}    
	\includegraphics[height=0.16\textheight]{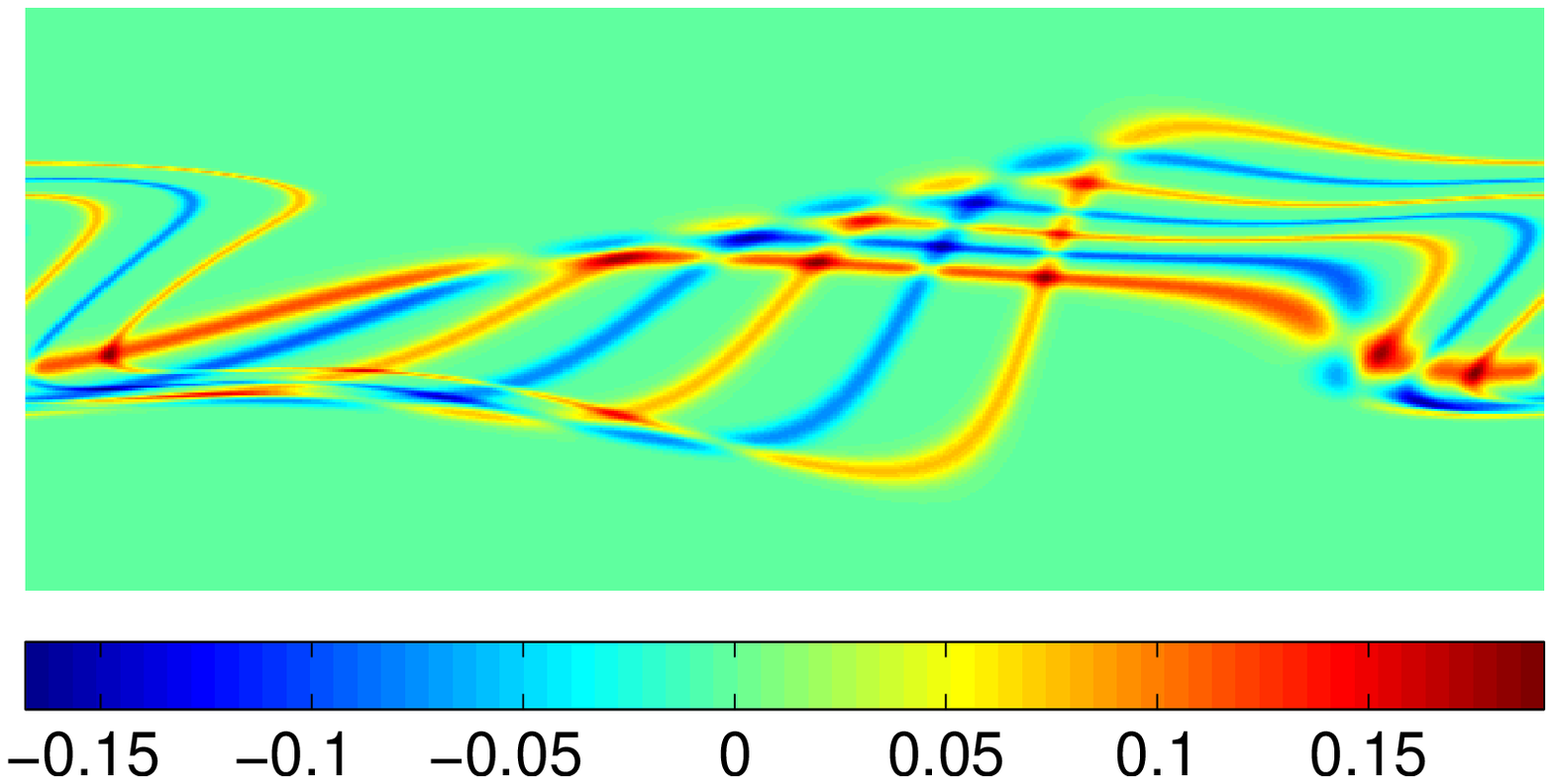}
	\includegraphics[height=0.16\textheight]{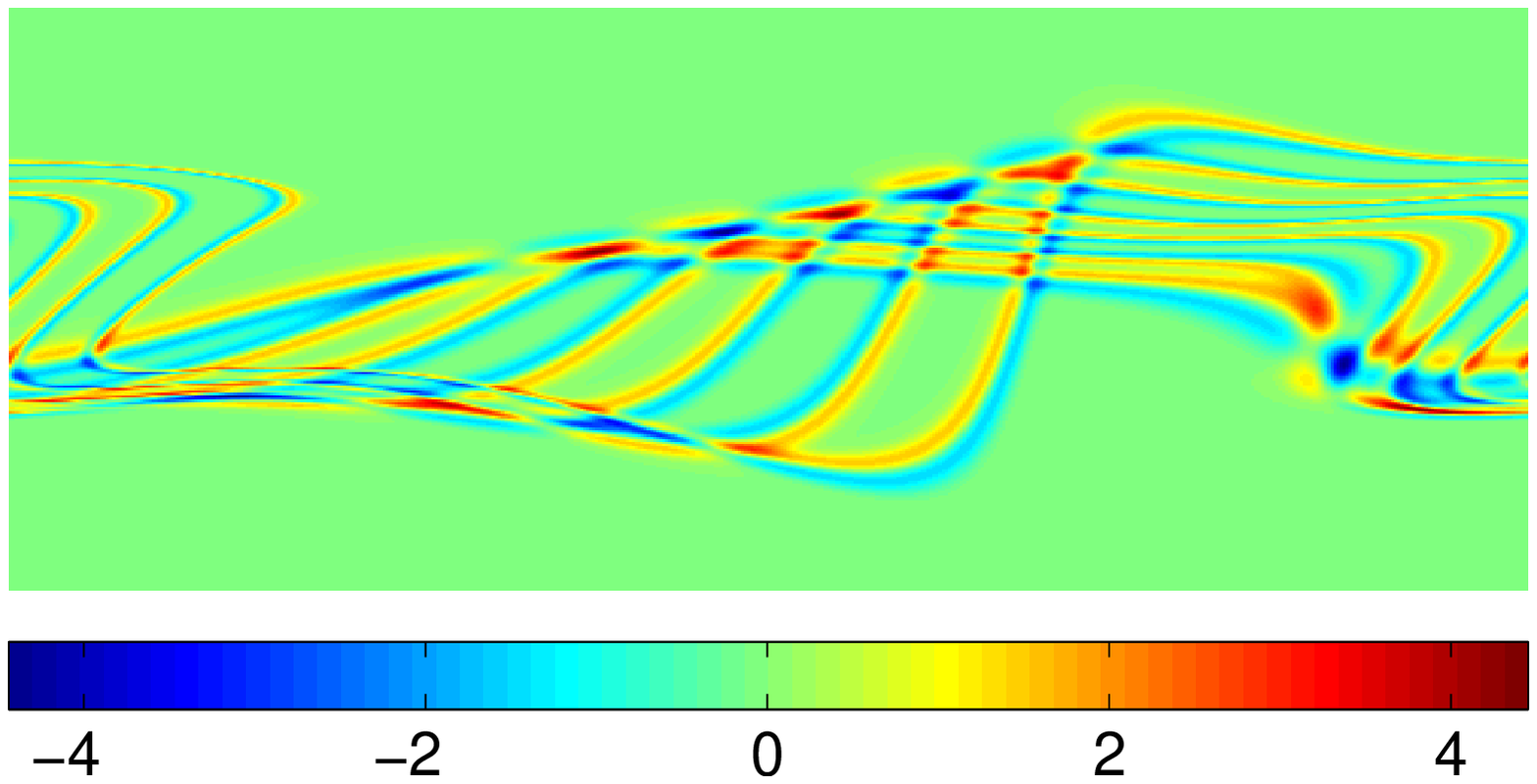}    
	\label{fig:lens4}
    }
    \caption{Experiments 1 and 2. Left to right: some geodesics for each metric used, forward data $I_0 f$ (Experiment 1) and $I_1(X_\perp f)$ (Experiment 2). }
    \label{fig:lensData}
\end{figure}

\begin{figure}[htpb]
    \centering
    \subfigure[$k=0.3$, iteration 1]{
	\includegraphics[height=0.22\textheight]{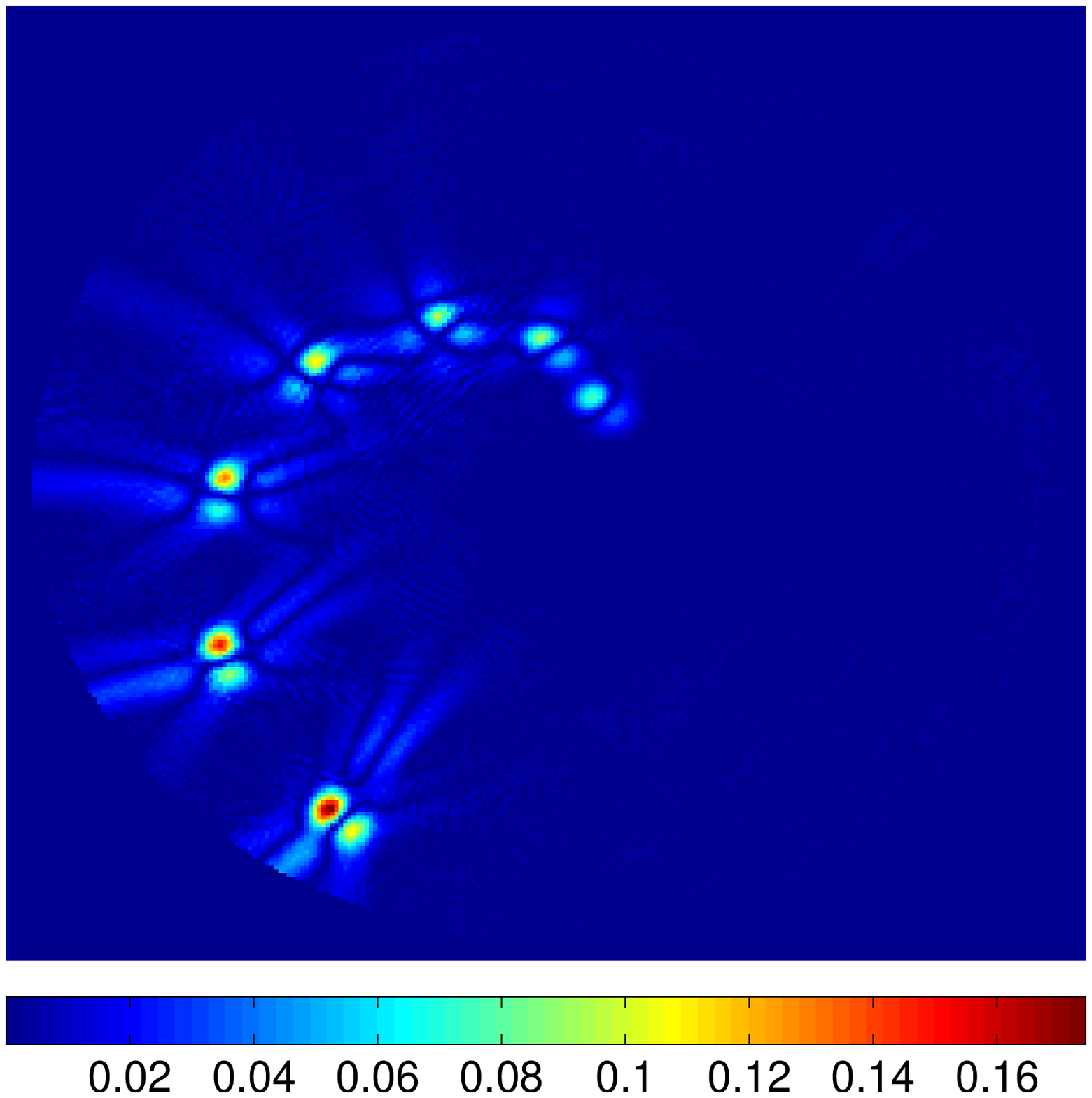}
	\label{fig:sample1}
    }
    \subfigure[$k=0.3$, iteration 9]{
	\includegraphics[height=0.22\textheight]{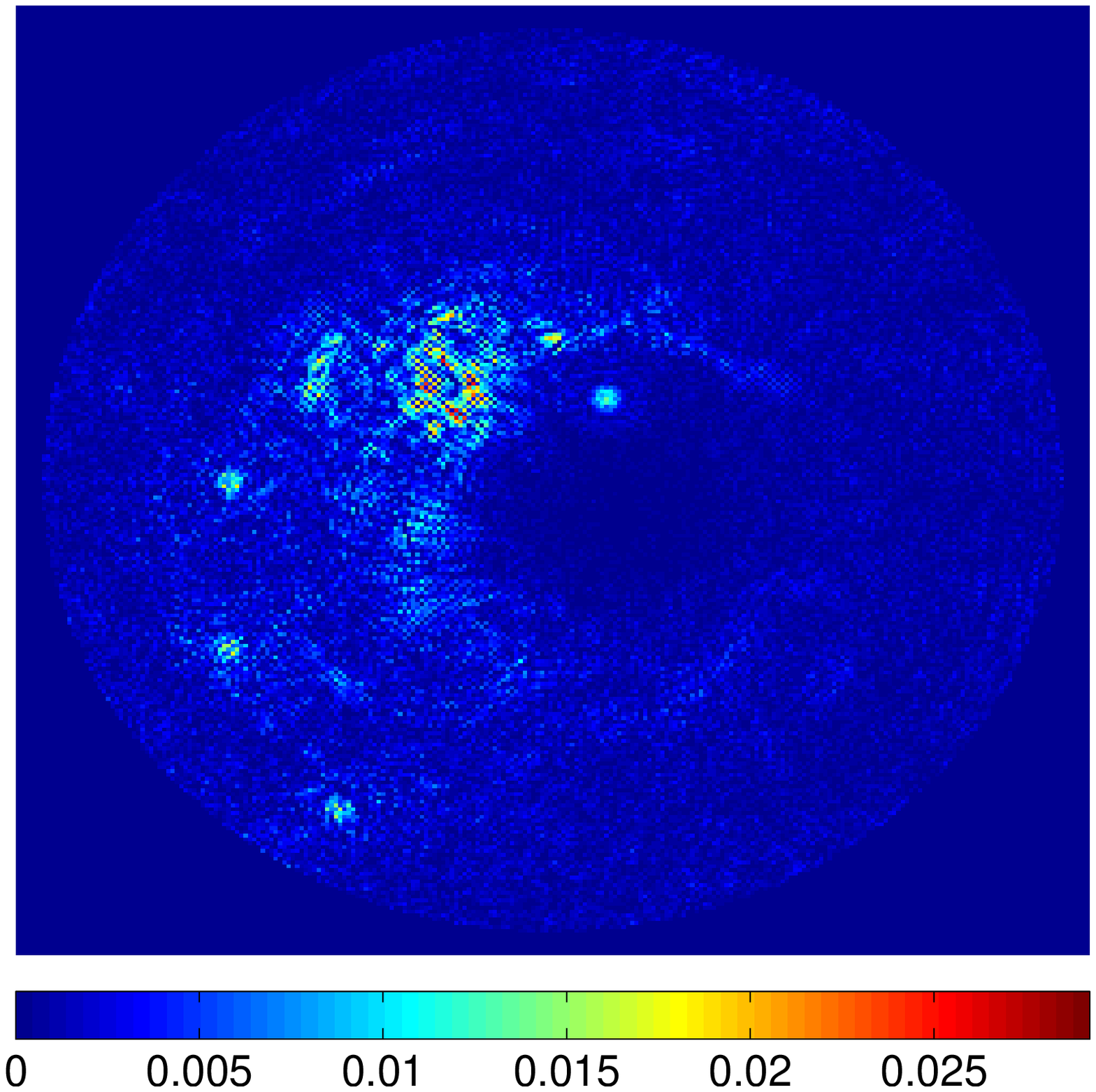} 
	\label{fig:sample2}	
    }
    \subfigure[$k=1.2$, iteration 9]{
	\includegraphics[height=0.22\textheight]{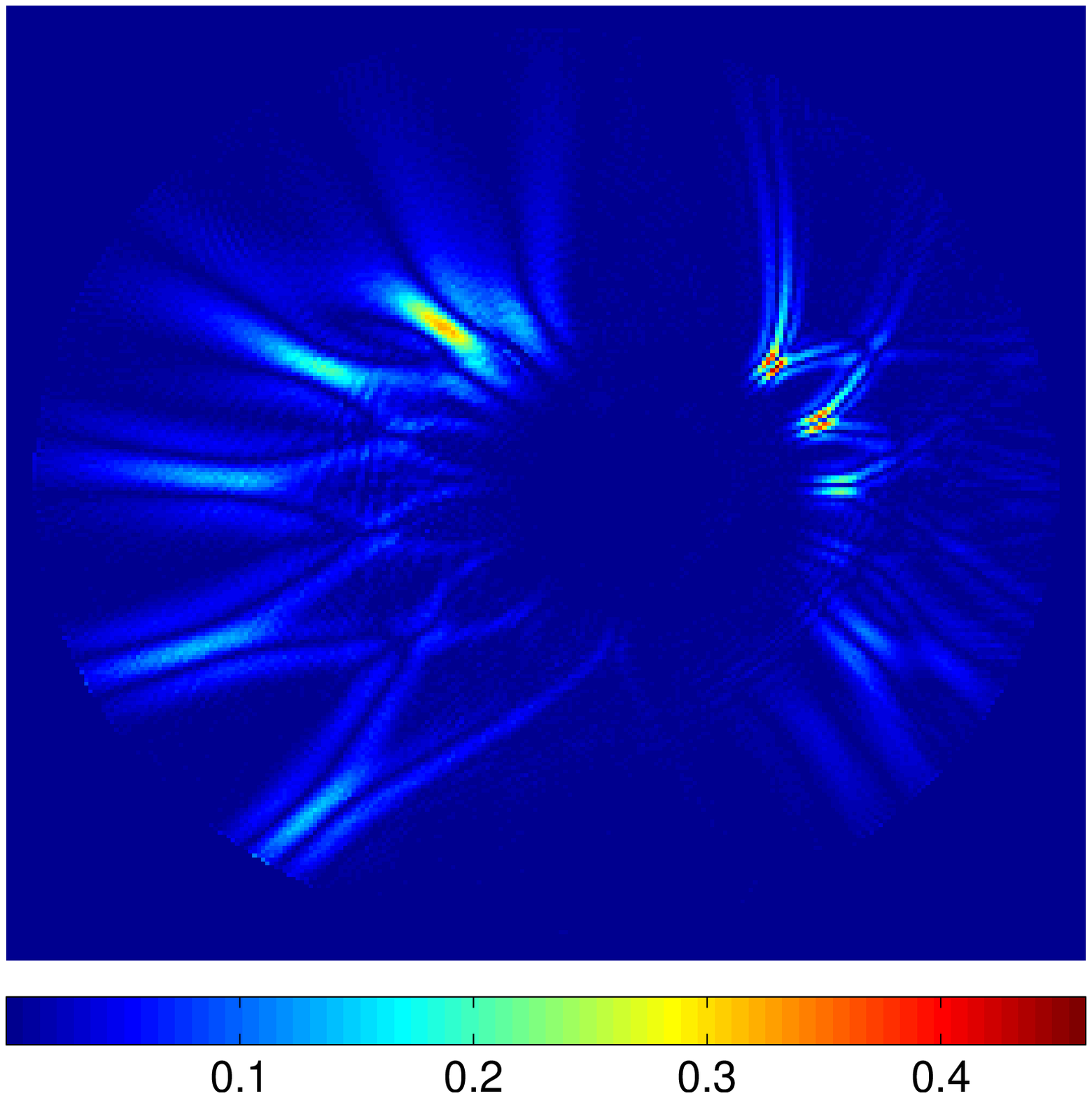}
	\label{fig:sample3}
    }

    \caption{Sample pointwise errors on $f$ for Experiment 1 (results for Experiment 2 are qualitatively similar). \subref{fig:sample1} and \subref{fig:sample2} show a case where the series converges. \subref{fig:sample3} displays the type of persistent artifact obtained in a non-simple case ($k=1.2$) where the series does not converge. Artifacts appear at the conjugate locus of each initial gaussian.}
    \label{fig:pwerror}
\end{figure}

\begin{figure}[htpb]
    \centering
    \includegraphics[height=0.17\textheight]{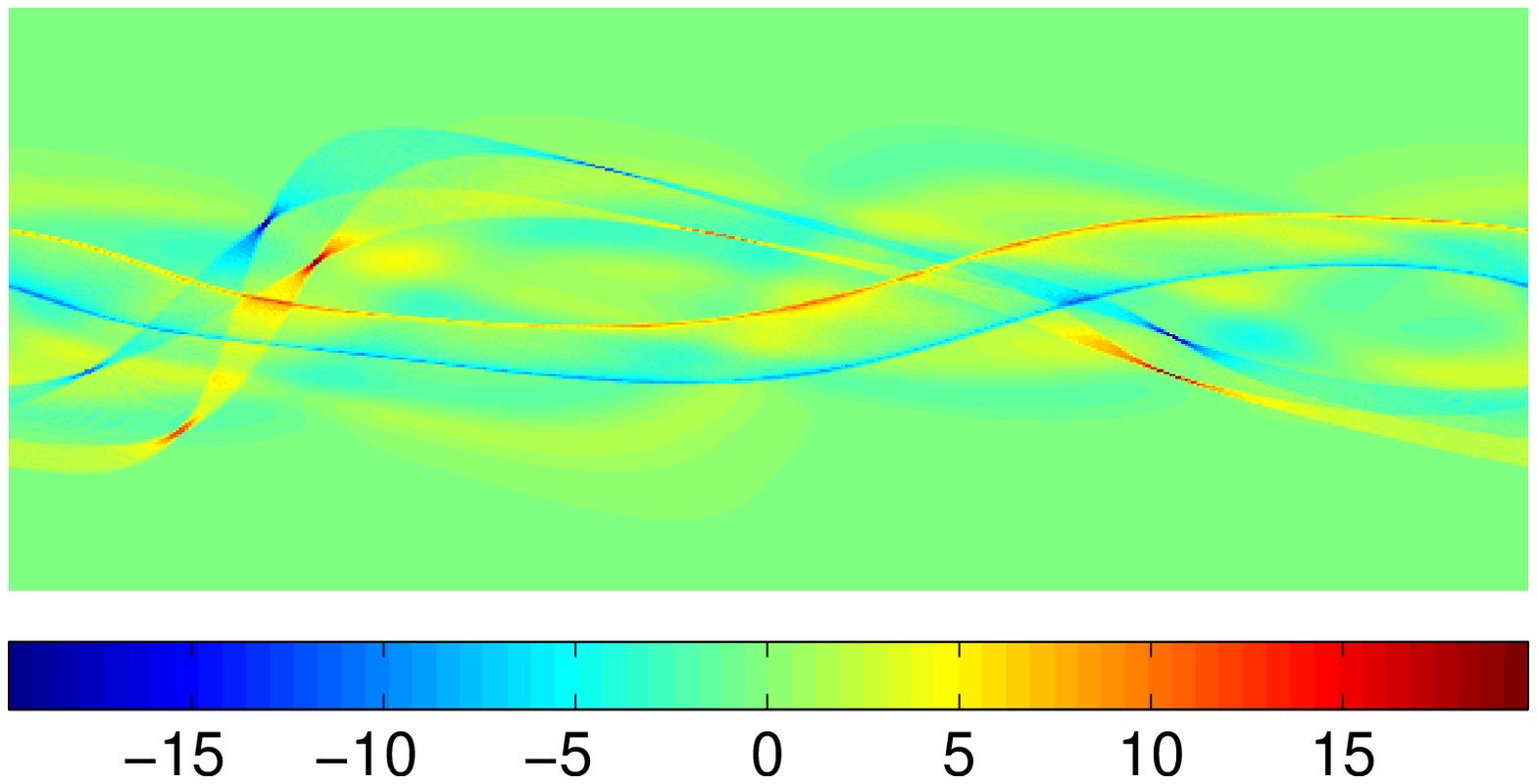}
    \includegraphics[height=0.17\textheight]{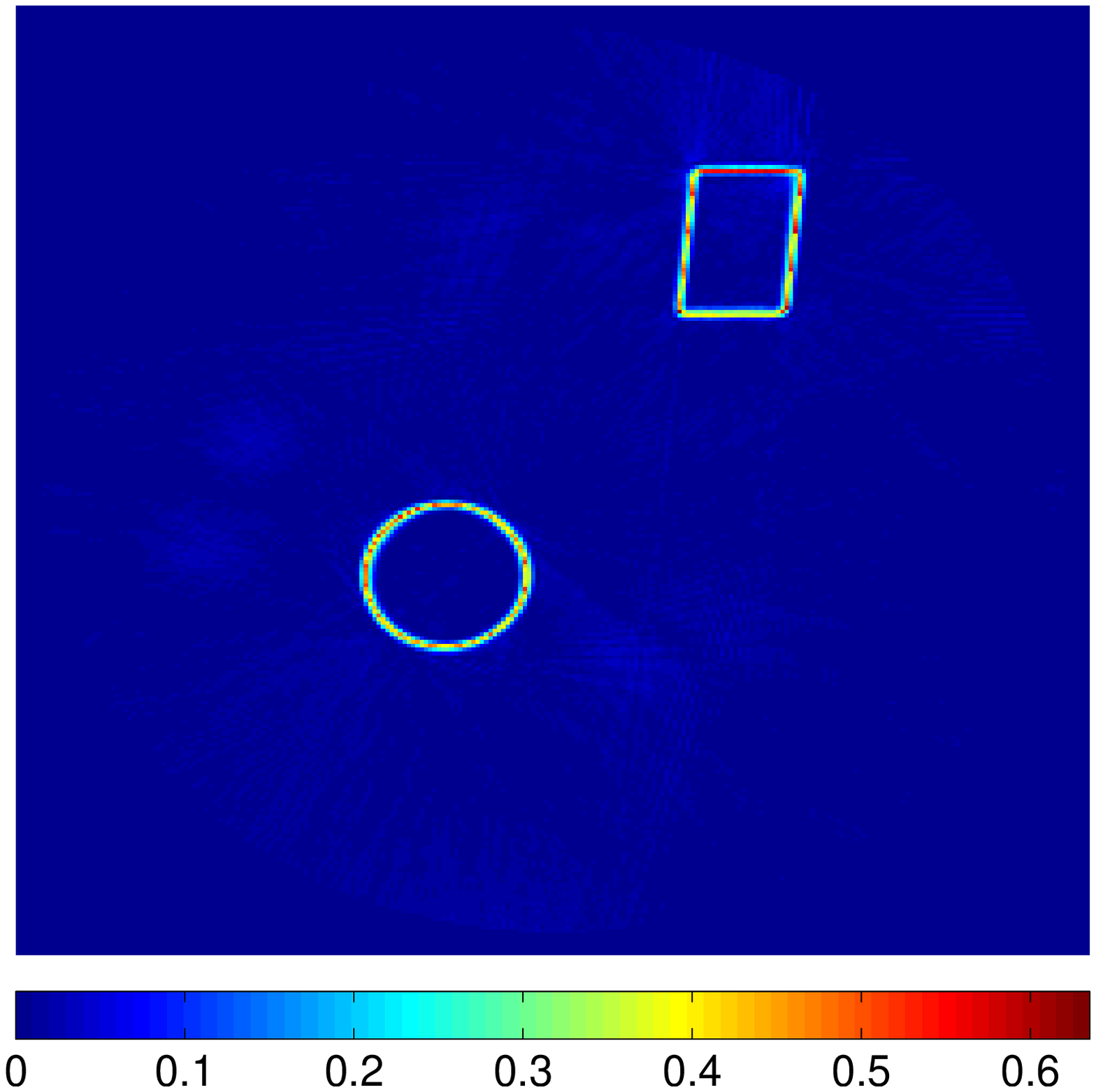}
    \includegraphics[height=0.17\textheight]{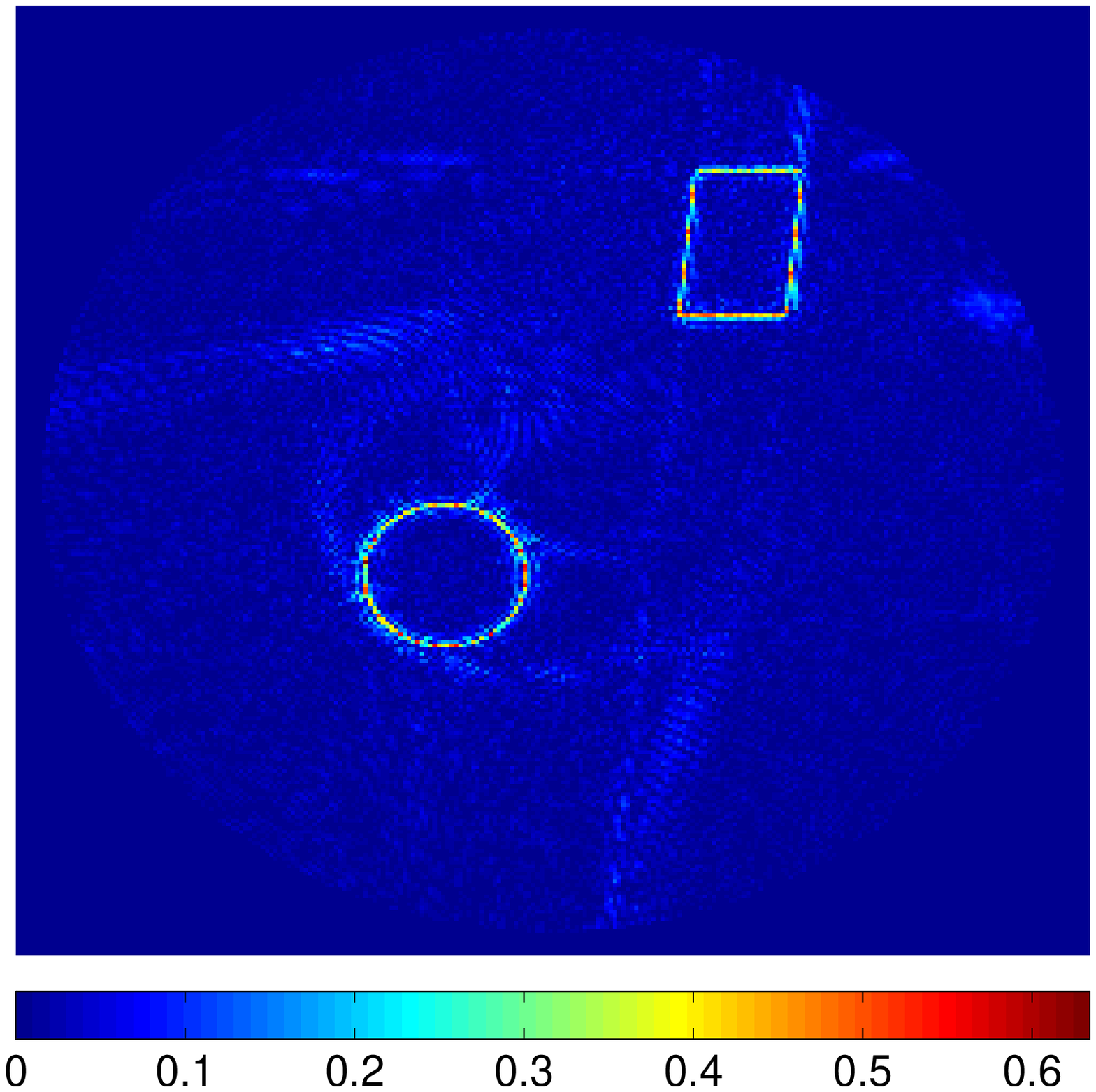} \\
    \includegraphics[height=0.17\textheight]{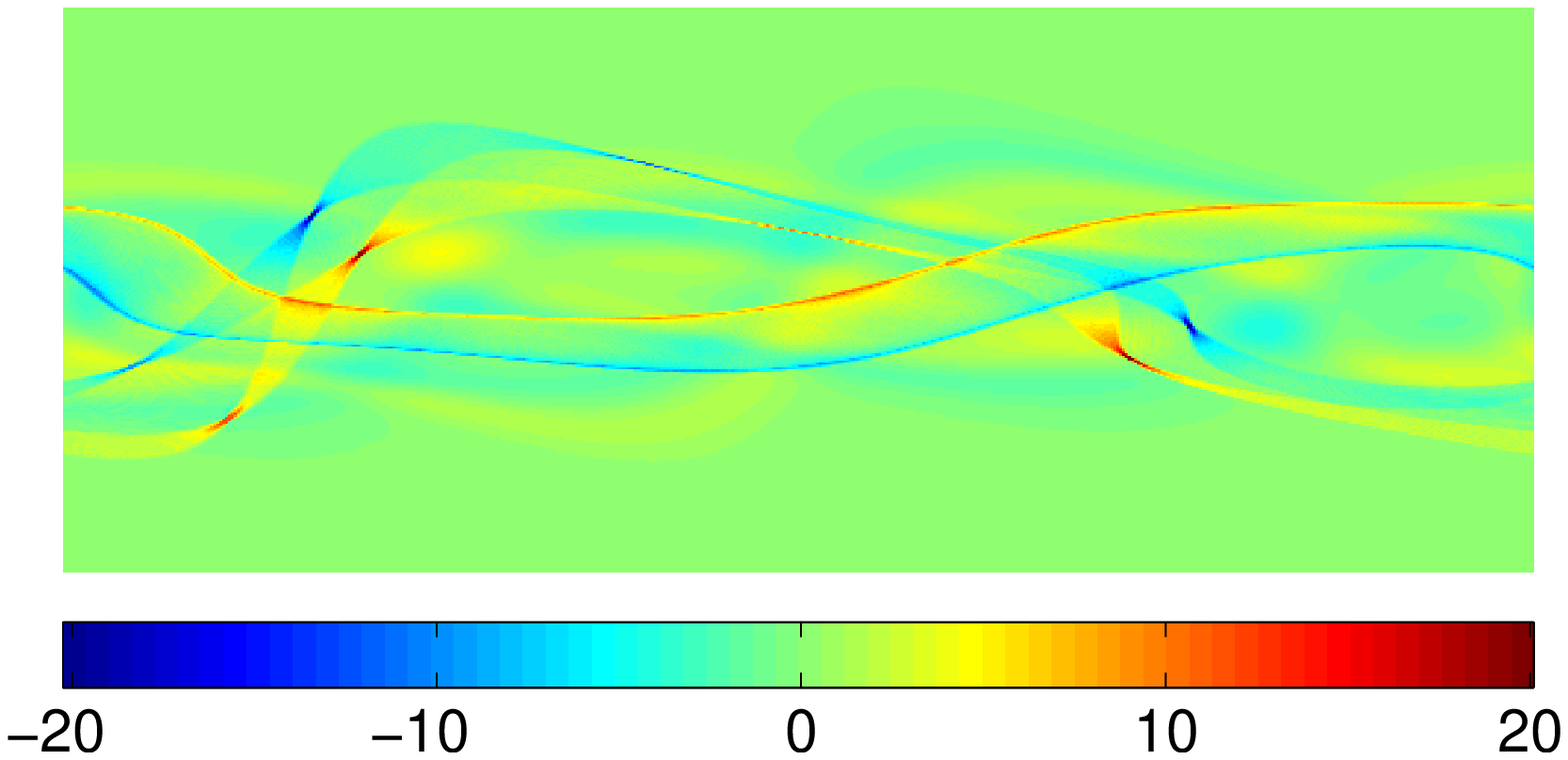}
    \includegraphics[height=0.17\textheight]{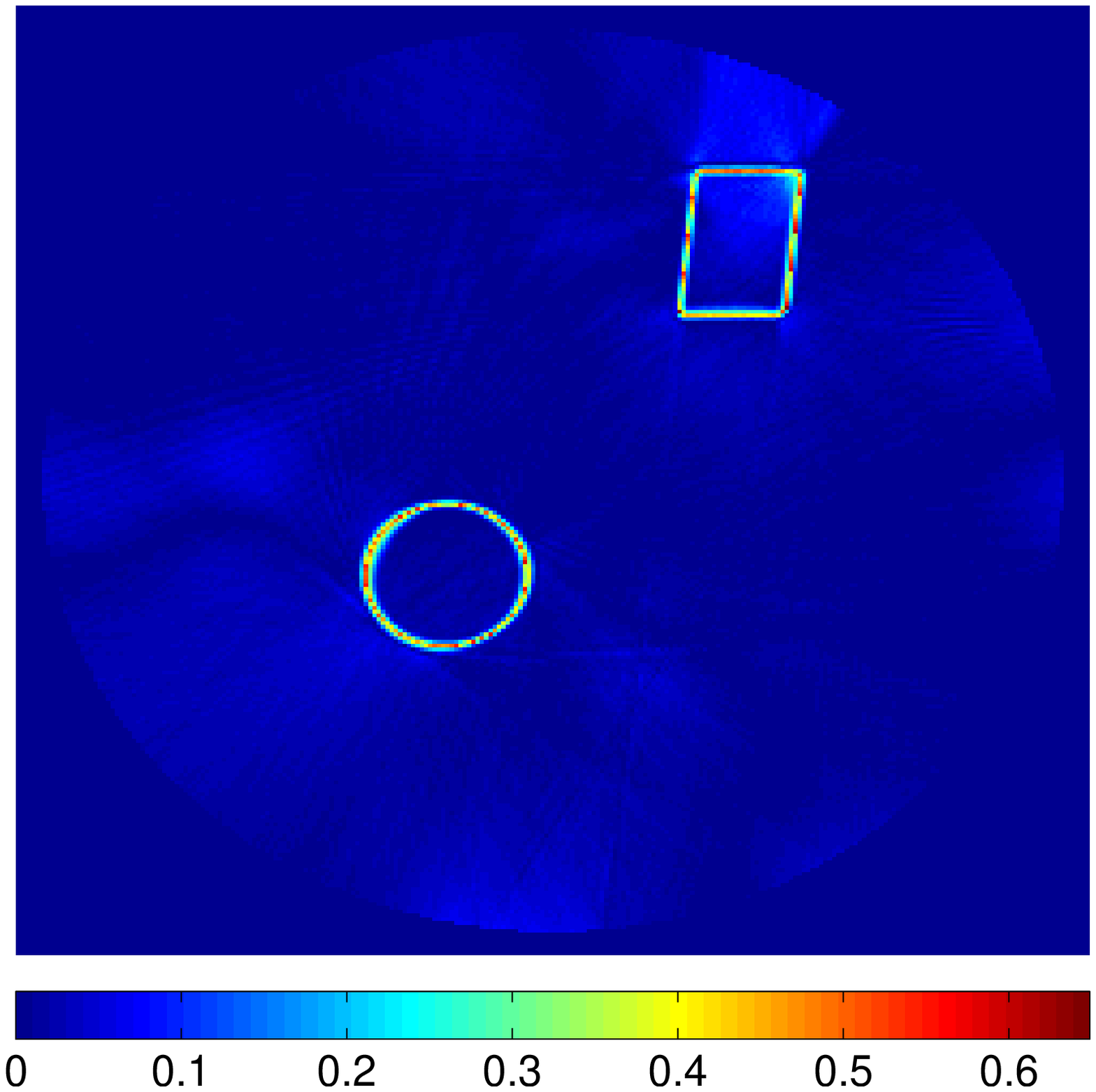}
    \includegraphics[height=0.17\textheight]{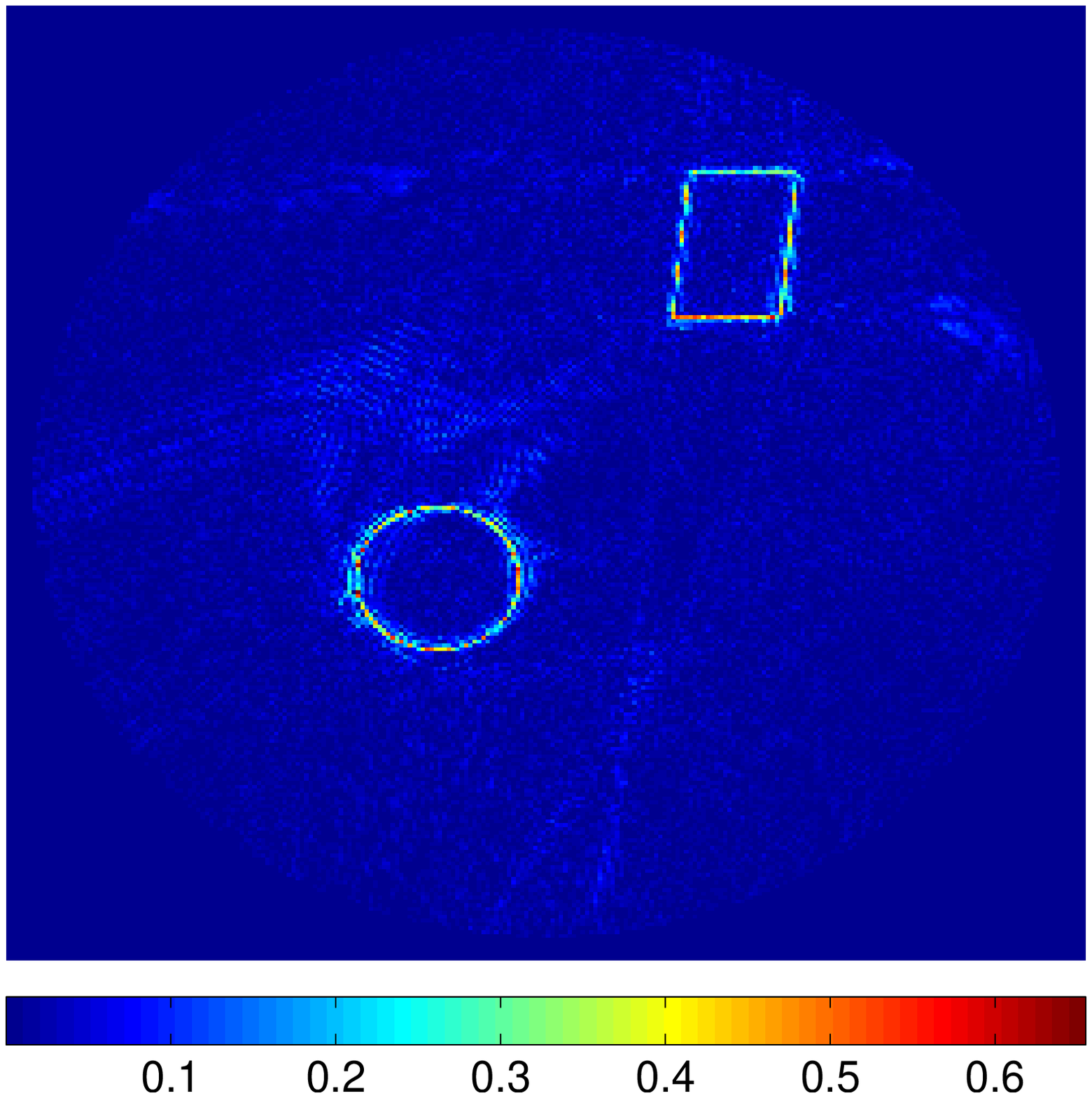} \\
    \includegraphics[height=0.17\textheight]{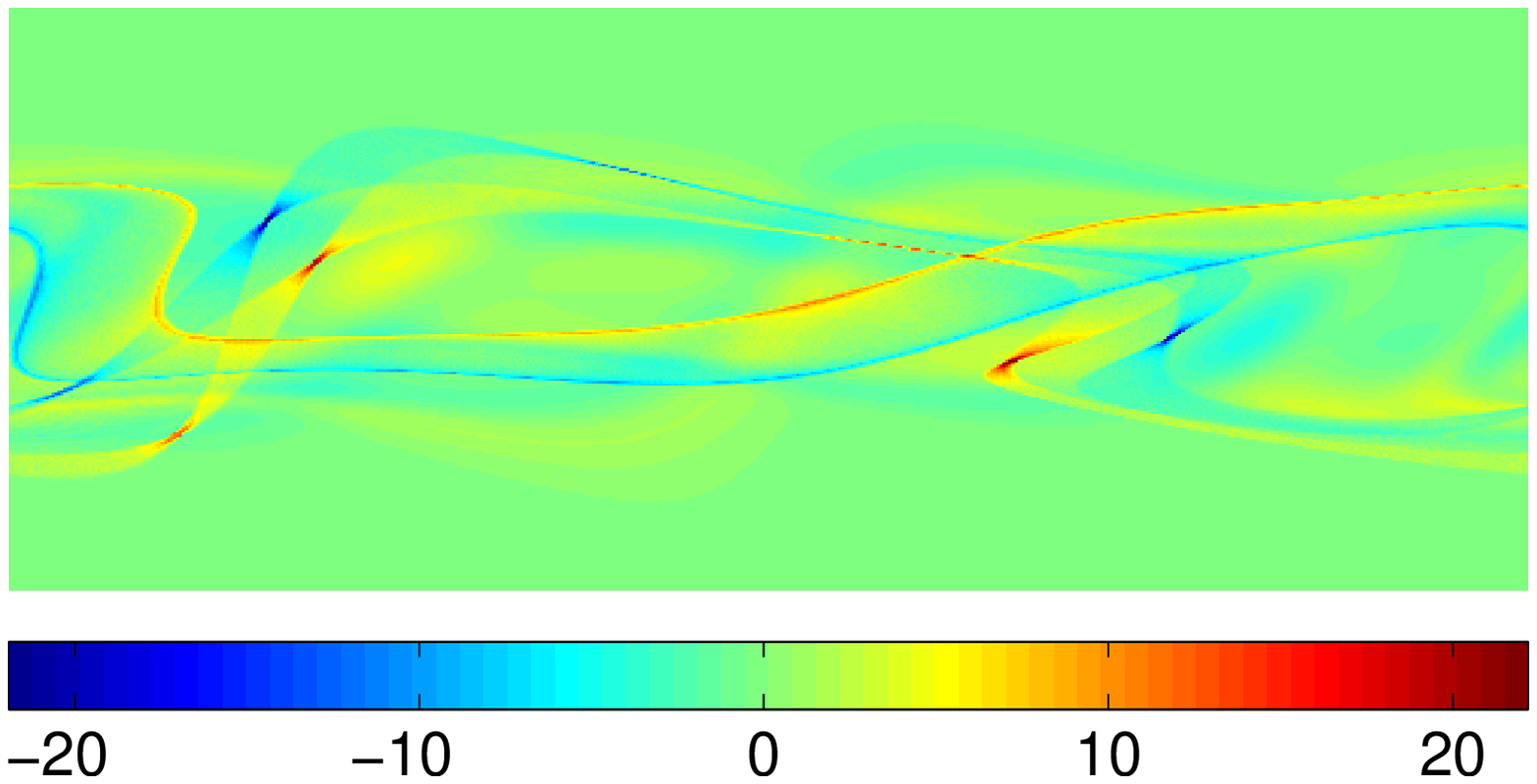}
    \includegraphics[height=0.17\textheight]{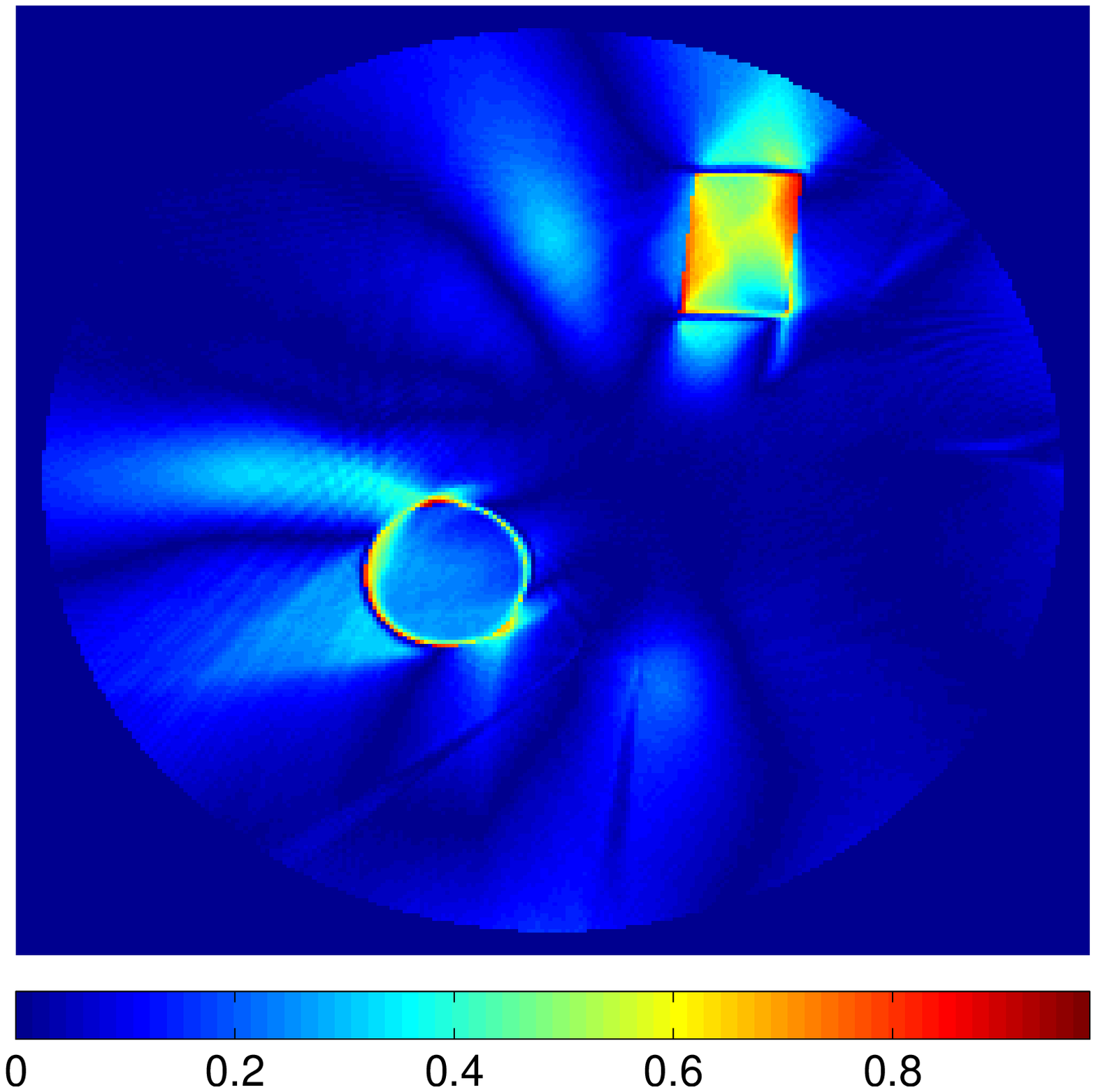}
    \includegraphics[height=0.17\textheight]{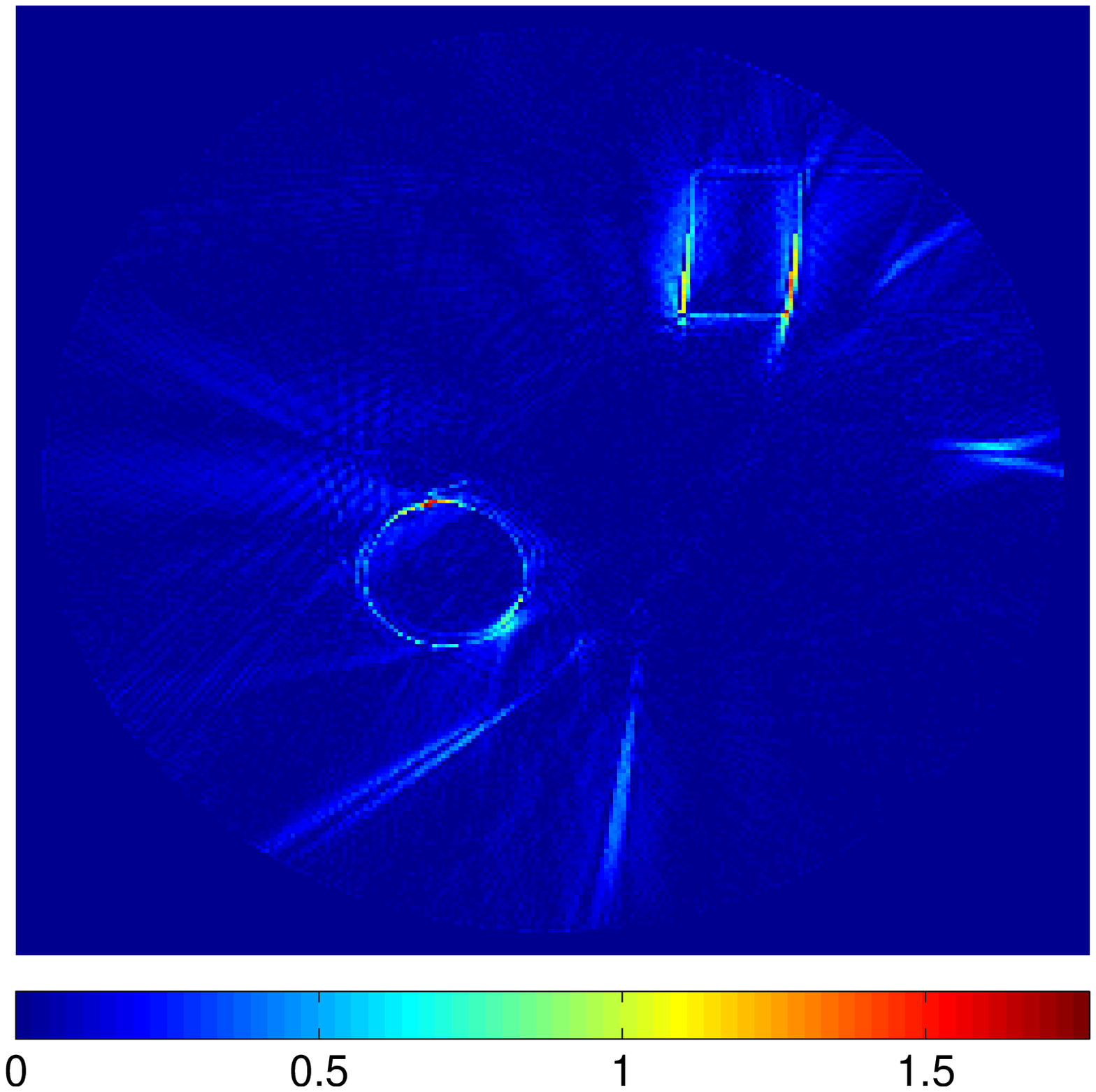}

    \caption{Experiment 3. Top to bottom: $k=0.3$ (simple), $k=0.6$ (not simple) and $k=1.2$ (not simple). Left to right: forward data, pointwise error on $f$ at iteration 1, then 9.}
    \label{fig:lenses_dual_blobs}
\end{figure}

For each experiment, we compute the forward data for all three values $k\in \{0.3,0.6,1.2\}$, and implement an iterative reconstruction following Algo. \ref{algo:neumann} for $9$ iterations. For experiments 1 and 2, forward data are shown side-by-side on Fig. \ref{fig:lensData} and some examples of pointwise errors are displayed on Fig. \ref{fig:pwerror} in a case where the Neumann series converges (Figs. \ref{fig:sample1} and \ref{fig:sample2}) and in a case where it does not (Fig. \ref{fig:sample3}). For Experiment 3, Fig. \ref{fig:lenses_dual_blobs} displays the forward data as well as the pointwise errors at first and last iterations. Notice that, again in the non-simple case, some artifacts appear at the conjugate loci of the conormal singularities. 

\begin{figure}[htpb]
    \centering 
    \subfigure[Exp. 1, error on $f$]{
		\includegraphics[width=0.3\textwidth]{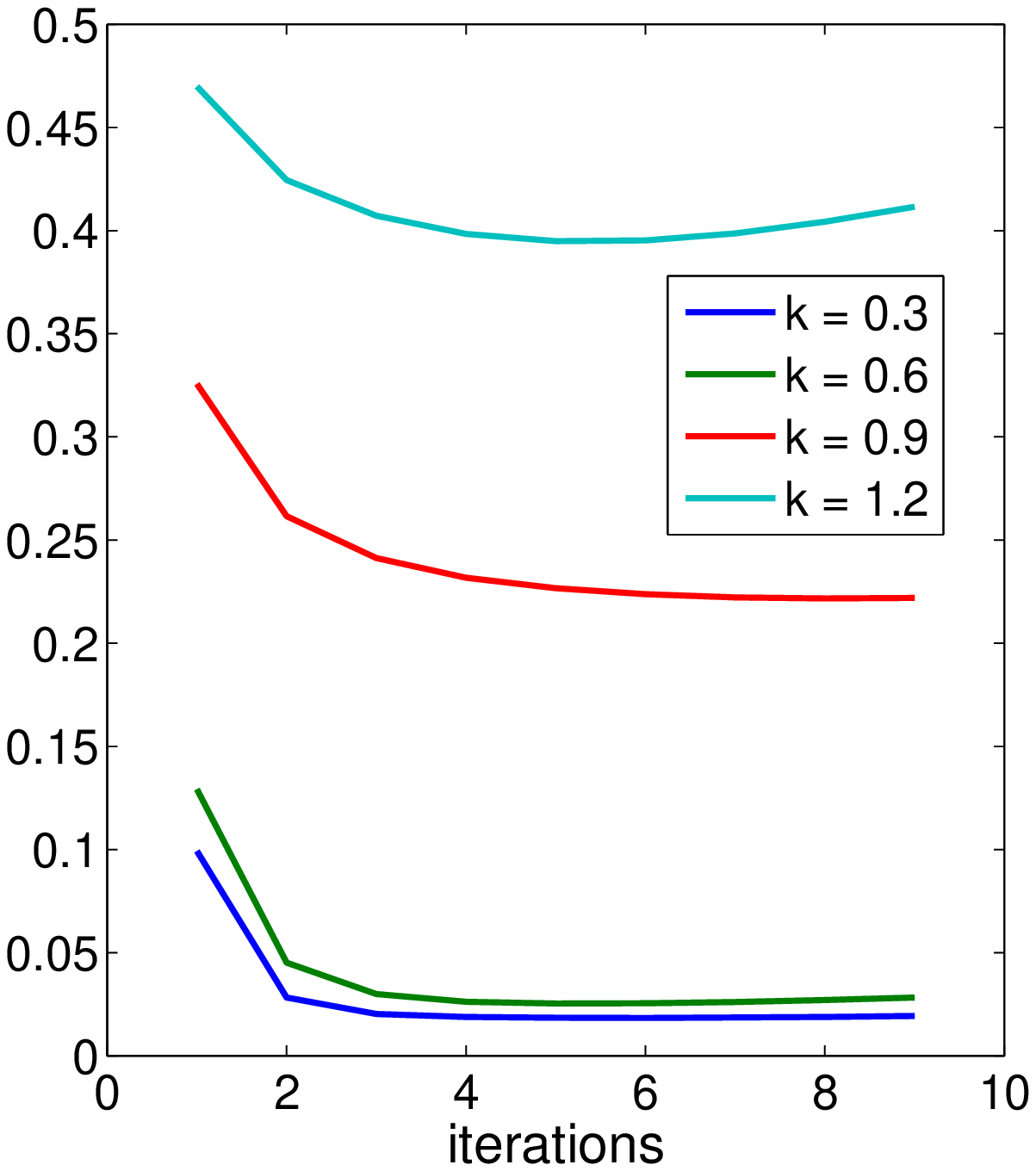}
		\label{fig:conv1}
		}
    \subfigure[Exp. 2, error on $f$]{
	  \includegraphics[width=0.3\textwidth]{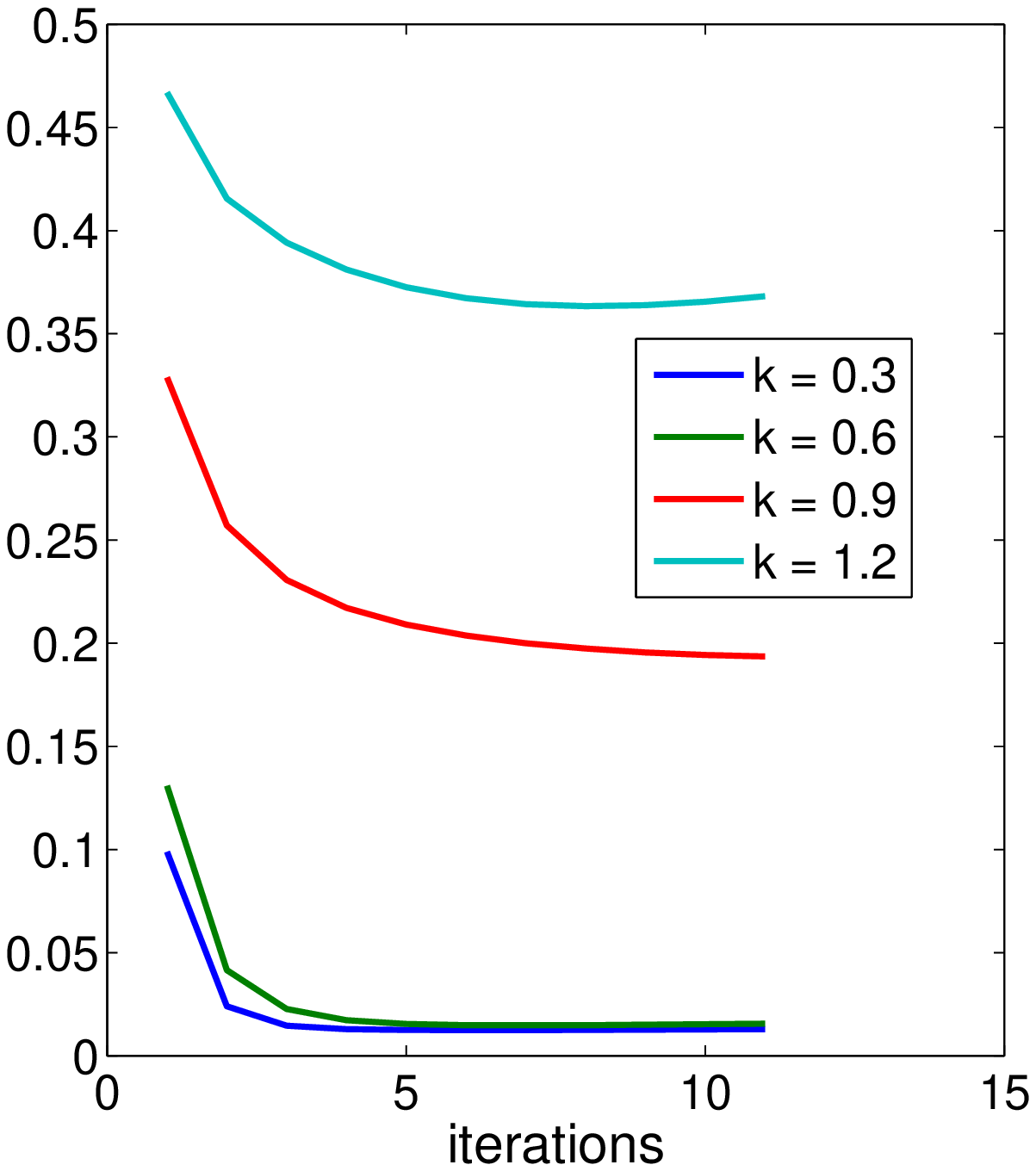}
	  \label{fig:conv2}
    }
    \subfigure[Exp. 3, error on $f$]{
	  \includegraphics[width=0.3\textwidth]{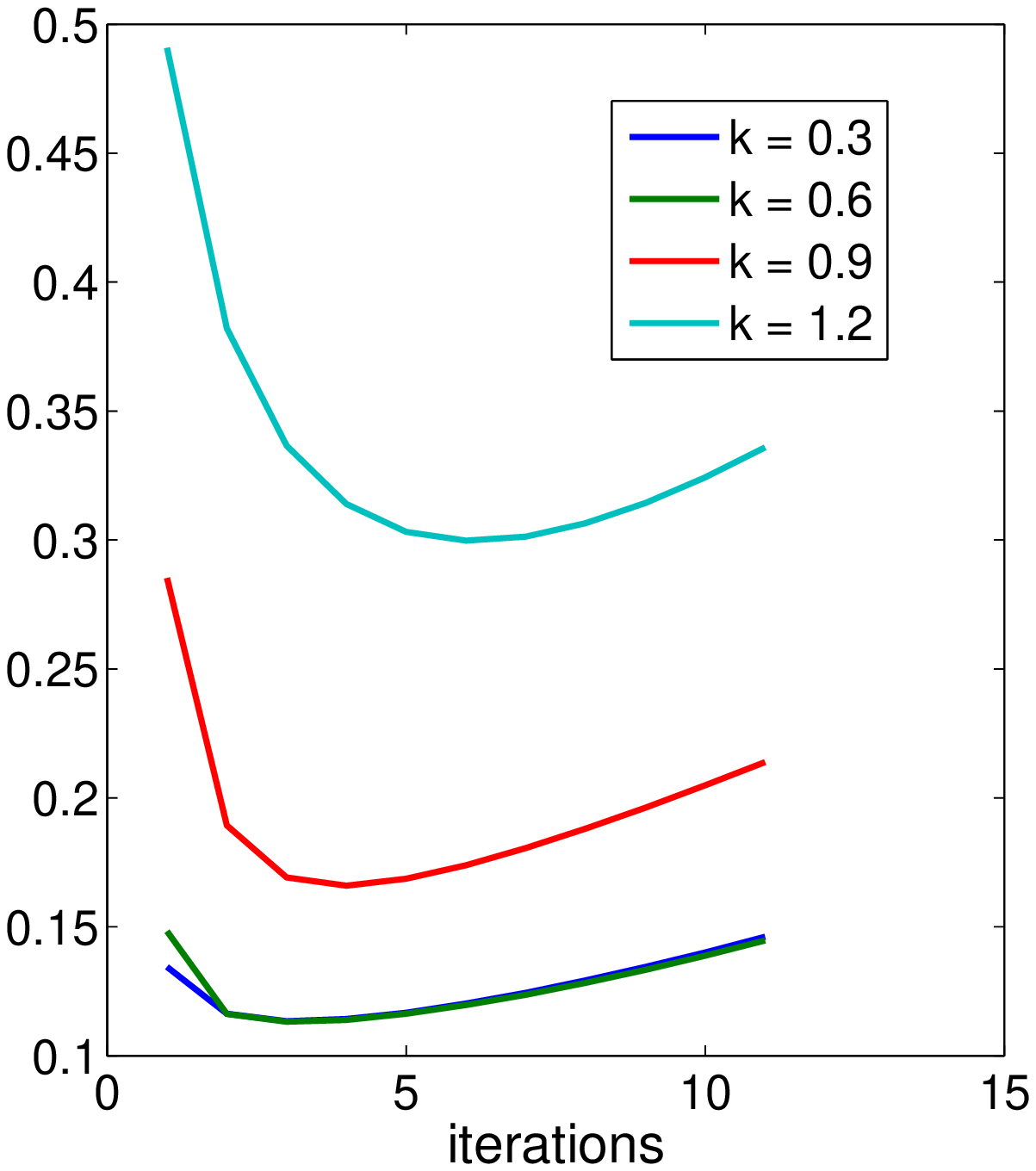}
	  \label{fig:conv3}
    }
    \subfigure[Exp. 1, error on data]{
	  \includegraphics[width=0.3\textwidth]{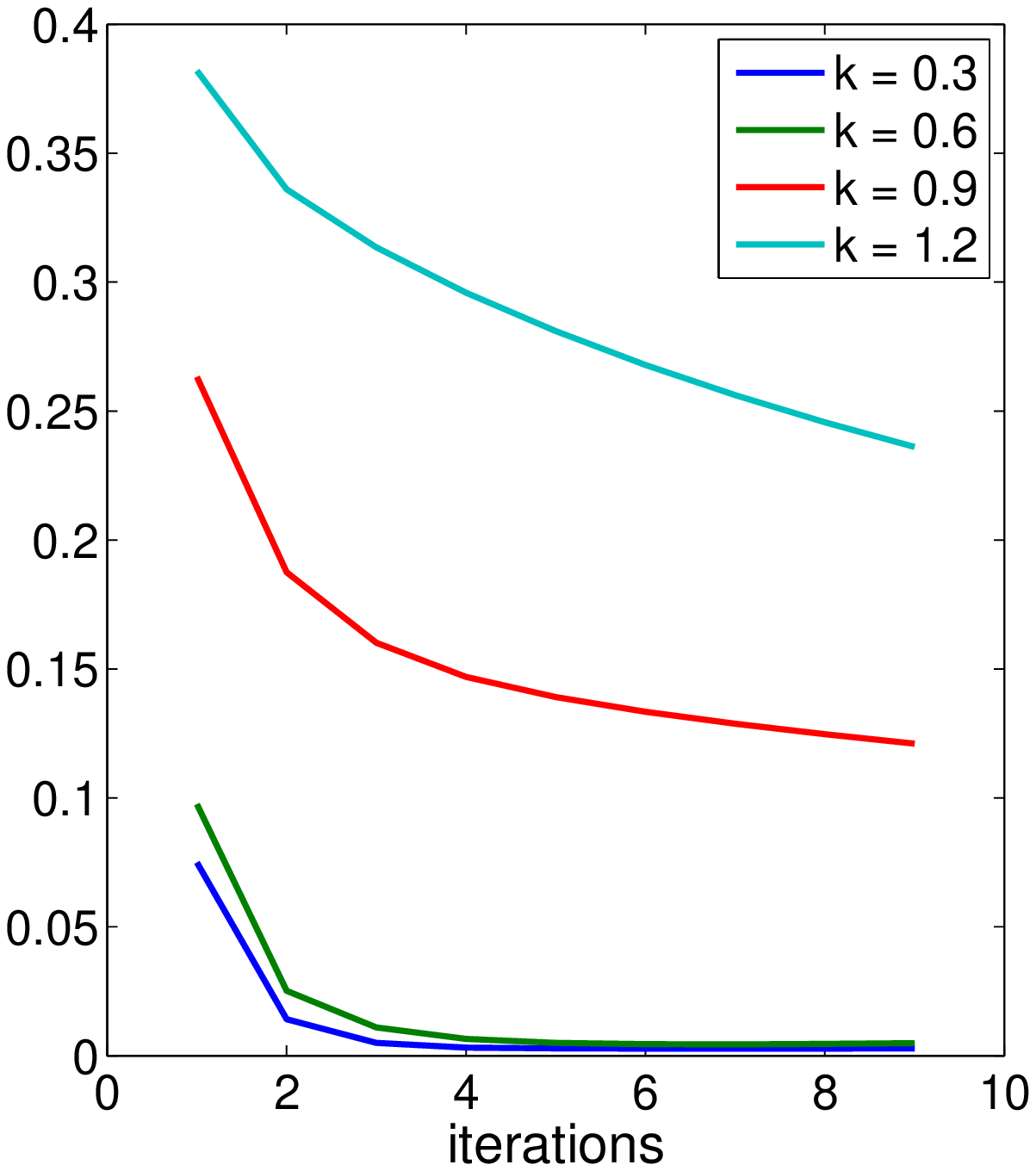}
	  \label{fig:conv4}
    }
    \subfigure[Exp. 2, error on data]{
	  \includegraphics[width=0.3\textwidth]{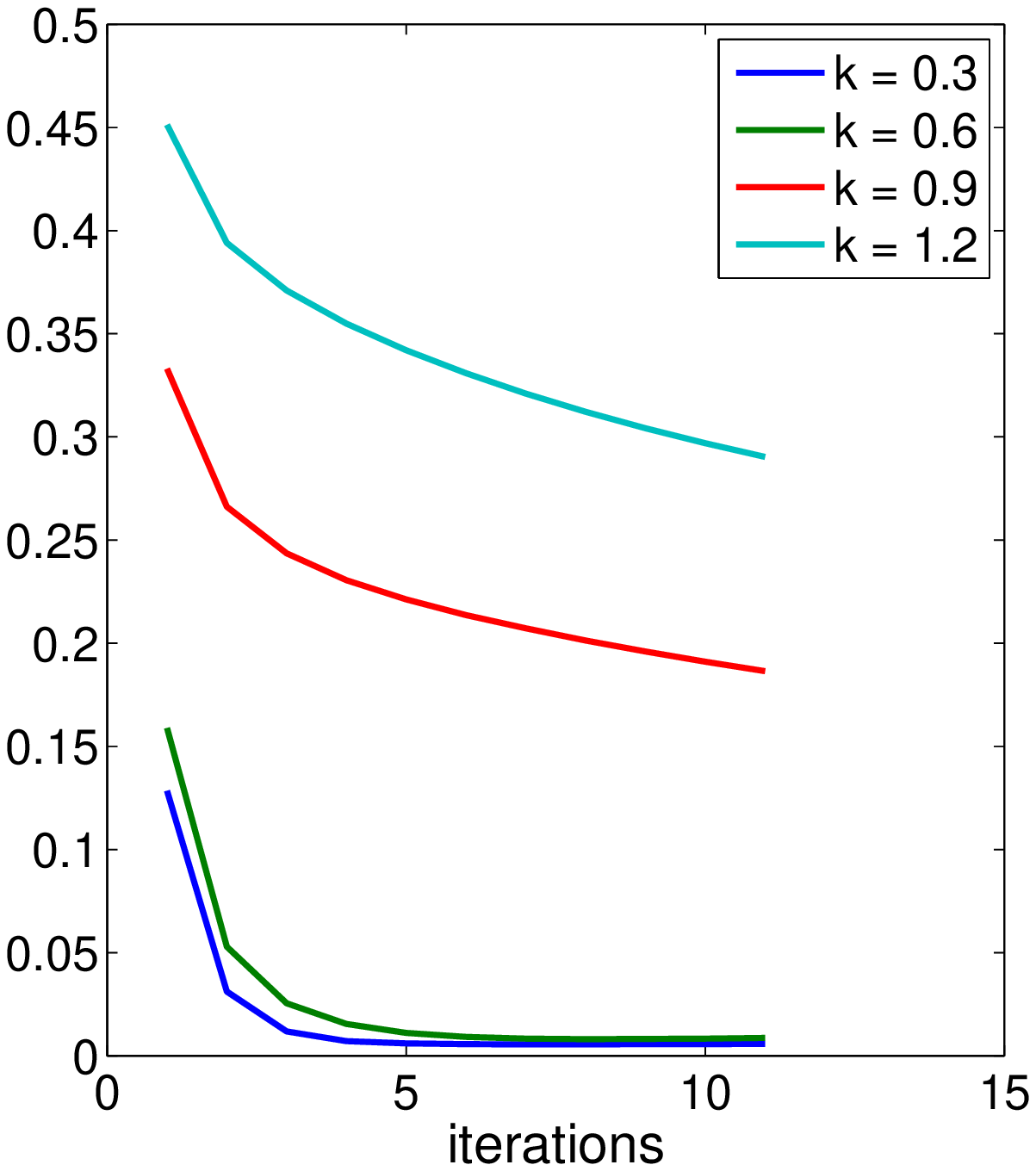}
	  \label{fig:conv5}
    }
    \subfigure[Exp. 3, error on data]{
	  \includegraphics[width=0.3\textwidth]{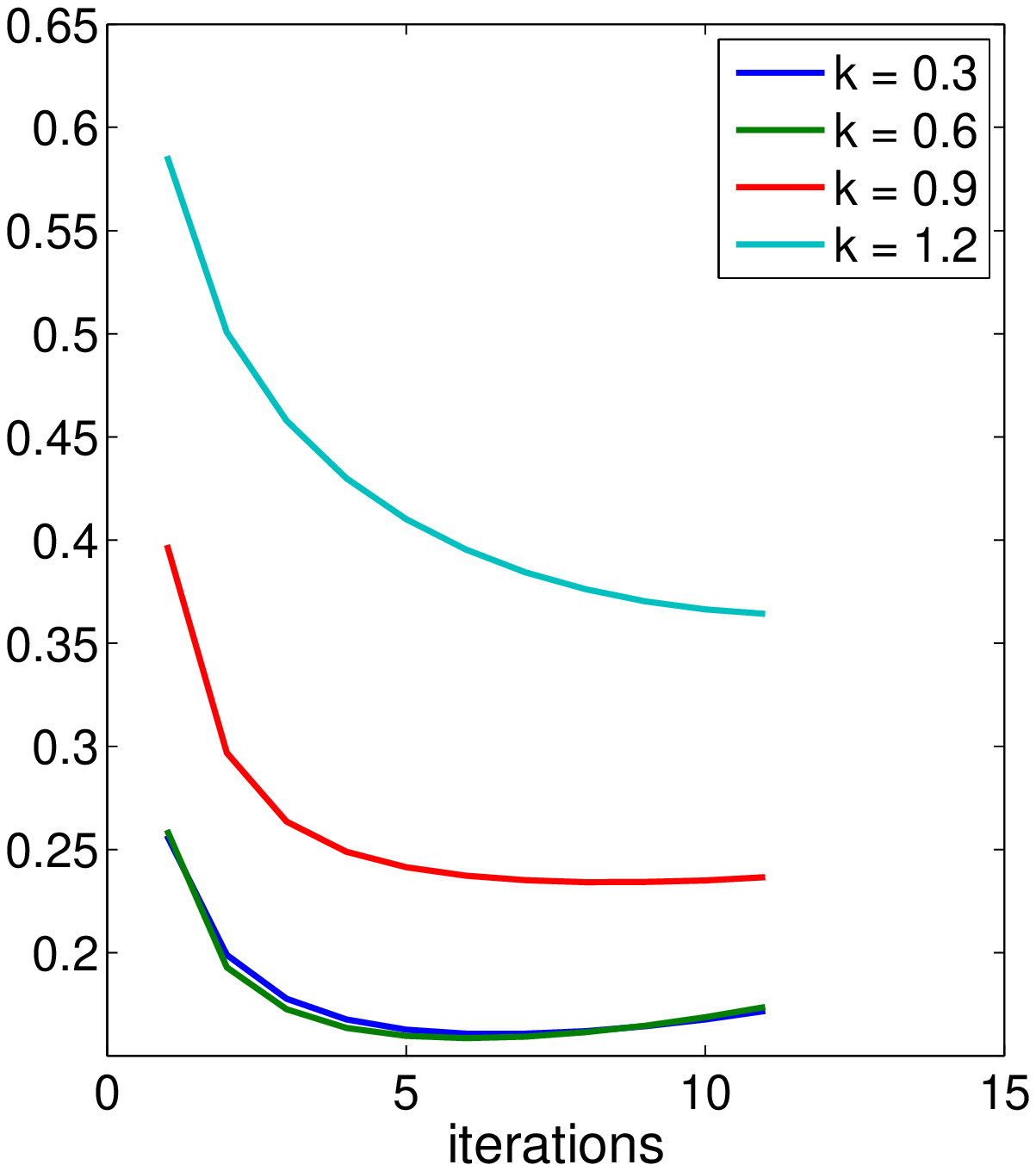}
	  \label{fig:conv6}
    }    
    \caption{Relative $L^2$ error convergence plots on the phantom and on the forward recomputed data ($I_0 f$ for Experiment 1 or $I_1(X_\perp f)$ for Experiments 2,3). In some cases, the reconstructed phantom $f_{rc}$ starts diverging while its ray transform is still converging. }
    \label{fig:convplots}
\end{figure}

Finally, for all three experiments, $L^2$ convergence plots on the unknown and on the data (i.e. we compare the ray transform of the reconstructed quantity with the initial data) are shown on Fig \ref{fig:convplots}. 

\paragraph{Comments}

In the case where $f$ is smooth, we notice rapid convergence (2 iterations) to the exact function in the simple case $k=0.3$ but also in the non-simple case $k=0.6$. In the last case $k=1.2$, some artifacts are noticed on the reconstruction, that do not attenuate as the iterations increase. Some of these artifacts are created near each ``source bump'' as well as their respective conjugate loci. This is an effect that will be discussed at length in a forthcoming work studying non-simple metrics theoretically and numerically, see \cite{Monard2013b}. 

In the case where $f$ has jump discontinuities, these discontinuities can never be exactly resolved no matter how fine the angular resolution is chosen for the backprojection. As a consequence, the error plots contain strong variations at the scale of the grid at the discontinuities, which in turn are amplified by the repeated differentiation that occurs in the iterated reconstruction procedure. This differentiation occurs when computing $X_\perp$, either in the inversion formula when inverting $I_0$, or in the forward operator when inverting $I_1 X_\perp$. Though the iterations improve at first, the repeated differentiations cause the iterations to diverge again, even in the case of simple metrics, as can be seen from the plots in Fig. \ref{fig:conv3}. In the simple cases, this non-convergence effect is presumably not due to the operator $W^2$ having eigenvalues of magnitude larger than $1$ here, as the corresponding eigenvectors, smooth if they existed, would stand out on the error plot.

The non-convergence of the series due to iterated differentiation is due to the fact that both inverse problems considered are ill-posed of order $\frac{1}{2}$ and therefore require regularization (leading to the so-called ``filtered-backprojection'' algorithm in the Euclidean case). In practical settings where measurements may be polluted by noise, regularizing an ill-posed inversion algorithm is also a crucial step to prevent high-frequency noise from overwhelming the reconstructions. While such regularized inversions are well-understood in the Euclidean case (see \cite{Natterer2001}), they rely on the homogeneity of Euclidean space (and the Fourier Slice Theorem that comes with it). It remains, however, an open question to find generalizations of these principles to general Riemannian settings, a question that will be the subject of future work by the author.

\section{Conclusion and remarks} \label{sec:conclu}

We have implemented a {\tt MatLab} code to extend the understanding of geodesic X-ray transforms, in particular their sensitivity to the metric and how injectivity and stability of the associated inverse problem are affected by that metric. The reconstruction algorithms derived in \cite{Pestov2004,Krishnan2010} were successfully implemented as one-shot inversions in the case of manifolds with constant curvature, and as an iterative algorithm when the curvature was close to constant, handling along the way a potentially large class of boundaries. After implementing the Neumann series in cases where it is not theoretically proved that the error operator is a contraction, one finds, aside from numerical instabilities due the to ill-posed nature of the inverse problem (see discussion below), that no smooth structure appears after iterating over the error operator, so that the Neumann series may converge in cases of metrics that are not necessarily close to constant curvature, or even simple, leaving potential room for theoretical improvements. 

The following open questions will be considered in future work. Some of these considerations come as a natural generalization of some issues treated at length in the Euclidean case in \cite{Natterer2001}.

\begin{remunerate}
    \item The accuracy of the transform and its inversion are highly sensitive to the sampling of geodesics at the influx boundary $\partial_+ SM$, and to how the geodesics emanating from this sampling in turn sample the manifold relatively well. It is clear that a uniform sampling of $\Sm^1\times [-\frac{\pi}{2},\frac{\pi}{2}]$ is not optimal in all cases, as negative curvature tends to make geodesics go away from each other, while positive curvature tends to make geodesics concentrate in some areas. As we assume to know the metric here, it is worth investigating how to find an appropriate sampling that mitigates the effect of curvature on the sampling of the manifold. 
    \item It seems more than necessary to generalize the {\it filtered}-backprojection algorithm to the case of non-Eulidean metrics, as the current inversion formulas are ill-posed of order $\frac{1}{2}$. Though the ill-posedness is very mild, the iterated differentiation of the noise will prevent the Neumann series to converge as the numerical errors at small scales will unavoidably be amplified by that differentiation. The main challenge here is to define a concept of regularization that is adapted to the geometry of the manifold and to the X-Ray transform itself. 
  \item Although this has not been observed numerically, if there are cases where $W^2$ or $(W^\star)^2$ is not a contraction yet is compact, the question of inverting $I+W^2$ when either of these operators has eigenvalues of magnitude larger than $1$ is unclear. Methods for doing this should be found.
  \item On the performance side, though the present code makes good use of {\tt Matlab}'s vectorization capabilities, it seems that both forward and inverse transforms are massively parallelizable, since all formulas rely on computing geodesics and base-points of geodesics that are all independent of one another. Accelerating the present code using multi-cores or GPUs will be the object of future work.
\end{remunerate}

\section*{Acknowledgements}
The author would like to thank Gunther Uhlmann and Plamen Stefanov for fruitful discussions, as well as Steve McDowall for helpful comments.  

\bibliographystyle{siam}

\begin{thebibliography}{10}

\bibitem{Anikonov1978}
{\sc Yu.~E. Anikonov}, {\em Some methods for the study of multidimensional
  inverse problems for differential equations}, Nauka Sibirsk. Otdel.,
  Novosibirsk, 1978.

\bibitem{Anikonov1997}
{\sc Yu.~E. Anikonov and V.~Romanov}, {\em On uniqueness of determination of a
  form of first degree by its integrals along geodesics}, J. Inverse Ill-Posed
  Probl., 5 (1997), pp.~487--490.

\bibitem{Bal2010a}
{\sc Guillaume Bal and Fran\c{c}ois Monard}, {\em An accurate solver for
  forward and inverse transport}, Journal of Comp. Phys., 229 (2010).

\bibitem{Dairbekov2006}
{\sc Nurlan~S. Dairbekov}, {\em Integral geometry problem for nontrapping
  manifolds}, Inverse Problems, 22 (2006), pp.~431--445.

\bibitem{Funk1916}
{\sc P.~Funk}, {\em {\"Uber eine geometrishe Anwendung der Abelschen
  Integralgleichnung}}, Math. Ann., 77 (1916), pp.~129--135.

\bibitem{Helgason1999}
{\sc Sigurdur Helgason}, {\em The Radon Transform}, Birk\"auser, second~ed.,
  1999.

\bibitem{Herglotz1905}
{\sc G.~Herglotz}, {\em \"uber die {E}lastizit\"at der {E}rde bei
  {B}er\"ucksichtigung ihrer variablen {D}ichte,}, Zeitschr. f\"ur Math. Phys.,
  52 (1905), pp.~275--299.

\bibitem{Krishnan2010}
{\sc Venky Krishnan}, {\em On the inversion formulas of {P}estov and {U}hlmann
  for the geodesic ray transform}, J. Inv. Ill-Posed Problems, 18 (2010),
  pp.~401--408.

\bibitem{Lee1997}
{\sc J.M. Lee}, {\em Riemannian Manifolds, An Introduction to Curvature.},
  vol.~176 of Graduate Texts in Mathematics, Springer, 1997.

\bibitem{Monard2013a}
{\sc Fran\c{c}ois Monard}, {\em On reconstruction formulas for the ray
  transform acting on symmetric differentials on surfaces}, Inverse 
	Problems, to appear (2014).
	
\bibitem{Monard2013b}
{\sc Fran\c{c}ois Monard, Plamen Stefanov, and Gunther Uhlmann}, {\em The
  geodesic ray transform on riemannian surfaces with conjugate points}, in 
	preparation (2014).

\bibitem{Mukhometov1977}
{\sc R.G. Mukhometov}, {\em The reconstruction problem of a two-dimensional
  riemannian metric, and integral geometry}, Dokl. Akad. Nauk. SSSR, 232
  (1977), pp.~32--35.
\newblock (Russian).

\bibitem{Natterer2001}
{\sc Frank Natterer}, {\em The Mathematics of Computerized Tomography}, SIAM,
  2001.

\bibitem{Paternain2012}
{\sc Gabriel Paternain, Mikko Salo, and Gunther Uhlmann}, {\em The attenuated
  ray transform for connections and higgs fields}, Geom. Funct. Anal. (GAFA),
  22 (2012), pp.~1460--1480.

\bibitem{Paternain2012a}
\leavevmode\vrule height 2pt depth -1.6pt width 23pt, {\em Spectral rigidity
  and invariant distributions on {A}nosov surfaces},  (2012).
\newblock arXiv:1208.4943.

\bibitem{Paternain2013a}
\leavevmode\vrule height 2pt depth -1.6pt width 23pt, {\em On the range of the
  attenuated ray transform for unitary connections},  (2013).
\newblock arXiv:1302.4880.

\bibitem{Paternain2011a}
\leavevmode\vrule height 2pt depth -1.6pt width 23pt, {\em Tensor tomography on
  surfaces}, Inventiones Math., 193 (2013), pp.~229--247.

\bibitem{Pestov2004}
{\sc Leonid Pestov and Gunther Uhlmann}, {\em On the characterization of the
  range and inversion formulas for the geodesic {X}-ray transform},
  International Math. Research Notices, 80 (2004), pp.~4331--4347.

\bibitem{Pestov2005}
{\sc Leonid Pestov and Gunther Uhlmann}, {\em Two-dimensional compact simple
  {R}iemannian manifolds are boundary distance rigid}, Annals of Mathematics,
  161 (2005), pp.~1093--1110.

\bibitem{Qian2011}
{\sc Jianliang Qian, Plamen Stefanov, Gunther Uhlmann, and Hong-Kai Zhao}, {\em
  An efficient {N}eumann-series based algorithm for the thermoacoustic and
  photoacoustic tomography with variable sound speed}, SIAM J. Imaging
  Sciences, 4 (2011), pp.~850--883.

\bibitem{Radon1917}
{\sc Johann Radon}, {\em {\"U}ber die {B}estimmung von {F}unktionen durch ihre
  {I}ntegralwerte l\"angs gewisser {M}annigfaltigkeiten}, Berichte \"uber die
  Verhandlungen der K\"oniglich-S\"achsischen Akademie der Wissenschaften zu
  Leipzig, Mathematisch-Physische Klasse, 69 (1917), pp.~262--277.

\bibitem{Salo2011}
{\sc Mikko Salo and Gunther Uhlmann}, {\em {T}he {A}ttenuated {R}ay {T}ranform
  on {S}imple {S}urfaces}, J. Diff. Geom., 88 (2011), pp.~161--187.

\bibitem{Sharafudtinov1994}
{\sc Vladimir Sharafudtinov}, {\em Integral geometry of tensor fields}, vsp,
  Utrecht, The Netherlands, 1994.

\bibitem{Sharafudtinov1997}
\leavevmode\vrule height 2pt depth -1.6pt width 23pt, {\em Integral geometry of
  tensor fields on a surface of revolution}, Siberian Math. J., 38 (1997).

\bibitem{Spivak1999}
{\sc M.~Spivak}, {\em A comprehensive Introduction to Differential Geometry},
  vol.~4, publish or Perish, 1999.

\bibitem{Stefanov2004}
{\sc Plamen Stefanov and Gunther Uhlmann}, {\em Stability estimates for the
  x-ray transform of tensor fields and boundary rigidity}, Duke Math. J., 3
  (2004), pp.~445--467.

\bibitem{Stefanov2008}
\leavevmode\vrule height 2pt depth -1.6pt width 23pt, {\em Integral geometry of
  tensor fields on a class of non-simple riemannian manifolds}, American J.
  Math, 130 (2008), pp.~239--268.

\bibitem{Stefanov2012a}
\leavevmode\vrule height 2pt depth -1.6pt width 23pt, {\em The geodesic {X}-ray
  transform with fold caustics}, Analysis and PDE, 5 (2012), pp.~219--260.

\bibitem{Uhlmanna}
{\sc Gunther Uhlmann and Andr\'as Vasy}, {\em The inverse problem for the local
  geodesic ray transform}, preprint,  (2012).
\newblock arXiv:1210.2084.

\end{thebibliography}

\end{document}